\documentclass[reqno,final,oneside]{amsart}
\usepackage[textwidth=16.1cm]{geometry}

\usepackage{mathrsfs}  
\usepackage{cancel} 
\usepackage[dvipsnames]{xcolor}
\usepackage{enumitem}
\usepackage{mathtools}
\usepackage[final]{hyperref}
\usepackage{bm}
\usepackage{subcaption}
\usepackage{comment}
\usepackage{xargs}
\usepackage{soul}
\usepackage{subcaption}
\usepackage[textsize=tiny,textwidth=2.05cm]{todonotes}
\newcommandx{\task}[2][1=]{\todo[linecolor=blue,backgroundcolor=blue!30,#1]{#2}}
\newcommandx{\note}[2][1=]{\todo[linecolor=green,backgroundcolor=green!30,#1]{#2}}
\newcommandx{\error}[2][1=]{\todo[linecolor=red,backgroundcolor=red!30,#1]{#2}}
\usepackage{marginnote}
\usepackage[dvipsnames]{xcolor}
\usepackage{mathrsfs}
\usepackage{mathtools}
\usepackage{esint}
\usepackage{stackengine}
\usepackage{hyperref}
\usepackage[marginal]{showlabels}

\usepackage{bbm}

\newtheorem{theorem}{Theorem}[section]
\newtheorem{lemma}[theorem]{Lemma}
\newtheorem{proposition}[theorem]{Proposition}
\newtheorem{corollary}[theorem]{Corollary}

\theoremstyle{definition}
  \newtheorem{definition}[theorem]{Definition}

\theoremstyle{remark}
  \newtheorem{remark}[theorem]{Remark}

\newcommand{\eps}{\varepsilon}
\newcommand{\N}{\mathbb{N}}
\newcommand{\Z}{\mathbb{Z}}
\newcommand{\R}{\mathbb{R}}

\newcommand{\defas}{\coloneqq}

\usepackage{stackengine}

\newcommand{\Me}[1]{M^{#1, \mathbf{v}}_{n,\triangle}}
\newcommand{\Met}[2]{M^{#1,\mathbf{v}}_{n, \triangle}(#2)}
\newcommand{\Cl}[2]{\mathcal{C}^{#1, {\rm L}}_n(#2)}

\newcommand{\Cn}[2]{\mathcal{C}^{#1}_n(#2)}

\newcommand{\Kl}[2]{{K}^{#1}_n(#2)}

\newcommand{\BBB}{\color{black}}

\newcommand{\EEE}{\color{black}}

\newcommand{\MMM}{\color{black}}

\newcommand{\OOO}{\color{black}}

\newcommand{\JJJ}{\color{black}}

\renewcommand{\L}{\mathcal{L}}

\newcommand{\T}{\mathcal{T}}

\newcommand{\wto}{\rightharpoonup}

\DeclarePairedDelimiterX\setof[1]\{\}{#1}
\DeclarePairedDelimiterX\abs[1]\lvert\rvert{#1}
\DeclarePairedDelimiterX\norm[1]\lVert\rVert{#1}
\DeclarePairedDelimiterX\sprod[2]\langle\rangle{#1, #2}

\usepackage{enumitem}
\numberwithin{equation}{section}

\setlist[enumerate,1]{font=\normalfont}
\setlist[itemize,1]{font=\normalfont}
\newlist{thmlist}{enumerate}{1}
\setlist[thmlist]{label=(\roman{thmlisti}),
	ref=(\roman{thmlisti}),font=\normalfont,
	noitemsep}

\title[Atomistic-to-continuum convergence for quasi-static crack growth]{Atomistic-to-continuum convergence for quasi-static crack growth in brittle materials}


\author[M.~Friedrich]{Manuel Friedrich} 
\address[Manuel Friedrich]{%
  Department of Mathematics \\
  Friedrich-Alexander Universit\"at Erlangen-N\"urnberg \\
   Cauerstr.~11, D-91058 Erlangen, Germany 
}
\email{manuel.friedrich@fau.de}

\author[J.~Seutter]{Joscha Seutter}
\address[Joscha Seutter]{
  Department of Mathematics \\
  Friedrich-Alexander Universit\"at Erlangen-N\"urnberg \\
  Cauerstr.~11, D-91058 Erlangen, Germany
}
\email{joscha.seutter@fau.de}

\begin{document}
\begin{abstract}
    We study the atomistic-to-continuum limit for a model of quasi-static crack evolution driven by time-dependent boundary conditions. We consider a two-dimensional atomic mass spring system whose interactions are modeled by classical interaction  potentials, supplemented by a suitable irreversibility condition accounting for the breaking of atomic bonding. In a simultaneous limit of vanishing interatomic distance and discretized time step, we identify  a continuum  model of quasi–static \BBB crack \EEE growth in brittle fracture \BBB \cite{frma98} \EEE featuring  an irreversibility condition, a global stability,  and  an energy balance. The proof of   global stability relies on a careful adaptation of the jump-transfer argument in \cite{Francfort-Larsen:2003} to the atomistic setting.  
\end{abstract}
\maketitle

\section{Introduction}
Since the pioneering work of {\sc  Griffith} \cite{griffith}  in the 1920's, brittle materials and their fracture behavior have been subject to intensive  research in mechanical engineering. The groundbreaking idea of Griffith's theory lies in regarding crack formation and  propagation as the result  of a competition between   elastic bulk energy and a surface energy related to the increase of the crack. In their seminal paper \cite{frma98}, {\sc Francfort and Marigo} adopted this viewpoint and paved the way to a variational approach to fracture, where displacements and crack paths are found by   the principle of energy minimization.  The proposed evolutionary model, called a \emph{quasi-static crack evolution},   is characterized  by three principles, namely an irreversibility condition, a static equilibrium at each time, and an energy balance reflecting the nondissipativity of the process. In contrast to the classical theory, this framework does not rely on prescribed crack paths and also effectively accounts for crack initiation.

Over the last two decades, the variational theory of brittle fracture has been studied in a remarkable depth, see  \cite{BFM}  for a broad  overview. Here, we mention only some of the numerous contributions \cite{Chambolle:2003, dMasoFranToad, Lazzaroni, DalMasoToader, toader,   Francfort-Larsen:2003,  FriedrichSolombrino} addressing the  mathematical well-posedness of the model  from \cite{frma98} in various settings.  In all results, the existence of continuous-time crack evolutions is established as the limit of suitably  time-discretized approximations resulting from an iterative minimization of so-called Griffith-type  functionals. Such  energies  comprise  elastic bulk contributions for the sound regions
of the body and surface terms along the crack, and  are typically defined on (generalized)  functions of bounded variation \cite{Ambrosio-Fusco-Pallara:2000}.  Among various technical difficulties that arise in the  discrete-to-continuous passage in time, the most severe one  is the \emph{stability of unilateral minimality properties} \cite{GiacPonsi} which ensures global minimality of displacements along the irreversible fracture process.

Besides being the limit of time-discretized solutions, quasi-static crack evolutions have been identified as effective variational limits of  sequences of problems in various settings of applicative relevance, including phase-field approximation schemes \cite{Giacomini:2005},  finite element approximation \cite{GP1, GP2},  homogenization of brittle composite materials \cite{GiacPonsi}, or a cohesive-to-brittle passage  in the limit of infinite specimen size  \cite{Giacomini:2005b}. The goal of this work is to advance the understanding of effective theories for  quasi-static crack growth    by establishing a rigorous connection between atomistic and continuum models in brittle fracture.

The passage from atomistic  systems to continuum models in the limit of vanishing interatomic distance has been a thriving field of research in the last years, see \cite{blancAtomisticContinuumLimits2007, braidesDiscretetoContinuumVariationalMethods} for an overview. Such an asymptotic analysis is frequently grounded on  the variational tool of $\Gamma$-convergence  \cite {braidesGammaconvergenceBeginners2002, dalmasoIntroductionGConvergence1993}, which ensures convergence of minimizers and hence allows for establishing a rigorous relation between microscopic and macroscopic models. Among the vast body of literature, this approach has also been used to derive and validate continuum theories in brittle and cohesive fracture\footnote{Already Griffith \cite[p.~165]{griffith}   motivated his approach by heuristic considerations on the microscopic nature of fracture: `It is known that, in the formation of a crack in a body composed of molecules which attract one another, work must be done against the cohesive forces of the molecules on either side of the crack. This work appears as potential surface energy.'}, both in one  \cite{braidesVariationalFormulationSoftening1999,  Braides-Lew-Ortiz:06,  schaffnerLennardJonesSystemsFinite2017} and higher dimensions    \cite{ Bach.Braides, Braides-Gelli:2002-2,  FS151,  rufApproxMumfordShah}. 
 Yet, the analysis of static problems has been predominant and the evolutionary nature of crack growth has largely been neglected. Indeed, for atomistic problems in fracture mechanics, available rigorous results seem to be limited to \cite{braidesVariationalEvolutionOnedimensional2014} where  a one-dimensional Lennard-Jones \BBB system \EEE  is investigated via \BBB a \EEE  minimizing movements scheme. The subject of this paper is to perform an atomistic-to-continuum analysis beyond the 1D  setting accounting for irreversibility,  and to derive a continuum  quasi-static crack evolution in the sense of \cite{frma98}.

%, which in this setting ensures convergence of minimizers of the atomistic system
%to the ones of the continuum problem.

We now describe our setting. We build upon the previous work in \cite{BFS} where an atomistic model for crack growth has been introduced including an irreversibility condition as a key feature. More precisely, we consider systems of interacting particles arranged in an atomic lattice occupying the bounded reference domain of a material. The interactions are modeled by classical potentials from Molecular Mechanics  \cite{Molecular, Lewars},  e.g., Lennard-Jones-type potentials.   Following an energy minimization principle, the evolution is driven by time-dependent boundary conditions. We implement an irreversibility constraint along the evolution which resembles the use of a maximal-opening memory variable in cohesive continuum models for crack evolution  \cite{Zanini, Giacomini:2005b}:   Whenever the distance between two atoms
does not exceed a certain threshold, interactions are considered to be `elastic' and can recover
completely after removal of loading. On the contrary, stretched beyond the threshold, interactions are supposed to be `damaged’, and
the system  should be affected at all future times. This is modeled by a  memory variable that tracks the deformation history of each interaction in the system by taking the supremum over all past time steps.  

While the modeling framework in \cite{BFS} is rather general,  in this work  we adopt a  specific  simple model for performing the passage to the continuum theory. More precisely, we assume that the atoms in the reference configuration are arranged in a triangular lattice and  interact via nearest-neighbor potentials. Moreover, we restrict ourselves to the case of anti-plane shear displacements. As the interatomic distance $\eps$ tends to zero,   the atomistic energies can then be related by means of $\Gamma$-convergence to a \emph{free-discontinuity  functional} of the form
\begin{align}\label{first en}
\int_{\Omega}\Phi(\nabla u)\,{\rm d}x+ \int_{S(u)}\varphi(\nu_u)\, {\rm d} \mathcal{H}^{1}\,,
\end{align}
where $\Omega \subset \R^2$ denotes the reference configuration, $u$ the displacement with gradient $\nabla u$, and $S(u)$ stands for the crack surface with normal $\nu_u$. This follows by adapting the result in \cite{FS151} where a vectorial variant of the present model has been investigated. In particular, the model features  an elastic part with energy density $\Phi$ derived from the \JJJ molecular interaction potentials \EEE and a surface term comparable to the one-dimensional Hausdorff measure of the crack set, i.e., $\mathcal{H}^{1}(S(u))$. A fundamental feature of the model is the anistropic density $\varphi$, which is induced by the triangular lattice and favors the formation of cracks along crystallographic lines, see   \cite{FS152, FS153} for rigorous results on crystal cleavage and   \cite{BFS}  for numerical evidence.

Although we believe that in principle  our analysis could  be carried out for far more general lattice models, including long-range and multi-body interactions (see e.g.\ \cite{Bach.Braides}), we prefer to consider here a simple mass-spring model, which already gives rise to a Griffith-type energy in the continuum limit with a non-degenerate elastic density and   anisotropic surface terms. This allows us to focus on the challenges due to the evolutionary nature of the fracture process, in particular issues related to irreversibility. 

Starting from an incremental minimization scheme with time-step size $\delta$ and corresponding time-discrete solutions for the atomistic model, we show in our main result (Theorem \ref{maintheorem}) that, in the simultaneous limit $\eps,\delta \to 0$, the atomistic evolutions converge to  an  irreversible  quasi-static crack evolution in the sense of \cite{frma98}. More precisely, we find  a pair $t\mapsto (u(t),K(t))$, where $u(t)$ lies in the space $SBV(\Omega)$, see \cite{Ambrosio-Fusco-Pallara:2000}, and $K(t)$ is a rectifiable set with $S(u(t)) \subset K(t)$, satisfying:  (1) an irreversibility condition,   (2) a global stability at all times,  and (3) an energy balance law, see \eqref{finalstability}--\eqref{energybalance} for details. As a byproduct, we show that  the energies converge at all times as $\eps \to 0$. Moreover, for suitable interpolations of the atomistic displacements that are exhibiting jumps, we obtain strong convergence of displacement gradients and convergence of crack sets in the sense of $\sigma$-convergence introduced in \cite{GiacPonsi}, see also Section \ref{set conv-sec}. 

 Brittle materials frequently  develop cracks already at moderately large strains without significant plastic deformation. As in \cite{Braides-Lew-Ortiz:06, FS151, FS152}, this observation is reflected in our model since either  elastic displacements are small, namely of order \JJJ $\eps^{1/2}$ (the exponent $1/2$ \EEE is related to the fact that $\Phi$ is positively  homogeneous of order \JJJ $2$, \EEE see \eqref{eq: lineraition} below), or springs leading to fracture in the limit are elongated by a factor \JJJ $\eps^{-1/2}$. \EEE As a consequence, the limit $\eps \to 0$ involves a linearization of the elastic energy and the fracture response becomes asymptotically  brittle despite the presence of cohesive forces on the atomistic level. In other words, while the surface energy in \eqref{first en} depends on the geometry of the crack, it is not affected by its opening. To our best knowledge, our main theorem represents the first result identifying a quasi-static crack evolution not only  in the atomistic-to-continuum limit but also in the limit from nonlinear-to-linearized elastic energies.  

Whereas the derivation of \eqref{first en} via $\Gamma$-convergence follows along the lines of \cite{FS151}, and the techniques to pass to the limit in energy balances are well established  \cite{dMasoFranToad, GiacPonsi}, the  essential novelty of the proof consists in showing the stability of the unilateral minimality property of the atomistic  evolutions. This is achieved by a suitable adaptation of the  jump-transfer argument from \cite{Francfort-Larsen:2003} to the atomistic setting. A discrete variant of this construction has indeed already been proposed in the literature in terms of discontinuous finite elements \cite{GP1, GP2}. However, the main challenge in the present atomistic setting lies in the fact that `broken' interactions do not contribute a fixed energy, but their contribution depends also on the elongation. (Only in the continuum limit this effect vanishes, as discussed above.) In this sense, our problem is related to the work by {\sc Giacomini} \cite{Giacomini:2005b}, where stability of unilateral minimality is shown in  a cohesive-to-brittle limit. In both problems, it is essential to  transfer the jump  {only} on the part of the crack where the memory variable exceeds a certain constant. We achieve this by splitting the broken interactions on the atomistic level into two parts $\mathcal{C}^{\rm small}$ and $\mathcal{C}^{\rm large}$, related to small and large values of the memory variable, respectively. Then, we carefully introduce interpolations exhibiting jumps only at sites associated to $\mathcal{C}^{\rm large}$.   With this, we can suitably adapt the original    jump-transfer argument and we identify  the continuum crack set by resorting  to the notion of \JJJ $\sigma^{p}$-convergence \EEE introduced in \cite{dMasoFranToad}. Eventually, although in general the concepts of \JJJ $\sigma^p$- \EEE and $\sigma$-convergence are different, we show that they coincide along the evolution, cf.\ also \cite{GiacPonsi}.     

Inspired by  \cite{FS153}, the arguments involve several types of interpolations,  ranging from piecewise affine interpolations to functions exhibiting jumps inside the triangles of the lattice. The use of the interpolations is vital to our reasoning since it allows us  to directly employ tools for $SBV$ functions, such as compactness and lower semicontinuity results as well as  the notion of \JJJ $\sigma^{ p } $-convergence\EEE. This is the main reason why we focus on a simple mass-spring model. In fact, discrete-to-continuum limits  for more general atomistic energies involving long-range interactions are usually based on slicing \cite{Braides-Gelli:2002-2} or on more abstract techniques, such as the localization method of $\Gamma$-convergence and integral representation  \cite{Bach.Braides}. Transferring our arguments  \BBB for \EEE the proof of stability of unilateral minimality to these frameworks appears to require nontrivial adaptations.

 \BBB
Let us mention that in our result  we start from time-discrete solutions for the atomistic evolutions, for existence of time-continuous evolutions  on the atomistic level is not guaranteed by \cite{BFS}. In fact, in \cite{BFS}  time-continuous solutions have only been established for a model with an additional rate-dependent  dissipation term.  \EEE
   
The paper is organized as follows.  In Section  \ref{sec: 2} we introduce the atomistic model for irreversible crack growth from \cite{BFS} and state our main result. Then, in Section \ref{sec: prel} we suitably split the atomistic energy in an elastic and a crack part, and we introduce   useful interpolations for atomistic displacements. Moreover, we prove compactness and lower semicontinuity results  for the elastic energy and deal with different notions of set-convergence. Section \ref{se:c main} is devoted to the proof of our main result,  where we first derive properties of  the  time-discrete atomistic problems and pass to the continuum limit afterwards. To ensure that the latter can be identified as   a quasi-static crack evolution, we crucially  rely on the above-mentioned stability of unilateral minimality whose proof is postponed to Section \ref{stabilitysection}.

\section{Model and main results}\label{sec: 2}

In this section we describe our model and present the main result on an atomistic-to continuum convergence for quasi-static crack growth.

\subsection*{Atomistic anti-plane shear model}
We consider a model for  anti-plane shear   displacements.  This corresponds to a reference
configuration of the form    $U \times \R$,   where $U \subset \R^2$ is a Lipschitz set, and the deformations for $(x_1,x_2) \in U$ and $x_3 \in \R$ are of the special form $\BBB (x_1,x_2,x_3) \mapsto \EEE (x_1,x_2,x_3 + w(x_1,x_2))$, where $w \colon U \to \R$ is a scalar-valued map. Since the behavior is identical for each cross section $\lbrace x_3 = c\rbrace$ for $c\in\R$, one can effectively reduce to deformations $y \colon U \to \R^3$ of the form $y(x_1,x_2) = (x_1,x_2, w(x_1,x_2))$, which corresponds to the cross section  $\lbrace x_3 = 0\rbrace$. We consider an atomistic model where in the reference configuration the particles in the cross section $\lbrace x_3 = 0\rbrace$ are given by the points of the scaled triangular lattice $\eps \L$ inside $U$. Here, 
  \[\mathcal{L}:= \begin{pmatrix}
    1 & 1/2 \\
    0 & \sqrt{3}/2
\end{pmatrix}\mathbb{Z}^2=\{\lambda_1 \mathbf{v}_1+ \lambda_2 \mathbf{v}_2 \colon \lambda_1,\lambda_2 \in \Z\}\,,\]  where  $\mathbf{v}_1:= e_1$ and $\mathbf{v}_2:= (1/2,\sqrt{3}/2)^{T}$ are the lattice vectors, and $\eps>0$ is a small parameter representing  the length scale of the typical interatomic distances. We collect all lattice vectors in the set $\mathcal{V}:=\{\mathbf{v}_1,\mathbf{v}_2,\mathbf{v}_2-\mathbf{v}_1\}$ and denote by $\mathcal{L}_\eps(U) :=  \eps\L \cap U$ the particles in the reference configuration $U$. Now, deformations  $y \colon \mathcal{L}_\eps(U) \to \R^3$  are given by $y(x) = (x,w(x))$ for all $x \in  \mathcal{L}_\eps(U)$, where $w \colon \mathcal{L}_\eps(U) \to \R$ is a scalar displacement field in vertical direction.  We model interactions between the particles by classical potentials from Molecular Mechanics \cite{Molecular, Lewars}. Following   \cite{FS151}, we consider an    associated \emph{phenomenological energy} featuring nearest-neighbor interactions, given by  
\begin{equation}\label{basicenergy}
\mathscr{F}_{\varepsilon}(y):=\frac{1}{2}\, \sum_{(x,x') \in \mathcal{\rm{NN}_{\varepsilon}}(U)} W\Bigl(\frac{|y(x)-y(x')|}{\varepsilon}\Bigr) \,,\end{equation}
where $\mathcal{\rm{NN}_{\varepsilon}}(U)$ denotes the\emph{ nearest neighbors } 
\[\mathcal{\rm{NN}_{\varepsilon}}(U):=\{ (x,x')\in \mathcal{L}_{\varepsilon}(U)\times\mathcal{L}_{\varepsilon}(U)\colon \, |x-x'|=\varepsilon \}\,  \]
(the factor $\frac{1}{2}$ accounts for double counting in the sum), \JJJ and $W$ represents a \emph{two-body interaction potential}  satisfying  
\begin{enumerate}[label=(\roman*)]
    \item $W\geq 0$ and $W(r)=0$ if and only if $r=1$. 
     \item   $W$ is  increasing   on $(1,\infty)$.

       \item   $W$ is Lipschitz on $[1,\infty)$   and $C^2$ on $(1,\infty)$ with $\lim_{r\searrow 1} W'(r)=\mu>0$ and $\lim_{r\searrow 1} |W''(r)| < \infty$.

    \item{\label{define-potential1}} $\lim_{r\to \infty}W(r)=\kappa $ with $\kappa>0$ constant.        
\end{enumerate}
Here,  on  $[1,\infty)$ the potential $W$ mimics a  Lennard-Jones-type potential up to a shift by $\kappa$, which is for convenience only to ensure that the potential is nonnegative. We emphasize that we use a linear instead of quadratic behavior close to the global minimizer $1$ as this allows to derive  an  effective  elastic energy density  in \eqref{eq: lineraition} below with \emph{quadratic} growth. \EEE    We  point out   that discrete energies of the form \eqref{basicenergy} and their relation to continuum energies of Griffith-type in the sense of $\Gamma$-convergence have been studied intensively in the literature in a general framework  including  multi-body  and multi-scale interactions, see e.g.\ \cite{Bach.Braides,  Braides-Gelli:2002-2, Braides-Lew-Ortiz:06}. As our focus lies on the irreversibility of crack growth, we consider here only a simple mass-spring model. 

\EEE

We impose Dirichlet boundary conditions  on a part $\partial_D U \subset \partial U$ of the boundary. Both in the framework of atomistic systems and of function spaces allowing for discontinuities, this is usually implemented by imposing boundary conditions in a \emph{neighborhood} of the boundary. More precisely, we suppose that there exists another Lipschitz set ${\Omega} \supset U$ with  $\partial_D U = \partial U \cap {\Omega}$. Then, for a given function $h \colon  {\Omega} \to  \R^3$, we  call deformations $y \colon \mathcal{L}_\eps( {\Omega}) \to \R^3$ \emph{admissible} with respect to $h$ if they satisfy
\begin{align}\label{y admi}
y(x) = h(x) \quad \text{ for all $x \in \L_\eps( {\Omega} \setminus U)$,     }
\end{align}
where similarly as before, we set $\L_\eps( {\Omega} \setminus U) = \eps\L \cap ( {\Omega} \setminus U)$ and $\L_\eps( {\Omega} )   = \eps\L \cap  {\Omega}$. From now on, deformations will always be defined on the enlarged set $ {\Omega}$, and accordingly the sum in \eqref{basicenergy} will run over pairs of  nearest-neighbor   points in $\mathcal{L}_\eps( {\Omega})$, denoted by ${\rm NN}_\eps( {\Omega})$. 

%\BBB Instead of (or additional to) Dirichlet boundary conditions, one could also consider body or boundary forces. However, as traction  tests with body forces are ill-posed for Griffith functionals in a continuum variational formulation, we prefer to focus on boundary value problems also on the atomistic level. For a  thorough discussion on the comparison of soft and hard devices in the context of the variational formulation of Griffith's fracture, we refer the reader to \cite[Section 3]{BFM}. \EEE

\subsection*{Evolutionary model: Irreversibility condition and memory variable}\label{sec: evo model}

Based on the energy \eqref{basicenergy}, we now   introduce an evolutionary model accounting for an irreversibility constraint along the fracture process in a time interval $[0,T]$. In analogy to \cite{BFS}, we implement such an irreversibility condition by incorporating a memory variable for each interacting pair of neighboring atoms $(x,x') \in {\rm NN}_\eps(\Omega)$. More precisely, suppose that for a fixed time step $t^k\in [0,T]$ a \BBB finite  \EEE family of deformations $(y^j)_{j<k}$ at previous time steps $(t^j)_{j<k}$ is given. \BBB Then we \EEE consider the  \emph{memory variable} related to $(x,x') \in {\rm NN}_\eps(\Omega)$ by
\begin{align}\label{MMM}
M_{x,x'}((y^j)_{j<k}):= \sup_{j < k}   \frac{|y^j(x') - y^j(x)|}{\varepsilon}\,.
\end{align}
The `spring' connecting a pair $(x,x')\in \mathcal{\rm{NN}_{\varepsilon}}(\Omega)$ is considered \emph{elastic} (or \emph{undamaged}) if $\BBB \varepsilon^{-1} \EEE |y^j(x)-y^j(x')|$ does not exceed a given  \emph{threshold} $\bar{R}>1$ for all past times $t^j<t^k$, i.e., if $M_{x,x'}((y^j)_{j<k})\le \bar{R}$. In this case, the pair interaction can recover completely after removal of loading and thus  the energy contribution of $(x,x')$ is simply given by the stored elastic energy $W\bigl(\BBB \varepsilon^{-1} \EEE |y(x)-y(x')|\bigr)$. By contrast, if $M_{x,x'}((y^j)_{j< k  })>\bar{R}$, the spring is supposed to be \emph{damaged} and the complete deformation history $(y^j)_{j<k}$  should be reflected at present time.  With    \eqref{MMM}, this leads to the energy contribution with memory
\[W\big( M_{x,x'}((y^j)_{j<k})\big)      \vee   W\big( \varepsilon^{-1} |y(x)-y(x')|\big)\,,\]
where $a \vee b := \max \lbrace a,b \rbrace$ for $a,b \in \R$. We refer to \cite[(2.4)]{BFS} for more details and mention that here we adopted as memory variable the \emph{maximal opening}, as implemented in continuum settings, see e.g.\   \cite{Zanini, Giacomini:2005b} or  \cite[Section 5.2]{BFM}.  
Note, that also other choices would be conceivable, such as that of a cumulative increment, see e.g. \cite{BFM, Vito}. This, however, is not discussed in this work.  Summarizing, we can express the energy contribution of a pair $(x,x')\in \mathcal{\rm{NN}_{\varepsilon}}(\Omega) $ by 
\[F_{x,x'}(y,(y^j)_{j<k})\BBB := \EEE \begin{cases}
    W\bigl(\varepsilon^{-1} |y(x)-y(x')|\bigr) & \text{ if } M_{x,x'}((y^j)_{j<k})\leq \bar{R} \\
W\big(  M_{x,x'}((y^j)_{j<k})\big)      \vee   W\big( \varepsilon^{-1} |y(x)-y(x')|\big)  & \text{ if } M_{x,x'}((y^j)_{j<k})>\bar{R}\,.
\end{cases}\] 
Then, the assigned energy for a given  history $(y^j)_{j<k}$ of deformations  and an admissible deformation $y\colon\mathcal{L}_\varepsilon(\Omega)\to \R^3$ satisfying \eqref{y admi} (for some time-dependent boundary condition) is given by 
\begin{equation}\label{F:standardenergy}
    \mathcal{F}_{\varepsilon}(y,(y^j)_{j<k}):=\frac{1}{2}\, \sum_{(x,x') \in \mathcal{\rm{NN}_{\varepsilon}}(\Omega)} F_{x,x'}(y,(y^j)_{j<k})\,.
\end{equation}

\subsection*{Formulation for rescaled displacement fields}

Since we restricted ourselves to the case of anti-plane shear, the model can be rewritten in terms of the vertical displacement   $w \colon \mathcal{L}_{\varepsilon}(\Omega)\to \mathbb{R}$. In particular, \BBB given $y(x_1,x_2) = (x_1,x_2, w(x_1,x_2))$, \EEE for each $(x,x') \in  {\rm NN}_\eps(\Omega)$, using that $|x-x'| = \eps$ and defining  $\Psi(s):=W(\sqrt{1+s^2})$ for $s \in \R$, we get   
\begin{align}\label{eq: phiii}
W\Bigl(\frac{|y(x)-y(x')|}{\varepsilon}\Bigr)=W\Bigl(\sqrt{1+\Bigl(\frac{|w(x)-w(x')|}{\varepsilon}\Bigr)^2}\Bigr)= \Psi\Bigl(\frac{|w(x)-w(x')|}{\varepsilon}\Bigr)\,,
\end{align}
where the assumptions on $W$ from above lead to the following properties for $\Psi$: 
\JJJ 
\begin{enumerate}[label=(\roman*')]
    \item\label{phiprop1}$\Psi\geq 0$ and $\Psi(s)=0$ if and only if $s=0$.
        \item \label{phiprop3}   $\Psi$ is increasing on $(0,\infty)$.  
    \item\label{phiprop2} $\Psi$ is $C^2$ and Lipschitz on $\R$ with $\Psi'(0) =0$ and $\Psi''(0) = \mu >0$. 
    \item\label{phiprop4} $\lim_{|s|\to \infty}\Psi(s)=\kappa$\,.        
\end{enumerate}
\EEE
We now discuss the scaling of the energy in terms of different displacements $w$. One can check that the energy of a displacement of the form $w_1 = t\chi_{V}$, for $t \in  \R\setminus \{0\}$ 
and $V \subset \Omega$ with $\partial V \cap \Omega$ being smooth, scales like $\sim \eps^{-1}$ as each spring crossing $\partial V  \cap \Omega$ contributes $\Psi(t/\eps)\sim \kappa$ to the energy and the number of springs crossing $\partial V  \cap \Omega$ scales like $\eps^{-1}$. On the other hand, the energy of an affine displacement $w_2(x) = a_\eps \cdot x$, for $a_\eps \in \R^2$ possibly depending on $\eps$, scales like $\eps^{-2} \Psi(a_\eps)$ since $\#  {\rm NN}_\eps(\Omega) \sim \eps^{-2}$.

The most interesting  regime is when the contributions of the displacement with crack $w_1$ and the elastic displacement $w_2$ are of the same order. Since for small values $s$ we have \JJJ $\Psi(s) \sim s^2$ \EEE (see \ref{phiprop2}), this suggests $a_\eps \sim  \JJJ \eps^{1/2}\EEE$, i.e., elastic  displacements should be of (infinitesimal) small order. To account for this, we pass to \emph{rescaled displacements}  
$$u(x) := \JJJ \eps^{-1/2} w(x) \EEE \quad \text{ for $x \in  \L_\eps (\Omega)$. }$$
We write the energy \eqref{basicenergy} in terms of $u$ and rescale by $\eps$ to ensure that for the displacements $w_1 $ and  $w_2 $ considered above the energy is of order $1$. Using also \eqref{eq: phiii}, this leads to 
\begin{equation}\label{basicenergy-neu}
\mathscr{E}_{\varepsilon}(u):=\frac{\eps}{2}\, \sum_{(x,x') \in \mathcal{\rm{NN}_{\varepsilon}}(\Omega)} \Psi\Bigl(\frac{ \JJJ \eps^{1/2} \EEE |u(x)-u(x')|}{\varepsilon}\Bigr) .\end{equation}
It remains to transfer the memory variable to the setting of rescaled displacements. To this end, we observe that for each  pair $(x,x')\in \mathcal{\rm{NN}_{\varepsilon}}(\Omega)$, setting $R:=(\bar{R}^2-1)^{1/2}$, we have 
\[\frac{|y^j(x)-y^j(x')|}{\varepsilon}> \bar{R} \iff \frac{|w^j(x)-w^j(x')|}{\varepsilon}> R  \iff \frac{|u^j(x)-u^j(x')|}{\varepsilon}> \JJJ \eps^{-1/2} \EEE R\,,\]  
where $(w^j)_j$ and $(u^j)_j$ denote the history of (rescaled) displacements related to the deformations $(y^j)_j$.  We now define the $u$-dependent memory variable as 
\begin{align}\label{eq: memory variableXXX}
M_{x,x'}^\eps((u^j)_{j<k}):=   \JJJ \varepsilon^{-1/2} \EEE \sup_{j<k} |u^j(x)-u^j(x')|\,,
\end{align}
where we include $\eps$ in the notation to highlight the rescaled nature of the quantity. Defining
\[   E^\eps_{x,x'}  (u,(u^j)_{j<k}) \BBB := \EEE \begin{cases}
    \Psi\bigl(\varepsilon^{\JJJ -1/2 \EEE} |u(x)-u(x')|\bigr) &  \text{ if }  M^\eps_{x,x'}( \BBB (u^j)_{j<k} \EEE )\leq  {R} \\
    \Psi\bigl( M_{x,x'}^\eps((u^j)_{j<k})\bigr) \vee    \Psi\bigl(\varepsilon^{\JJJ -1/2 \EEE} |u(x)-u(x')|\bigr)  & \text{ if } M^\eps_{x,x'}((u^j)_{j<k})> {R},
\end{cases}  \] we can formulate the energy   \eqref{F:standardenergy}   by referring only to the rescaled displacement $u \colon \L_\eps(\Omega) \to \R$, namely  
\begin{equation}\label{energy-in-u}
\mathcal{E}_{\varepsilon}(u,(u^j)_{j<k}) \BBB := \EEE \frac{\varepsilon}{2}\, \sum_{ {\rm NN}_\eps(\Omega)} E^\eps_{x,x'} (u,(u^j)_{j<k})\,.\end{equation}
In this setting, the boundary condition in \eqref{y admi} transfers to requiring that displacements $u$ lie in the set
$$\mathcal{A}_\eps(g) := \big\{   u \colon \L_\eps(\Omega) \to \R\; \colon \, u(x) = g(x) \quad \text{ for all $x \in \L_\eps(\Omega \setminus U)$}  \big\},  $$
where $g :=  \eps^{\JJJ -1/2 \EEE} h \cdot e_3$.

\subsection{Main result: Convergence of quasi-static crack growth}

In this subsection, we state our main result on the atomistic-to-continuum convergence of crack growth. We start by introducing evolutions on the atomistic level.  To this end, we fix an arbitrary sequence $(\eps_n)_n \subset (0,\infty)$ with $\eps_n \to 0$ as $n \to \infty$. For simplicity, in the following we  write $\mathscr{E}_n$,  $\mathcal{E}_n$, $\L_n(\Omega)$, $\mathcal{\rm{NN}}_n(\Omega)$,  and $\mathcal{A}_{n}$   in place of $\mathscr{E}_{\eps_n}$, $\mathcal{E}_{\eps_n}$,   $\L_{\varepsilon_n}(\Omega)$, $\mathcal{\rm{NN}}_{\eps_n}(\Omega)$,  and $\mathcal{A}_{\varepsilon_n}$.  

%(Note that the regularity of $g$ is assumed for simplicity and could be relaxed at the expense of more refined estimates. We do not dwell on this point.)

We introduce a time discrete evolution \BBB which is \EEE driven by time-dependent boundary conditions $g\in W^{1,1}(0,T;W^{2,\infty}(\Omega))$.    We choose a sequence $(\delta_n)_n \subset (0,\infty)$ with $\delta_n \to 0$ and for each $\delta_n$ we consider the subdivision $0=t^0_n <\dots<t_n^{ T /\delta_n}=T$  of the interval $[0,T]$ with step size $\delta_n$. (Without restriction, we assume that $T /\delta_n \in \N$.)  Correspondingly, let $(g(t^k_n))_k$ be the sequence of boundary data at different time steps $k\in \{0,\dots, T/\delta_n \}$.   We suppose that the \emph{initial value} $u_n^0 \in \mathcal{A}_n(g(0) )$ is a minimum configuration in the sense that 
\begin{equation}\label{minimizing-scheme0}
u_n^{0}\in  {\rm  argmin} \big\{ \mathscr{E}_n(v) \colon v \in \mathcal{A}_n(g(0) ) \big\} \,,
\end{equation}
with $\mathscr{E}_n$ as given in \eqref{basicenergy-neu}. We inductively define an evolution as follows: given $(u_n^j)_{0 \le j \le k-1}$,   we  let 
\begin{equation}\label{minimizing-scheme}
    u_n^{k}\in {\rm argmin}\, \big\{\mathcal
{E}_n(v, (u_{n}^{j})_{j <  k})\colon \, v \in\mathcal{A}_n (g(t^{k}_n))   \big\}  \,,
\end{equation}  
where the minimization problem is influenced by $(u_n^{j})_{j< k}$ in terms of the memory variable given in \eqref{eq: memory variableXXX}, see \eqref{energy-in-u}.   We then define the evolution $u_n\colon [0,T]\times  \mathcal{L}_n(\Omega) \to \R$, discrete  \BBB  in space and piecewise constant in time, \EEE by
\begin{equation}\label{definitionofun}
u_n(t):= u^k_n \; \text{for}\; t\in[t^{k}_n,t^{k+1}_n)\,.
\end{equation}
The existence of minimizers in \eqref{minimizing-scheme} follows from the Direct Method and is briefly discussed in Lemma \ref{lemma: dm} below. We emphasize that we consider here an evolution discretized in time which is close to the spirit of the approach by {\sc Francfort and Marigo} \cite{frma98}. On the contrary, the existence of a continuous-in-time solution for the energy \eqref{energy-in-u} has been established in \cite{BFS} with an additional \BBB rate-dependent \EEE dissipation term in \eqref{minimizing-scheme} which allows to pass to the time-continuous limit for fixed $\eps_n$. 

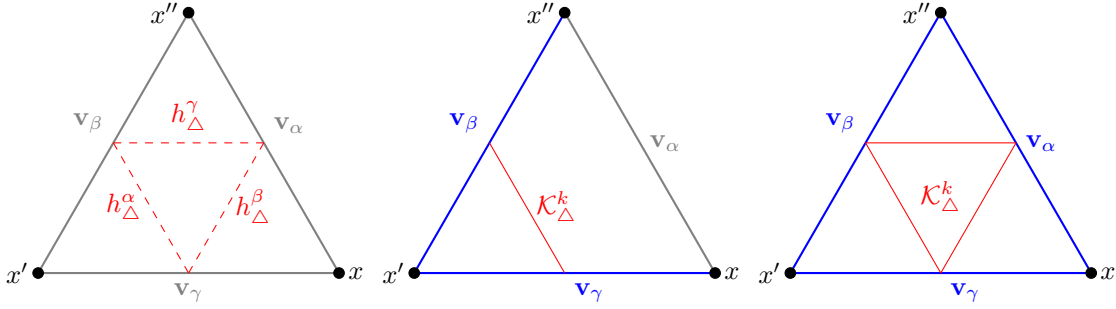
\begin{figure}
    \centering
    \begin{tikzpicture}
        \draw [gray, thick] (0:0)  -- node[anchor=south east] {$\mathbf{v}_\beta$}++(60:4) --node[anchor=south west] {$\mathbf{v}_\alpha$}++(-60:4) --node[anchor=north ] {$\mathbf{v}_\gamma$} cycle;
       \draw[dashed,red](0:2)--node[anchor=  east] {$h_{\triangle}^{\alpha}$}++(120:2) ;
       \draw[dashed,red](0:2)--node[anchor=  west] {$h_{\triangle}^{\beta}$}++(60:2)  ;
       \draw[dashed,red](60:2)--node[anchor= south] {$h_{\triangle}^{\gamma}$}++(0:2)  ;
            \filldraw[black] (0,0) circle (2pt) node[anchor=east]{$x'$};
       \filldraw[black] (0,0)+(60:4) circle (2pt) node[anchor=east]{$x''$};
       \filldraw[black] (0,0)+(0:4) circle (2pt) node[anchor=west]{$x$};
           \draw [blue, thick] (5,0)  -- node[anchor=south east] {$\mathbf{v}_\beta$}+(60:4);
       \draw [blue, thick] (5,0) --node[anchor=north west] {$\mathbf{v}_\gamma$}+(0:4);
       \draw [gray,thick] (5,0)++(60:4) --node[anchor=west ] {$\mathbf{v}_\alpha$} ++(-60:4); 
       \draw[red](7,0)--node[anchor=  west] { $\mathcal{K}^k_\triangle$}++(120:2);
       \filldraw[black] (5,0) circle (2pt) node[anchor=east]{$x'$};
       \filldraw[black] (5,0)+(60:4) circle (2pt) node[anchor=east]{$x''$};
       \filldraw[black] (5,0)+(0:4) circle (2pt) node[anchor=west]{$x$};
    \draw [blue, thick] (10,0)  -- node[anchor=south east] {$\mathbf{v}_\beta$}+(60:4);
       \draw [blue, thick] (10,0) --node[anchor=north west] {$\mathbf{v}_\gamma$}+(0:4);
       \draw [blue,thick] (10,0)++(60:4) --node[anchor=west ] {$\mathbf{v}_\alpha$} ++(-60:4); 
       \node[anchor= center, red] at (12,1) { $\mathcal{K}^k_\triangle$};
       \draw[red](12,0)-- ++(120:2);
        \draw[red](12,0)--%node[anchor=  west] { $\mathcal{K}^k_\triangle$}%
        ++(60:2);
        \draw[red](12,0)+(60:2)-- ++(120:2);
       \filldraw[black] (10,0) circle (2pt) node[anchor=east]{$x'$};
       \filldraw[black] (10,0)+(60:4) circle (2pt) node[anchor=east]{$x''$};
       \filldraw[black] (10,0)+(0:4) circle (2pt) node[anchor=west]{$x$};
   \end{tikzpicture}
   \caption{In the first figure, the segments $h_{\scalebox{1.0}{$\scriptscriptstyle \triangle$}}^\alpha$, $h_{\scalebox{1.0}{$\scriptscriptstyle \triangle$}}^\beta$, $h_{\scalebox{1.0}{$\scriptscriptstyle \triangle$}}^\gamma$  are depicted, where $\alpha,\beta,\gamma$ always denote pairwise distinct indices in $\lbrace 1,2,3 \rbrace$. In the second  and third   figure, the blue lines denote the `broken springs' and the red line the corresponding 'crack set' $\mathcal{K}^k_{\scalebox{1.0}{$\scriptscriptstyle \triangle$}}$. } 
\label{figure1}
\end{figure}

Next, we associate a `crack set` to the atomistic evolution by collecting the `broken springs'. Roughly speaking, in each triangle we draw a line between the centers of two 'broken springs'. \EEE To this end, we let $\T_n$ be the set of all \JJJ  closed \EEE triangles $\triangle$ contained in  $\L_n(\Omega)$ and consider $\triangle \in \mathcal{T}_n$ with vertices $x,x',x''$. For $\alpha,\beta,\gamma=1,2,3$ pairwise distinct, we let $h_{\scalebox{1.0}{$\scriptscriptstyle \triangle$}}^\alpha$ denote the segment in $\triangle$ between the centers of the two sides in $\mathbf{v}_\beta$ and $\mathbf{v}_\gamma$ direction, see Figure~\ref{figure1}. Given $t \in [t_n^k,t_n^{k+1})$, we say that $h_{\scalebox{1.0}{$\scriptscriptstyle \triangle$}}$ lies in $\mathcal{K}^k_{\scalebox{1.0}{$\scriptscriptstyle \triangle$}}$ if  $ h_{\scalebox{1.0}{$\scriptscriptstyle \triangle$}}  = h_{\scalebox{1.0}{$\scriptscriptstyle \triangle$}}^\alpha$ for some $\alpha=1,2,3$ and $\min\lbrace M_{x,x'}^\eps((u^j)_{j<k+1}), M_{x',x''}^\eps((u^j)_{j<k+1})  \rbrace >R$, where $x,x',x''$ are labeled such that $x-x' = \pm \eps_n \mathbf{v}_\beta$ and $x''-x' = \pm \eps_n \mathbf{v}_\gamma$.  Then, for each $t \in [t_n^k,t_n^{k+1})$ we define   
\begin{align}\label{eq: first K def}
K_n(t) :=  \bigcup_{\triangle \in \mathcal{T}_n} \bigcup_{h_{\scalebox{1.0}{$\scriptscriptstyle \triangle$}} \in  \mathcal{K}^k_{\scalebox{1.0}{$\scriptscriptstyle \triangle$}} }  h_{\scalebox{1.0}{$\scriptscriptstyle \triangle$}}.  
\end{align}

To formulate the main result, we introduce the \emph{limiting continuum energy}
\begin{align}\label{eq: lim-en}
\mathcal{E}(u,K):= \int_{\Omega}\Phi(\nabla u)\,{\rm d}x+ \int_{K}\varphi(\nu_K)\, {\rm d} \mathcal{H}^{1}\,,
\end{align}
for each $u \in SBV^{\JJJ 2 \EEE}(\Omega)$ and each rectifiable set $K\subset \R^2$ with $\mathcal{H}^1(K) < + \infty$, where $\nabla u$ denotes the approximate differential and  $S(u)$ the jump set of $u$, which is subject to the constraint $S(u) \, \tilde{\subset} \,  K$. (Here and in the following, $\, \tilde{\subset} \, $ stands for inclusions up to $\mathcal{H}^1$-negligible sets.) For basic notation and properties of $SBV$ functions, we refer the reader to \cite{Ambrosio-Fusco-Pallara:2000}. The energy features an \emph{elastic term} with energy density 
\begin{align}\label{eq: lineraition}
\Phi(z):=\JJJ\frac{ \mu }{\sqrt{3} }\EEE  \sum_{\mathbf{v}\in \mathcal{V}}|z  \cdot \mathbf{v} |^{\JJJ 2} \quad \text{ for $z\in \R^2$}, 
\end{align}
resulting from a suitable linearization of the potential $W$, and an anisotropic \emph{surface energy} with density 
\begin{align}\label{varphidef}
\varphi(\nu):=\frac{2\kappa}{\sqrt{3}}\sum_{\mathbf{v}\in \mathcal{V}}|  \nu  \cdot \mathbf{v}| \quad \text{ for $\nu \in  \mathbb{S}^1 := \lbrace x \in \R^2 \colon |x| = 1\rbrace$}
\end{align}
\BBB where we recall $\mathcal{V}:=\{\mathbf{v}_1,\mathbf{v}_2,\mathbf{v}_2-\mathbf{v}_1\}$. \EEE Note that the latter has \BBB already appeared \EEE  in the atomistic-to-continuum result \cite{FS153}.   Moreover,  by $AD(g,H)$ we denote all functions $v\in SBV^{\JJJ 2 \EEE}(\Omega)$ such that $v=g$ on $\Omega \setminus \overline{U}$ and $S(v)\, \tilde{\subset} \,  H$. In \cite{Francfort-Larsen:2003}  (see \cite{dMasoFranToad} for the case of anisotropic densities), the existence of an \emph{irreversible quasi-static crack evolution} with respect to the boundary displacement $g$ has been shown. This is a mapping  $t\to (u(t),K(t))$ with $u(t) \in AD(g(t),K(t))$ \BBB for all $t \in [0,T]$ \EEE such that the following three conditions hold:

\begin{itemize}
\item[(a)] \emph{Irreversibility}:   $K(t_1) \, \tilde{\subset} \,  K(t_2)$   for all $0\leq t_1\leq t_2\leq T$.
\item[(b)] \emph{Global stability}:  \BBB For every $t \in [0,T]$, given any  \EEE $K(t)\, \tilde{\subset} \,  H$ and $v\in AD(g(t),H)$     it holds that 
    \begin{equation}\label{finalstability}
    \mathcal{E}(u(t), K(t))\leq \mathcal{E}(v,H)\,.
    \end{equation}
    \item[(c)]   \emph{Energy balance}:    The function $t\mapsto \mathcal{E}(u(t), K(t))$ is absolutely continuous and it holds that
    \begin{equation}\label{energybalance}
        \frac{\rm d}{{\rm d}t}\mathcal{E}(u(t), K(t))=\int_{\Omega}D\Phi(\nabla u) \cdot \nabla \partial_{t}g(t)\, {\rm d}x \BBB \quad \text{for a.e.\ $t \in [0,T]$}\,, \EEE 
    \end{equation} 
    where by $\partial_t$ we denote the time derivative of $g$.
\end{itemize}
Eventually,  we need to define a notion of convergence. \BBB We define \EEE
 \begin{align}\label{tildeom}
\Omega_n:=\bigcup_{\triangle\in \mathcal{T}_n}\triangle.
\end{align}
For a sequence of discrete  displacements $u_n\colon \L_n(\Omega)\to \R$, \BBB we \EEE denote by $\tilde{u}_n\colon \Omega_n \to \R$ the interpolation of $u_n$ \BBB  which satisfies $\tilde{u}_n = {u}_n$ on $\mathcal{L}_n(\Omega)$ \EEE and is affine on each $\triangle\in \mathcal{T}_n$. 

\begin{definition}\label{define-convergence}
    Let $u_n \colon  \mathcal{L}_n(\Omega) \to \R$ be a sequence of atomistic displacements. We say that $u_n$ ${\rm AC}$-converges to $u\in SBV^{\MMM 2 \EEE}(\Omega)$ if $\chi_{\Omega_n}\tilde{u}_n\to u$ in $L^1(\Omega)$.
\end{definition} 

For the crack sets, we use the notion of \emph{$\sigma$-convergence} recalled in Definition \ref{def: sigma conv} below, which is a suitable notion of convergence for sets related to jump sets of $SBV$ functions. Our main result reads as follows.

\begin{theorem}[Atomistic-to-continuum crack growth]\label{maintheorem}
    There exists a quasi-static crack growth $t\to (u(t),K(t))$ with respect to the boundary condition $g$ such that, up to a subsequence, we have 
\begin{align}\label{eq: sigma coniiiii}
\text{$K_n(t)$ $\sigma$-converges to $K(t)$ for all $t \in [0,T]$ as $n \to \infty$},
\end{align}
    and for each $t \in [0,T]$ there exists a $t$-dependent subsequence $(n_k)_k$ such that  
    \begin{align}\label{eq: l1convi}
    \text{$u_{n_k}(t)$ AC-converges to $u(t)$ as $k \to \infty$.}
    \end{align}
    Moreover, for all $t \in [0,T]$ we have
\begin{align}\label{energ convi}
\mathcal{E}_n(u_n(t), (u_{n}^{j})_{j <  k(t)}) \to   \mathcal{E}(u(t),K(t)) \quad \text{as $n \to \infty$},    
\end{align}
    where for each $n$ the \BBB ($n$-dependent)  \EEE index  $k(t)$ is chosen such that $t \in [t^{k(t)}_n,t^{k(t)+1}_n)$.

\end{theorem}

\begin{remark}
{\normalfont
(i) In addition to the convergence of the energy in \eqref{energ convi}, it holds a separate convergence of suitably defined elastic and crack energies, see \eqref{en convo1}.  

(ii)  Besides the AC-convergence stated in \eqref{eq: l1convi}, also strong convergence of derivatives for \BBB suitably defined \emph{jump interpolations} \EEE introduced in Section \ref{sec: jump interp} holds. We refer to Lemma \ref{rem: the new one} for details.

(iii)  Note that the exact choice of the threshold $R$ does not affect the limiting densities $\Phi$ and $\varphi$  in  \eqref{eq: lim-en}. However, due to the memory variable on the atomistic level, the limiting quasi-static crack evolution might be sensitive to the choice of $R$.  

}
\end{remark}

The proof of this result is given in the next three sections and is structured as follows.  In Section~\ref{sec: prel}, we reformulate our problem in terms of  the triangles  of the lattice, introduce a suitable notion of elastic and crack energy on the atomistic level, and  define auxiliary interpolations. We also introduce the necessary notions of convergence  for crack sets,  and prove some preliminary results on compactness and lower semicontinuity. Subsequently, Section \ref{se:c main} contains the proof of Theorem \ref{maintheorem}. The most technical part consists in verifying a \emph{global stability property}, which will be analyzed in detail in Section \ref{stabilitysection}.

\section{Preliminaries}\label{sec: prel}

In this section \BBB we \EEE  collect some preliminaries that are needed  for our analysis. First, we will reformulate the energy \eqref{energy-in-u} in terms of triangles  and introduce an alternative interpolation for atomistic displacements that exhibits jumps \BBB relating \EEE to the 'atomistic crack set' defined in \eqref{eq: first K def}. We also recall compactness results for discrete displacements and lower semicontinuity results  for elastic energies which are essentially adaptations of statements in \cite{FS151}. Eventually, we \BBB recall \EEE two types of convergence of sets and discuss their main properties.

\subsection{Reformulation of the energy in terms of triangles}

We consider the set of all triangles $\T_n$ and the piecewise affine interpolation $\tilde{v}_n\colon \Omega_n \to \R$ of an atomistic displacement $v_n\colon\mathcal{L}_n(\Omega)\to \R$ introduced below \eqref{tildeom}. We also recall the set $\mathcal{V} = \{\mathbf{v}_1,\mathbf{v}_2,\mathbf{v}_2-\mathbf{v}_1\}$, representing the three directions of the triangular lattice. The gradient $\nabla \tilde{v}_n$ is constant on each $\triangle\in \T_n$, and for convenience we  denote this value by $({v}_n)_{\triangle}\in \R^2$. The rescaled distance between the two atoms $x,x'\in \mathcal{L}_n(\Omega)$ in the deformed configuration can then be expressed by
\[\varepsilon_n^{\MMM -1/2 \EEE}|v_n(x)-v_n(x')|=\varepsilon_n^{\MMM 1/2 \EEE}\Big|({v}_n)_{\triangle}\cdot \frac{x-x'}{\varepsilon_n}\Big|=\varepsilon_n^{\MMM 1/2 \EEE}|({v}_n)_{\triangle} \cdot \mathbf{v}|\,,\]
 where $\mathbf{v}\in\mathcal{V}$ is the  unit vector parallel to  $x-x'$. In fact, we can find a pair $(\triangle, \BBB \mathbf{v} \EEE )\in \T_n\times \mathcal{V}$, which unambiguously corresponds to a nearest-neighbor pair $\{x,x'\}$. \JJJ More precisely, there is a surjective mapping 
 $\Xi: \T_n\times \mathcal{V}  \times \lbrace -1,1\rbrace \to \mathcal{\rm{NN}}_n(\Omega) $, where $\Xi(\triangle, \mathbf{v}, s)$ is given by the unique pair $(x,x') \in \mathcal{\rm{NN}}_n(\Omega)$ such that $x,x'$ are vertices of $\triangle$ and $x-x'=s \mathbf{v}$\EEE. Note, however, except for springs at $\partial \Omega_n$, each pair $(x,x')\in \mathcal{\rm{NN}}_n(\Omega)$  is associated to two different triangles.

Next, we  express the memory variable in terms of  $(\triangle, \BBB  \mathbf{v} \EEE )\in \mathcal{T}_n\times \mathcal{V}$.  To this end,  given the time discretization $\{0=t^{0}_n< t^{1}_n<\dots< t^{T /\delta_n}_n =T\}$ of the interval $[0,T]$ with step size $\delta_n$, we suppose that for  $t^k_n\in [0,T]$ the displacements $(u^j_{n})_{j<k}$ at previous time steps $(t^j_n)_{j<k}$ have already been found. Then, recalling \eqref{eq: memory variableXXX}, we define
\begin{align}\label{eq: noacons}
   \Me{k}:=M_{x,x'}^{\eps_n} ((u^j_n)_{j<k}) =\varepsilon_n^{\MMM -1/2 \EEE}\,\sup_{j< k} \,|u^j_{n}(x)-u^{j}_n(x')|\,.
   \end{align}
In the sequel, we frequently need to keep track of the memory variable including competitors in the next time step. To this end, for an arbitrary displacement $v_n \colon \L_n(\Omega) \to \R$ we \BBB define \EEE 
\begin{align}\label{eq: noacons-neu}
\Met{k}{v_n}  :=\Me{k}  \vee \varepsilon_n^{\MMM 1/2 \EEE}| ({v}_n)_{\triangle} \cdot \mathbf{v}| 
\end{align}
 which includes  the impact of $v_n$ on the spring in $\triangle$ pointing in direction $\mathbf{v}\in\mathcal{V}$ into the memory variable.

%  \[\mathcal{C}_n^{k}:=\{\triangle \in \mathcal{T}_{n}\colon \Me{k} >R  \quad \text{for at least one}\, \mathbf{v}\in\mathcal{V}\}\].
 
We now use this to define the notion of a 'broken triangle'. We \BBB regard \EEE a triangle $\triangle \in \mathcal{T}_n$ as `broken' if at least one edge has been ``broken'', i.e., stretched above the threshold $R$. This means, for a given displacement history  
$(u^j_{n})_{j < k}$ and an arbitrary displacement $v_{n} \colon \mathcal{L}_{n}(\Omega)\to \R$, we define
\begin{align}\label{nom jump tri}
\mathcal{C}^{k}_{n}&:= \big\{\triangle\in \T_{n}\colon \Me{k} > R \quad \text{for at least one}\; \mathbf{v}\in \BBB \mathcal{V} \EEE \big\}\,, \notag \\  
\mathcal{C}^{k}_{n}(v_n)&:= \big\{\triangle\in \T_{n}\colon \Met{k}{v_n} >R\quad \text{for at least one}\; \mathbf{v}\in \mathcal{V} \big\}\,.
\end{align}  
Note that, whenever  a triangle is not broken, i.e., $\triangle \notin \mathcal{C}^{k}_{n}(v_n)$, we have by \eqref{eq: noacons-neu} \BBB that \EEE
 \begin{align}\label{eq: not so uniform}
 |( v_{n})_\triangle| \BBB = \big(  |( v_{n})_\triangle \cdot \mathbf{v}_1   |^2 +  |( v_{n})_\triangle \cdot \tfrac{1}{\sqrt{3}}(\mathbf{v}_2 + \mathbf{v}_3)    |^2   \big)^{1/2}         \EEE  \le 2 R \eps_n^{\MMM -1/2 \EEE}.
\end{align}         
Now, we can express the energy \eqref{energy-in-u} at time $t^{k}_n$ in terms of  triangles and their corresponding edges. In fact, we can separate the `broken' triangles $\mathcal{C}_{n}^{k}(v_n)$ from the `intact triangles' $ \mathcal{T}_{n} \setminus \mathcal{C}_{n}^{k}(v_n)$, and obtain  
\begin{equation}
    \begin{aligned}\label{energy-triangle}
    \mathcal{E}_{n}(v_{n}; (u^j_n)_{j < k}) = &\frac{\varepsilon_n}{2}\sum_{\BBB \triangle \in \EEE\T_{n}\setminus \Cn{k}{v_n}} \sum_{\mathbf{v}\in\mathcal{V}} \Psi(\varepsilon_n^{\MMM 1/2 \EEE}|({v}_n)_{\triangle} \cdot \mathbf{v}|) \;+ %\frac{\varepsilon}{2}\sum_{\Cn{k}{v_n}} \sum_{\mathbf{v}\in\mathcal{V}} \Psi(\Met{k}{v_n})
        \frac{\varepsilon_n}{2}\sum_{\BBB \triangle \in \EEE \Cn{k}{v_n}} \sum_{\substack{\mathbf{v}\in\mathcal{V} \\ \Met{k}{v_n}  \leq R}} \Psi(\varepsilon_n^{\MMM 1/2 \EEE}\,|({v}_n)_{\triangle}  \cdot \mathbf{v}|) + \\ & \ \ \ +  \frac{\varepsilon_n}{2}\sum_{\BBB \triangle \in \EEE\mathcal{C}_n^{k}(v_n)} \sum_{\substack{\mathbf{v}\in\mathcal{V}\\ \Met{k}{v_n}>R}} \Psi( \Met{k}{v_n})  
+  \mathcal{E}_n^{k, {\rm{bdy}}}(v_n),
    \end{aligned}
\end{equation}
where $\mathcal{E}_n^{k,{\rm{bdy}}} $ contains the energy contributions of nearest-neighbor pairs $(x,x')$ that are part of only one triangle and consequently only reflected with a prefactor $1/2$ in the other terms of \eqref{energy-triangle}.  More precisely, we have
\begin{align}\label{bdy energy}
  \mathcal{E}_n^{k, {\rm{bdy}}}(v_n ) \EEE  = \frac{\varepsilon_n}{2} \sum_{\substack{ (x,x') \in {\rm NN}_n(\Omega), \\  x,x' \in  \partial \Omega_n}}  E^{\eps_n}_{x,x'} (v_n,(u^j_n)_{j< k}), 
\end{align}
\JJJ where we recall the definition in \eqref{tildeom}. \EEE Note that the prefactor $\frac{1}{2}$ appears again, because in the above sums all springs that do not lie at the boundary are counted twice, being associated to two different triangles. For lower bounds, we can ignore $\mathcal{E}_n^{ {\rm{bdy}}}$ as it is clearly nonnegative. For upper bounds, however, we will need to show explicity that the term vanishes in the limit $n \to \infty$, see Theorem \ref{stability} for details. 

From now on, we simplify the notation in the following sense, whenever no confusion arises: instead of writing $\triangle \in \mathcal{T}_n'$ below sums for a subset $\mathcal{T}_n' \subset \mathcal{T}_n$, we simply write $\mathcal{T}_n'$ indicating that the sum is taken over all triangles in $\mathcal{T}_n'$. In a similar fashion, in sums  over $\mathcal{V}$ involving $\Met{k}{v_n}$, we do not include $ \mathbf{v} \in \mathcal{V}$ in the notation.

In  \eqref{energy-triangle}, we have split the energy in terms of intact and broken triangles. We further split the broken triangles in a crack part, which will lead to jumps in the limiting description as $n \to \infty$, and a remainder term. First, we choose a sequence of diverging thresholds $(R_n)_n$ satisfying
\begin{align}\label{rnepsn}
R_n \to + \infty \quad \quad \text{and} \quad \quad \JJJ R_n^{3/2 }\varepsilon_n^{1/4} \to 0  \EEE ,
\end{align}
and consider the set of triangles
\begin{align}\label{big jump tri}
\Cl{k}{v_n} := \{\triangle\in  \mathcal{C}^{k}_n(v_n)\colon \Met{k}{v_n}  > 2 R_n   \; \text{for at least one $\mathbf{v}\in\mathcal{V}$}\}.
\end{align}
Here, the supercript L indicates that the memory variable is `large'. By elementary geometry this implies $ \Met{k}{v_n}  >  R_n$ for at least two $ \mathbf{v} \in \mathcal{V}$, and we split the set into
\begin{align}\label{2,3}
 \Cl{k}{v_n}= \mathcal{C}^{k,2}_n(v_n)  \,\dot{\cup}\, \mathcal{C}^{k,3}_n(v_n) \,,
 \end{align}
corresponding to triangles where exactly two or three $\mathbf{v}  \in \mathcal{V}$ satisfy $ \Met{k}{v_n}  >  R_n$. Then, we can write  
\begin{equation}\label{energy:main}
\begin{aligned}
    \mathcal{E}_{n}(v_n; (u^j_n)_{j < k}) &= \mathcal{E}_n^{k, {\rm{ela}}}(v_n) + \mathcal{E}_n^{k, {\rm{cra}}}(v_n) + \mathcal{E}_n^{k, {\rm{rem}}}(v_n) +\mathcal{E}_n^{k, {\rm{bdy}}}(v_n),  
\end{aligned}
\end{equation}
where $\mathcal{E}_n^{k, {\rm{bdy}}}$ is defined in \eqref{bdy energy}, and we set 
\begin{align}\label{eq: different enegies}
\mathcal{E}_n^{k, {\rm{ela}}}(v_n) & =    \frac{\varepsilon_n}{2}\sum_{\T_{n}\setminus \Cn{k}{v_n}} \sum_{\mathbf{v}\in\mathcal{V}} \Psi(\varepsilon_n^{\MMM 1/2 \EEE}|({v}_n)_{\triangle} \cdot \mathbf{v}|),\notag\\
\mathcal{E}_n^{k, {\rm{cra}}}(v_n) & = \frac{\varepsilon_n}{2}\sum_{\Cl{k}{v_n}} \ \ \  \sum_{\Met{k}{v_n}>R_n} \Psi(\Met{k}{v_n}), \notag \\
\mathcal{E}_n^{k, {\rm{rem}}}(v_n) & =\frac{\varepsilon_n}{2}\sum_{\Cn{k}{v_n}}\Big( \sum_{R< \Met{k}{v_n}\le R_n} \Psi(\Met{k}{v_n}) + \sum_{\Met{k}{v_n}  \leq R} \Psi(\varepsilon_n^{\MMM 1/2 \EEE}\,|({v}_n)_{\triangle}  \cdot \mathbf{v}|) \Big)\notag \\
& \ \ \ + \frac{\varepsilon_n}{2}\sum_{\Cn{k}{v_n} \setminus \Cl{k}{v_n}}   \sum_{\Met{k}{v_n}> R_n} \Psi(\Met{k}{v_n}).  
\end{align}
We will see that the contributions of the \emph{remainder term} $\mathcal{E}_n^{k, {\rm{rem}}}(v_n)$   will eventually vanish in the limit $n \to \infty$. Still, from a technical point of view, we need delicate arguments to control \BBB its \EEE contribution in all estimates below.   Our next goal is to rewrite the \emph{elastic energy} $\mathcal{E}_n^{k, {\rm{ela}}}(v_n)$ and the \emph{crack energy}  $\mathcal{E}_n^{k, {\rm{cra}}}(v_n)$ as an integral functional which will allow us to apply lower semicontinuity results in Sobolev and $SBV$ spaces.  This will be achieved by introducing a jump interpolation which is \BBB inspired by \EEE  the one introduced in \cite{FS151}.

\subsection{Jump interpolation and representation of the energy}\label{sec: jump interp}

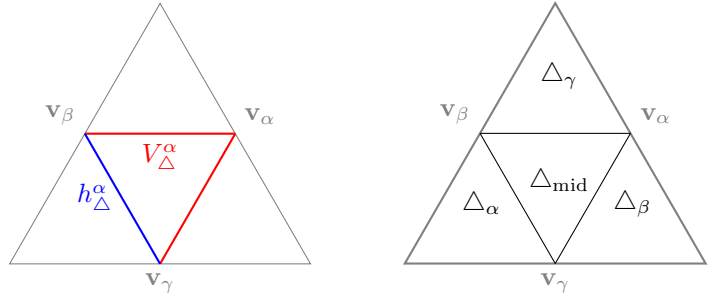
\begin{figure}
    \begin{tikzpicture}
        \draw [gray] (0,0)  -- node[anchor=south east] {$\mathbf{v}_\beta$}++(60:4) --node[anchor=south west] {$\mathbf{v}_\alpha$}++(-60:4) --node[anchor=north ] {$\mathbf{v}_\gamma$} cycle;
        \draw[red,thick](2,0)--++(60:2)  ;
        \draw[red,thick](2,0)+(60:2)-- node[anchor=  north] {$V^\alpha_{\triangle}$} ++(120:2);
        \draw[blue,thick](2,0)-- node[anchor=  east] {$h^\alpha_{\triangle}$}++(120:2);
    \end{tikzpicture}
    \hspace*{1cm}
    \begin{tikzpicture}
        \draw [gray, thick] (0:0)  -- node[anchor=south east] {$\mathbf{v}_\beta$}++(60:4) --node[anchor=south west] {$\mathbf{v}_\alpha$}++(-60:4) --node[anchor=north ] {$\mathbf{v}_\gamma$} cycle;
        \draw[black](0:2)--node[anchor=  east] {}++(120:2) ;
        \draw[black](0:2)--node[anchor=  west] {}++(60:2)  ;
        \draw[black](60:2)--node[anchor= south] {}++(0:2)  ;
        --cycle; 
        \node[anchor=center] at (2.05, 1.1) {$\triangle_{\rm mid}$};
        \node[anchor= east] at (1.4,0.8) {$\triangle_\alpha$};
        \node[anchor= east] at (3.4,0.8) {$\triangle_\beta$};
        \node[anchor=center] at (2.05,2.5) {$\triangle_\gamma$};
    \end{tikzpicture}

    \caption{The first picture shows the sets $V^\alpha_{\scalebox{1.0}{$\scriptscriptstyle \triangle$}}$  and $h^\alpha_{\scalebox{1.0}{$\scriptscriptstyle \triangle$}}$ \EEE which are used for the \BBB definition \EEE of $K_n^{k,i}$. The second one  depicts the subdivision of a triangle into sub-triangles. }\label{different-triangles}
\end{figure}

We now define (three versions of) a  `crack set' related to $\Cl{k}{v_n}$. To this end, as before, given $\triangle \in \mathcal{T}_n$ and  $\alpha,\beta,\gamma=1,2,3$ pairwise distinct,  we let $h^\alpha_{\scalebox{1.0}{$\scriptscriptstyle \triangle$}}$ denote the segment in $\triangle$  between the centers of the both sides in direction  $\mathbf{v}_\beta$ and $\mathbf{v}_\gamma$, see Figure~\ref{figure1}. Note that $h^\alpha_{\scalebox{1.0}{$\scriptscriptstyle \triangle$}}$ is parallel to $\mathbf{v}_\alpha$. We define the set $V^\alpha_{\scalebox{1.0}{$\scriptscriptstyle \triangle$}} =h^{\beta}_{\scalebox{1.0}{$\scriptscriptstyle \triangle$}} \cup h^{\gamma}_{\scalebox{1.0}{$\scriptscriptstyle \triangle$}}$, \BBB see Figure~\ref{different-triangles}, \EEE and  for each $i=1,2,3$, we let
\begin{align}\label{eq: hat sets}
\Kl{k,i}{v_n} := \bigcup_{\triangle \in \mathcal{C}^{k,3}_n(v_n)}  V^i_{\scalebox{1.0}{$\scriptscriptstyle \triangle$}}  \ \cup \  \bigcup_{\triangle \in \mathcal{C}^{k,2}_n(v_n)}  h_{\scalebox{1.0}{$\scriptscriptstyle \triangle$}},  \quad \quad \quad  K^{k,{\rm L}}_n(v_n) :=  \bigcup\nolimits_{i=1}^3 \Kl{k,i}{v_n},
\end{align}
where for each $\triangle \in  \mathcal{C}^{k,2}_n(v_n)$  we choose $h_{\scalebox{1.0}{$\scriptscriptstyle \triangle$}} = h_{\scalebox{1.0}{$\scriptscriptstyle \triangle$}}^m$, where $m$ is the unique index with $ M^{k, \mathbf{v}_m}_{ n, \triangle}(v_n) \le  R_n$.  The reason for using three different variants lies in the derivation of a lower semicontinuity result for the length of the sets, and follows the ideas in \cite{FS151}.   Note carefully that on triangles $\triangle \in \mathcal{C}^{k,2}_n(v_n)$ all variants of $\Kl{k,i}{v_n}$ coincide.  
 
\begin{remark}\label{remark-about-2Ks}
{\normalfont
Recalling \eqref{eq: first K def}, one can check that the solution $u^k_n$ at time step $t_n^k$, see \eqref{minimizing-scheme}, fulfills
\begin{align*}
 K^{k,{\rm L}}_n(u^k_n)   \subset K_n(t)\quad \text{for $t\in [t_n^k,t_n^{k+1})$}\,,
\end{align*} 
and the difference of the sets is due to the triangles  $\triangle \in \Cn{k}{u^{k}_n}$ where the numbers $\mu_{n,\triangle}^k = \# \lbrace \mathbf{v} \in \mathcal{V} \colon  M_{n,\triangle}^{k,\mathbf{v}} (u^k_n) > R\rbrace$  and $\nu_{n,\triangle}^k = \# \lbrace \mathbf{v} \in \mathcal{V} \colon  M_{n,\triangle}^{k,\mathbf{v}}  (u^k_n) > R_n\rbrace$ differ. Indeed, if $\mu_{n,\triangle}^k = \nu_{n,\triangle}^k =3$, then both  $K_n(t)$ and $K^{k,{\rm L}}_n(u^k_n)$ coincide with $\bigcup_{i=1}^3 \BBB h_{\scalebox{1.0}{$\scriptscriptstyle \triangle$}}^{i} \EEE $ on $\triangle$. If $\mu_{n,\triangle}^k = \nu_{n,\triangle}^k =2$, then there exists one spring in direction $\mathbf{v}_i$ such that $M_{n,\triangle}^{k,\mathbf{v_i}} (u^k_n) \le R$, and both  $K_n(t)$ and $K^{k,{\rm L}}_n(u^k_n)$ coincide with $h_{\scalebox{1.0}{$\scriptscriptstyle \triangle$}}^i$ on $\triangle$. Along the proof, we will see that the difference of $K_n(t)$ and $K^{k,{\rm L}}_n(u^k_n) $ is asymptotically vanishing for $n \to \infty$.
}
\end{remark}

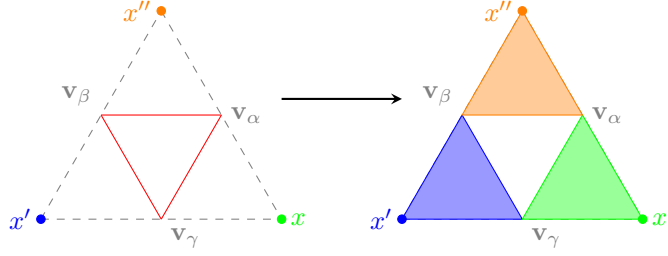
\begin{figure}
    \begin{tikzpicture}[scale=0.8]
 
        \draw [gray, dashed] (0,0)  -- node[anchor=south east] {$\mathbf{v}_\beta$}+(60:4);
        \draw [gray, dashed] (0,0) --node[anchor=north west] {$\mathbf{v}_\gamma$}+(0:4);
        \draw [gray,dashed] (0,0)++(60:4) --node[anchor=west ] {$\mathbf{v}_\alpha$} ++(-60:4); 
        \draw[red](2,0)--++(60:2)  ;
        \draw[red](2,0)+(60:2)-- ++(120:2)  ;
        \draw[red](2,0)-- ++(120:2)  ;
        \filldraw[blue] (0,0) circle (2pt) node[anchor=east]{$x'$};
        \filldraw[orange] (0,0)+(60:4) circle (2pt) node[anchor=east]{$x''$};
        \filldraw[green] (0,0)+(0:4) circle (2pt) node[anchor=west]{$x$};
        \draw [thick,-stealth](4,2) -- (6,2);
        
        \draw [gray, dashed] (6,0)  -- node[anchor=south east] {$\mathbf{v}_\beta$}+(60:4);
        \draw [gray, dashed] (6,0) --node[anchor=north west] {$\mathbf{v}_\gamma$}+(0:4);
        \draw [gray,dashed] (6,0)++(60:4) --node[anchor=west ] {$\mathbf{v}_\alpha$} ++(-60:4);

    \filldraw[fill=blue!40!white,draw=blue]{} (6,0)--++(60:2) --++ (-60:2) --cycle;
    \filldraw[fill=green!40!white,draw=green]{} (8,0)--++(60:2) --++ (-60:2)--cycle;
    \filldraw[fill=orange!40!white,,draw=orange]{} (10,0)++(120:4)--++(-60:2) --++ (0:-2) --cycle;
        \filldraw[blue] (6,0) circle (2pt) node[anchor=east]{$x'$};
        \filldraw[orange] (6,0)+(60:4) circle (2pt) node[anchor=east]{$x''$};
        \filldraw[green] (6,0)+(0:4) circle (2pt) node[anchor=west]{$x$};
    \end{tikzpicture}
    \caption{Construction of $\hat{v}_n^{k}$ for a triangle $\triangle \in \mathcal{C}^{k,3}_n(v_n)$ with three broken springs.} \label{3crack-construction}
\end{figure}

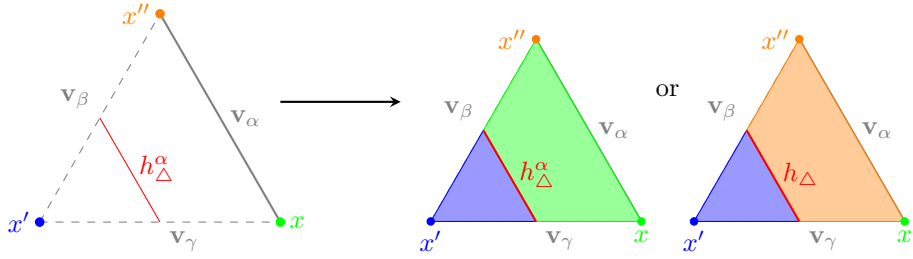
\begin{figure}
\begin{tikzpicture}[scale=0.8]
    \draw [gray, dashed] (0,0)  -- node[anchor=south east] {$\mathbf{v}_\beta$}+(60:4);
    \draw [gray, dashed] (0,0) --node[anchor=north west] {$\mathbf{v}_\gamma$}+(0:4);
    \draw [gray,thick] (0,0)++(60:4) --node[anchor=west ] {$\mathbf{v}_\alpha$} ++(-60:4); 
    \draw[red](2,0)--node[anchor=  west] {$h^\alpha_{\triangle}$}++(120:2)  ;
    \filldraw[blue] (0,0) circle (2pt) node[anchor=east]{$x'$};
    \filldraw[orange] (0,0)+(60:4) circle (2pt) node[anchor=east]{$x''$};
    \filldraw[green] (0,0)+(0:4) circle (2pt) node[anchor=west]{$x$};

    % Pfeil zwischen den Dreicken
    \draw [thick,-stealth](4,2) -- (6,2);
\end{tikzpicture}
% Erste Option der Interpolation
\begin{tikzpicture}[scale=0.7]
    \centering
    \draw [gray, dashed] (6,0)  -- node[anchor=south east] {$\mathbf{v}_\beta$}+(60:4);
    \draw [gray, dashed] (6,0) --node[anchor=north west] {$\mathbf{v}_\gamma$}+(0:4);
    \draw [gray,thick] (6,0)++(60:4) --node[anchor=west ] {$\mathbf{v}_\alpha$} ++(-60:4); 
    \draw[red](8,0)--node[anchor=  west] {$h^\alpha_\triangle$}++(120:2)  ;
 
\filldraw[fill=blue!40!white,draw=blue]{} (6,0)--++(60:2) --++ (-60:2)--cycle;
\filldraw[fill=green!40!white, draw=green]{} (8,0)--++(0:2) --++ (120:4) --++ (-120:2)  --cycle;
\draw[red,thick](8,0)--node[anchor=  west] {$h^\alpha_\triangle$}++(120:2)  ;

\filldraw[blue] (6,0) circle (2pt) node[anchor=north]{$x'$};
\filldraw[orange] (6,0)+(60:4) circle (2pt) node[anchor=east]{$x''$};
\filldraw[green] (6,0)+(0:4) circle (2pt) node[anchor=north]{$x$};
    
% Second option for interpolation in triangles with two broken springs: 
\node at (10.5,2.5) {or};

\draw [gray, dashed] (11,0)  -- node[anchor=south east] {$\mathbf{v}_\beta$}+(60:4);
\draw [gray, dashed] (11,0) --node[anchor=north west] {$\mathbf{v}_\gamma$}+(0:4);
\draw [gray,thick] (11,0)++(60:4) --node[anchor=west ] {$\mathbf{v}_\alpha$} ++(-60:4); 
    \filldraw[fill=blue!40!white,draw=blue]{} (11,0)--++(60:2) --++ (-60:2)--cycle;
    \filldraw[fill=orange!40!white, draw=orange]{} (13,0)--++(0:2) --++ (120:4) --++ (-120:2)  --cycle;
    \draw[red,thick](13,0)--node[anchor=  west] {$h_\triangle$}++(120:2)  ;

    \filldraw[blue] (11,0) circle (2pt) node[anchor=north]{$x'$};
    \filldraw[orange] (11,0)+(60:4) circle (2pt) node[anchor=east]{$x''$};
    \filldraw[green] (11,0)+(0:4) circle (2pt) node[anchor=north]{$x$};

\end{tikzpicture}
\caption{Construction of $\hat{v}_n^{k}$ in a triangle $\triangle \in \mathcal{C}^{k,2}_n(v_n)$ with $M_{n,\triangle}^{k,\mathbf{v}_\alpha}(v_n)\leq R_n$.} \label{2crack-construction}
\end{figure}

We now introduce a jump interpolation $\hat{v}_n^{k} \colon \Omega_n \to \R$, that exhibits jumps on triangles $\triangle\in \mathcal{C}^{k}_n(v_n)$ and leaves $\tilde{v}_n$ unchanged on $\T_n\setminus \mathcal{C}^{k}_n(v_n)$. On $\mathcal{T}_n\setminus \Cn{k}{v_n}$, we simply set $\hat{v}^{k}_n=\tilde{v}_n$. Next, we consider triangles $\triangle \in \Cl{k}{v_n}$. We divide $\triangle$ in four equilateral sub-triangles, where we denote by $\triangle_{\alpha}$ the triangle that is formed by taking $h_{\scalebox{1.0}{$\scriptscriptstyle \triangle$}}^\alpha$ as one side, and the middle triangle with sides $h_{\scalebox{1.0}{$\scriptscriptstyle \triangle$}}^\alpha$, $h_{\scalebox{1.0}{$\scriptscriptstyle \triangle$}}^\beta$, and $h_{\scalebox{1.0}{$\scriptscriptstyle \triangle$}}^\gamma$ is denoted by $\triangle_{\rm mid}$, see Figure \ref{different-triangles}. Now, if $\triangle \in \mathcal{C}^{k,3}_n(v_n)$,  we define $\hat{v}_n^{k}$ such that $\hat{v}_n^{k}=v_n$ on the three vertices of $\triangle$ and $\nabla \hat{v}_n^{k} \equiv 0$ on $\bigcup_{i=1}^{3} \triangle_{i}$,  see   Figure \ref{3crack-construction}. On $\triangle_{\rm mid}$, let $\hat{v}_n^{k}$ be constant, attaining an arbitrary value.  For  $\triangle\in \mathcal{C}^{k,2}_n(v_n)$, assuming $M^{k,\mathbf{v}_\alpha}_{n,\triangle}(v_n) \le  R_n$, we define $\hat{v}^{k}_n$ in such a way that   $\hat{v}^{k}_n$ is continuous on ${\rm int}(\triangle) \setminus h^\alpha_{\scalebox{1.0}{$\scriptscriptstyle \triangle$}}$, $\nabla \hat{v}^{k}_n\equiv0$ in $\triangle$, and  $\hat{v}^{k}_n$ coincides with $v_n$ on two of the vertices, see Figure \ref{2crack-construction} (we choose one of the two options.)  For  $\triangle\in  \Cn{k}{v_n}  \setminus  \Cl{k}{v_n}$,  we  define  $\hat{v}^{k}_n$ such that $\hat{v}^{k}_n$ is constant on ${\rm int}(\triangle) $ and coincides with $v_n$ in one of the three vertices (we choose one of the three options).   Note that the interpolation exhibits jumps also at the boundary of broken triangles, both on $ \Cl{k}{v_n}$ \BBB and on \EEE $\Cn{k}{v_n}  \setminus  \Cl{k}{v_n}$. We \BBB cover \EEE this set by
\begin{align}\label{eq: small jumpi}
{K}^{k,{\rm S}}_n(v_n) := \bigcup_{\Cn{k}{v_n}}  \partial \triangle,
\end{align}
where we use the superscript S to indicate that, as we will see in Remark \ref{good remark} below, the jump height on this set is small. We also observe that, in view of \eqref{eq: not so uniform},  we have 
\begin{align}\label{eq: not so uniform, but neu}
\Vert \nabla  \hat{v}^{k}_{n}\Vert_\infty \le  2R \eps_n^{\MMM -1/2 \EEE}.
\end{align} In fact, either $\nabla \hat{v}^{k}_{n} = 0$, namely on broken triangles, or we have $\nabla  \hat{v}^{k}_{n} = ( v_n)_\triangle$ constant on a triangle $\triangle$.

We can now use this interpolation to obtain a convenient representation of the elastic energy. By the construction of $\hat{v}^{k}_{n}$, we can take the sum over all  $\T_{n}$ in the {elastic energy}, see \eqref{eq: different enegies}, and write this term as an  integral over $\Omega_{n}$. Precisely, we have  
\begin{align}\label{eq: ela en}
 \mathcal{E}_n^{k, {\rm{ela}}}(v_{n}) =\frac{4}{\sqrt{3}\varepsilon_n}\int_{\Omega_n} \frac{1}{2} \sum_{\mathbf{v}\in\mathcal{V}}\Psi(\varepsilon_n^{\MMM 1/2 \EEE}|\nabla \hat{v}^{k}_n(x)\cdot \mathbf{v}|)\,{\rm d}x =   \int_{\Omega_{n}}  \Psi_{n}^{\rm cell}(  \nabla \hat{v}^{k}_n \EEE (x)) \,{\rm d}x\,,
 \end{align}
where we used that the area of each triangle is $\frac{\sqrt{3}}{4}\eps_n^2$ and we defined  the \emph{cell energy density} $ \Psi_{n}^{\rm cell}$   by 
 \begin{align}\label{psitriangle}
 \Psi_{n}^{\rm cell}(z):= \frac{2}{\sqrt{3}\eps_n}\sum_{\mathbf{v}\in\mathcal{V}}\Psi\big(\eps_n^{\MMM 1/2 \EEE}|z\cdot \mathbf{v}|\big) \quad \text{ for all $z \in \R^2$}.
 \end{align}
Concerning the  {crack energy} defined in \eqref{eq: different enegies}, a crucial point in our analysis is the characterization of $\mathcal{E}_n^{k, {\rm{cra}}}(v_{n})$ in terms of the sets $\Kl{k,i}{v_n}$ introduced in \eqref{eq: hat sets}. To formulate this relation, we define $\varphi_\alpha(\nu):=(|\mathbf{v}_\beta\cdot\nu|+|\mathbf{v}_\gamma\cdot \nu|)$ for $\alpha,\beta,\gamma \in \lbrace 1,2,3 \rbrace $ pairwise distinct  and for $\nu \in \mathbb{S}^1$. Also recall the definition of $\kappa$ \JJJ in \ref{define-potential1}. \EEE

\begin{proposition}\label{prop: repri}
 There exists a universal constant $C>0$ such that, for given displacements $(u^j_{n})_{j<k}$ and an arbitrary 
 $v_n \colon \L_n(\Omega) \to \R$,   the crack energy $\mathcal{E}_n^{k, {\rm{cra}}}(v_n)$ satisfies 
\begin{equation}\label{boundary-jump0}
\Big| \mathcal{E}_n^{k, {\rm{cra}}}(v_{n}) -   \frac{\kappa}{\sqrt{3}}\sum_{i=1}^{3}\int_{{K}^{k,i}_n(v_n)}  \varphi_i(\nu_K) \,{\rm d}\mathcal{H}^1 \Big|  \le      C|\kappa - \Psi(R_n)| \eps_n \# \mathcal{C}^{k,L}_n(v_n),
 \end{equation}  
 where $\nu_K$ denotes a unit normal to ${K}^{k,i}_n(v_n)$. 
% Moreover, the jump interpolation $\hat{v}_n^{k}$ satisfies 
\end{proposition} 

\begin{remark}\label{good remark}
{\normalfont
(i)  In view of \eqref{rnepsn} and the properties of $\Psi$, the right-hand side of \eqref{boundary-jump0} will vanish for $n \to \infty$ as we will indeed have a bound on $\# \mathcal{C}^k_n(v_n)$ of the order $\eps_n^{-1}$. \\
(ii) One can check that the interpolation $\hat{v}_n^{k}$  and the crack set are related in the sense that $S(\hat{v}_n^{k}) \subset  {K}^{k,{\rm L}}_n(v_n) \cup \BBB {K}^{k,{\rm S}}_n(v_n) \EEE$ with 
\begin{equation}\label{boundary-jump}
|[\hat{v}_n^{k}(x)] |\leq  C\varepsilon_n^{\MMM 1/2 \EEE} R_n  \quad \quad \text{for all $x \in {K}^{k,{\rm S}}_n \BBB (v_n) \EEE $}.
 \end{equation}   
Indeed, the inclusion follows directly from the construction of $\hat{v}_n^{k}$. Property  \eqref{boundary-jump} can be derived by using the continuity of $\hat{v}_n^{k}$ at (some of) the vertices of the triangles along with  \eqref{eq: not so uniform, but neu},   and the fact that for every pair $x,x' \in \triangle$ such that the segment between $x$ and $x'$ does not intersect $ {K}^{k,{\rm L}}_n(v_n)$ we have $\eps_n^{\MMM -1/2 \EEE}|v_n(x) - v_n(x')| \le 2R_n$, see \eqref{big jump tri}.  

(iii) The choice in \eqref{rnepsn} also shows that \JJJ  $R_n\varepsilon_n^{ 1/2 }\to 0$ \EEE and hence \eqref{boundary-jump} implies that the jumps in ${K}^{k,{\rm S}}_n$ are asymptotically of vanishing height. \BBB (This is the reason why this set is indicated with label `small', see below  \eqref{eq: small jumpi}.) \EEE Therefore, these jumps, artificially introduced by our interpolation, will not affect the continuum limit.

}
\end{remark}

\begin{proof}[Proof of Proposition \ref{prop: repri}]
    The proof consists in showing that
    \begin{align}\label{eq: FFFF}
        \mathcal{G}_n^{k, {\rm{cra}}}(v_n):= \kappa  \frac{\varepsilon_n }{2}    \big( 2 \# \mathcal{C}_n^{k,2}(v_n) +   3\# \mathcal{C}_n^{k,3}(v_n)\big)= \frac{\kappa}{\sqrt{3}}\sum_{i=1}^{3}\int_{K^{k,i}_n(v_n) }  \varphi_i(\nu_K) \,{\rm d}\mathcal{H}^1      \,.
    \end{align}
    Recalling  \eqref{2,3},  \eqref{eq: different enegies}, and \ref{phiprop3}, the statement then follows from the fact that
    \[\mathcal{E}_n^{k, {\rm{cra}}}(v_{n})\geq \Psi(R_n)  \frac{\varepsilon_n }{2}    \big( 2 \# \mathcal{C}_n^{k,2}(v_n) +   3 \# \mathcal{C}_n^{k,3}(v_n)\big),\] and therefore
    $$|\mathcal{G}_n^{k, {\rm{cra}}}(v_n)  - \mathcal{E}_n^{k, {\rm{cra}}}(v_{n})   | \le C|\kappa - \Psi(R_n)| \eps_n \# \mathcal{C}^{k,L}_n(v_n).$$
 To show \eqref{eq: FFFF}, we will consider the sets $\mathcal{C}^{k,2}_{n}(v_n)$ and $\mathcal{C}^{k,3}_{n}(v_n)$ separately.
     
    Firstly, consider  $\triangle\in \mathcal{C}^{k,3}_{n}(v_n)$.  As before, we denote by $h^{\alpha}_{\scalebox{1.0}{$\scriptscriptstyle \triangle$}}$ the segment between the midpoints of the  sides in directions $\mathbf{v}_\beta, \mathbf{v}_\gamma\in \mathcal{V}$  for pairwise distinct $\alpha,\beta,\gamma\in \{1,2,3\}$. We express the contribution inside $\triangle$ \BBB as \EEE an integral. 
    Recall that $\mathcal{H}^1(h^{\alpha}_{\scalebox{1.0}{$\scriptscriptstyle \triangle$}})   =\frac{\varepsilon_n}{2}$ for all $\alpha=1,2,3$, and    note that for the normal $\nu_\alpha$ of $h^{\alpha}_{\scalebox{1.0}{$\scriptscriptstyle \triangle$}}$ we have $|\mathbf{v}_\alpha\cdot \nu_\alpha|=0$,  and $|
        \mathbf{v}_\beta\cdot \nu_\alpha|=|
        \mathbf{v}_\gamma\cdot \nu_\alpha|=\frac{\sqrt{3}}{2}$.  Given $\alpha=1,2,3$,   the set $h^{\beta}_{\scalebox{1.0}{$\scriptscriptstyle \triangle$}}   \cup h^{\gamma}_{\scalebox{1.0}{$\scriptscriptstyle \triangle$}}  = V^{\alpha}_{\scalebox{1.0}{$\scriptscriptstyle \triangle$}}\EEE$ is exactly the part of  $K^{k,\alpha}_n(v_n)$ inside  $\triangle$. Denoting the normal to ${K}^{k,\alpha}_n(v_n)$ again by  $\nu_{K}$, we can rewrite  
    \begin{align}\label{ident1}
    \frac{3}{2}  =\frac{1}{\sqrt{3}\varepsilon_n}\sum_{\alpha = 1}^{3}\int_{{K}^{k,\alpha}_n(v_n) \cap \mathrm{int}(\triangle)}  (|\mathbf{v}_\beta\cdot\nu_K|+|\mathbf{v}_\gamma\cdot \nu_K| ) \,{\rm d}\mathcal{H}^1 = \frac{1}{\sqrt{3}\varepsilon_n}\sum_{\alpha=1}^{3}\int_{{K}^{k,\alpha}_n(v_n) \cap \mathrm{int}(\triangle)}  \varphi_\alpha(\nu_K) \,{\rm d}\mathcal{H}^1 \,.
        \end{align}
    Now consider $\triangle\in \mathcal{C}^{k,2}_{n}(v_n)$.   We choose $\mathbf{v}_m \in \mathcal{V}$ with  $M^{k,\mathbf{v}_m}_{n,\triangle}(v_n) \le R_n$. Then, by construction  $K^{k,\alpha}_n(v_n) \cap {\rm int}(\triangle)$ coincides with $h^{m}_{\scalebox{1.0}{$\scriptscriptstyle \triangle$}}$ for each $\alpha=1,2,3$. Therefore, we have $\mathcal{H}^1(K^{k,\alpha}_n(v_n) \cap {\rm int}(\triangle)) = \frac{\eps_n}{2}$ for $\alpha=1,2,3$ and, since   $|\nu_K\cdot \mathbf{v}_m|=0$ and $|\nu_K\cdot \mathbf{v}_\alpha|=\frac{\sqrt{3}}{2}$ for $\alpha\in \lbrace 1,2,3\rbrace \setminus \lbrace m \rbrace$, we can check  that 
    \begin{align}\label{ident2}
    \frac{1}{\sqrt{3}\varepsilon_n}\sum_{\alpha=1}^{3}\int_{{K}^{k,\alpha}_n(v_n) \cap \mathrm{int}(\triangle)}  \varphi_\alpha(\nu_K) \,{\rm d}\mathcal{H}^1 = \frac{1}{\sqrt{3}\varepsilon_n}\sum_{\alpha=1}^{3}\int_{{K}^{k,\alpha}_n(v_n) \cap \mathrm{int}(\triangle)}  (|\mathbf{v}_\beta\cdot\nu_K|+|\mathbf{v}_\gamma\cdot \nu_K| ) \,{\rm d}\mathcal{H}^1  = 1. 
     \end{align}
    Collecting   \eqref{ident1}--\eqref{ident2}   we find \eqref{eq: FFFF}. 
        \end{proof}

\subsection{Compactness and lower semicontinuity of elastic energies}\label{sec: comp}

In this subsection, we present a compactness result for sequences with bounded energy $\mathcal{E}_{n}$ (see \eqref{energy-triangle})  with respect to the notion of  AC-convergence introduced in  Definition \ref{define-convergence}. Afterwards, we present a lower semicontinuity result for the elastic energy. We start with elementary properties of $\Psi_{n}^{\rm cell}$ defined in \eqref{psitriangle} and $\Phi$ given in \eqref{eq: lineraition}. 

\begin{lemma}[Properties of $\Psi_{n}^{\rm cell}$]\label{lem:help2}
   {\rm (i)}  There exists a function  $\omega:\R^2 \to \R$ with  $\frac{\omega(y)}{|y|^{\JJJ 2 \EEE}}\to 0$ for $|y|\to 0$ such that 
    \begin{equation*}
    \Psi_{n}^{\rm cell}(y)=\Phi(y)+\eps_n^{-1}\omega\big(\eps_n^{\MMM 1/2 \EEE} \BBB y \EEE \big) \quad \quad \text{for all $y \in \R^2$ and all $n \in \N$}\,.
        \end{equation*}
 {\rm (ii)}   Let $r>0$. Then, there exists a constant $c_r>0$ depending on $r$ such that for all $n \in \N$  we have   
    \[ c_r |y|^{\MMM 2}\leq   \Psi_{n}^{\rm cell}(y), \quad \quad  c_r |y| \leq   |D\Psi_{n}^{\rm cell}(y)|  \quad \quad \text{for all  $y \in \R^2$ with  $|y|\leq r\eps_n^{\MMM -1/2 \EEE }$.} \] 
{\rm (iii)}  Then, there exists a constant $C>0$ such that for all $n \in \N$  we have 
$$ |  \Psi_{n}^{\rm cell}(y)| \le C|y|^{\MMM 2 \EEE },\quad \quad  | D\Psi_{n}^{\rm cell}(y)| \le C \, |y|  \quad  \quad   \quad  \text{for all $y \in \R^2$}. $$      
{\rm (iv)} For each $M >0$ and  sequences $(y_n^1)_n, (y_n^2)_n \subset \R^2$ with $|y_n^1| \le M$, $|y_n^2| \le M$, and       $|y_n^1 - y_n^2| \to 0$, it holds that $|D\Psi_{n}^{\rm cell}(y_n^1) - D\Psi_{n}^{\rm cell}(y_n^2)| \to 0$ as $n \to \infty$. 
\end{lemma}  

\begin{proof}
    We start with a preliminary computation.  Recalling the definition  $\Psi(s) =  W(\sqrt{1+s^2})$ for $s \in \R$ and the fact that $W$ is minimized at $1$ with $W(1) = 0$ \JJJ and $\lim_{r \searrow 1 }W'(r)=\mu$, \EEE a Taylor expansion yields   
        \begin{align}  \label{lineariii}\JJJ
        \Psi(s)= & \JJJ W\big(\sqrt{1+s^2}\big)=W(1)+\mu \frac{s^2}{2}   + \mathcal{O}(s^3) \,, \EEE
        \end{align}
    as $s \to 0$, and similarly \JJJ 
        \begin{align}\label{lineariii2}
      \Psi'(s) =  \mu\, s+ \BBB \mathcal{O}(s^2) \,. \EEE
      \end{align}    \EEE 
        Thus, by \eqref{lineariii} there exists $\eta >0$ such that \JJJ 
    \begin{align}\label{prelliii}
    \Psi(s)  \ge \frac{\mu}{4}s^2, \quad \Psi'(s)  \ge \frac{\mu}{2}\, s  \quad  \text{ for all $s\in \R$ with $|s| \le \eta$ }. 
    \end{align}    \EEE 
    First, putting together \eqref{lineariii} and \eqref{psitriangle}, we get \JJJ 
\begin{align}\label{repriiii}
\Psi_{n}^{\rm cell}(y):= \frac{\mu}{\sqrt{3}} \sum_{\mathbf{v}\in\mathcal{V}}
    |y\cdot \mathbf{v}|^2+ \varepsilon^{-1}_n \sum_{\mathbf{v}\in\mathcal{V}}\mathcal{O}\big( (\eps_n^{1/2}   |y\cdot \mathbf{v}|)^3 \big),
    \end{align} \EEE
     and thus 
    (i) is an immediate consequence of the definition of $\Phi$ in \eqref{eq: lineraition}.  
    % \eqref{lineariii} along with the definitions of $\Phi$ and $\Psi_{n}^{\rm cell}$ in \eqref{eq: lineraition} and \eqref{psitriangle}.
    
   We now come to the proof of (ii).  Recalling that $\mathbf{v}_1 = (1,0)^T$, $\mathbf{v}_2 = (\frac{1}{2},\frac{\sqrt{3}}{2})$,  and  $\mathbf{v}_3 = (-\frac{1}{2},\frac{\sqrt{3}}{2})$  we calculate for each $y \in \R^2$  
        \[y_1^2=|y \cdot \mathbf{v}_1 |^2, \quad   y_2^2 =     \frac{1}{3}|  y \cdot (\mathbf{v}_2 + \mathbf{v}_3)  |^2 \le \frac{2}{3} \big( |  y \cdot \mathbf{v}_2  |^2 + |  y \cdot \mathbf{v}_3  |^2 \big)\,, \]\EEE 
        which yields \JJJ $|y|^2\leq C(|y\cdot \mathbf{v}_1|^2+|y\cdot \mathbf{v}_2|^2+|y\cdot \mathbf{v}_3|^2)$ for a constant $C>0$. \EEE  \\
        Now let $y \in \R^2$ with $|y| \le r\eps_n^{\MMM -1/2 \EEE }$. Without restriction we assume that $\eta \le r$ and set $y' = \frac{\eta}{r} y$. We then have $|\varepsilon_n^{\MMM 1/2 \EEE}y'\cdot \mathbf{v}_i|\leq \eta$ for all $i\in \{1,2,3\}$, which by \eqref{prelliii} yields $|\varepsilon_n^{\MMM 1/2 \EEE}y'\cdot \mathbf{v}_i|^{\MMM 2 \EEE}\leq C \Psi(|\varepsilon_n^{\MMM 1/2 \EEE}y'\cdot \mathbf{v}_i|)$ for a constant $C>0$. Then, the previous estimate along with the definition of $\Psi_{n}^{\rm{cell}}$ gives us 
        % and \eqref{prelliii} yields for a constant $C>0$
        \begin{equation*}
        |y'|^{\MMM 2 \EEE}\leq C\eps_n^{-1} ( |\eps_n^{\MMM 1/2 \EEE} y' \cdot \mathbf{v}_1|^{\MMM 2 \EEE}+ |\eps_n^{\MMM 1/2 \EEE} y' \cdot \mathbf{v}_2|^{\MMM 2 \EEE}+|\eps_n^{ \MMM 1/2 \EEE} y' \cdot \mathbf{v}_3|^{\MMM 2 \EEE}) \le C\Psi_n^{\rm cell }( y') \le C\Psi_n^{\rm cell }(y) \,, 
        \end{equation*}        
        where the last step follows from the fact that $\Psi$ is increasing, see \ref{phiprop3}. This shows the first part of (ii) for a constant depending on $r$. The second part of (ii) can be shown along similar lines, using the second estimate in \eqref{prelliii}.  
        
     Now, we prove (iii).  We compute \JJJ   
    \begin{align}\label{eq: Dpsi}
     D\Psi_{n}^{\rm cell}(y)=  \frac{2}{\sqrt{3}}\eps_n^{  -1/2 }\sum_{\mathbf{v}\in\mathcal{V}}\Psi'\big(\eps_n^{1/2   }|y \cdot \mathbf{v}|\big)\frac{y \cdot \mathbf{v}}{|y \cdot \mathbf{v}|}\,  \mathbf{v} 
     \end{align} \EEE 
      for all $y \in \R^2$.  Thus, using \eqref{lineariii2}, there exist $C>0$ and $\eta >0$ such that \JJJ $|D\Psi_{n}^{\rm cell}(y)| \le C|y|$ for all $y \in \R^2$  \EEE with $|y| \le  \eta \eps_n^{\MMM -1/2 \EEE }$. If instead \BBB $\eta \eps_n^{\MMM -1/2 \EEE } < |y|$,  we can use $\Psi' \in L^\infty([0,\infty))$ by the Lipschitz continuity of $\Psi$  \EEE to derive that \JJJ 
     $$ |D\Psi_{n}^{\rm cell}(y)|  \le C \eps_n^{ -1/2  } =    C (|y|^{-1}\eps_n^{-1/2})  |y| \le C |y|.  $$ \EEE
The bound on $\Psi_{n}^{\rm cell}$ follows in a similar way, using \eqref{psitriangle}, \eqref{repriiii}, and $\Psi \in L^\infty([0,\infty))$.     Eventually, (iv) follows directly from \eqref{lineariii2}  and \eqref{eq: Dpsi}. This concludes the proof.   
    \end{proof}

We now proceed with the compactness result.

\begin{proposition}[Compactness]\label{bounded-jump}
Let $(v_n)_n$ be a sequence of displacements with $v_n \in \mathcal{A}(g_n)$ for a sequence of boundary conditions  $(g_n)_n$ such that $g_n$ AC-converges to some $g \in W^{2,\infty}(\Omega)$. For each $n \in \N$, let $(u^j_{n})_{j < k_n}$  be a displacement history for $k_n \in \N$  depending on $n$. Suppose that $\mathcal{E}_{n}(v_n; (u^j_{n})_{j < k_n}) + \|\tilde{v}_n\|_{L^{\infty}(\Omega_n)} \le M$ for some $M>0$. Then, there exists $C>0$ only depending on $M$ such that 
\begin{align}\label{eq: jump bound}
\#\mathcal{C}_{n}^{k_n}(v_n) \leq \frac{C}{\varepsilon_n} \quad \quad \text{ and } \quad \quad \mathcal{H}^{1}(S({\hat{v}^{k_n}_{n}}))\leq C\,.
\end{align}
Moreover,  there exists a function $v\in SBV^{\MMM 2 \EEE}(\Omega)$ with $v=g$ on $\Omega \setminus \overline{U}$ such that, up to a subsequence (not relabeled), $v_n$  AC-converges to $v$ and  for the jump interpolations  $\hat{v}^{k_n}_n$  we have 
\begin{align}\label{eq: jump bound2}
\chi_{\Omega_{n}} \hat{v}^{k_n}_n \to v \; \text{in}\;\; L^1(\Omega) \quad \text{and}\quad \chi_{\Omega_{n}}\nabla \hat{v}^{k_n}_{n} \rightharpoonup \nabla v \;\;in\; L^{\MMM 2 \EEE }(\Omega;\R^2)\,.
\end{align}   
 \end{proposition} 

\begin{proof}
For convenience, we simply write $k$ in place of $k_n$ in the proof. In view of \eqref{energy-triangle} and \ref{phiprop3}, we have  
\[\frac{\varepsilon_n}{2}\,\# \mathcal{C}_{n}^{k}(v_n)\, \Psi(R) \leq \mathcal{E}_{n}(v_n; (u^j_{n})_{j < k_n}) \le M < + \infty\,.\] 
Hence,  there exists some $C>0$ such that for all $\varepsilon_n >0$ it holds that 
\[\#\mathcal{C}_{n}^{k}(v_n) \leq \frac{C}{\varepsilon_n} \,.\]
 Recalling that by construction $\hat{v}^{k}_n$ jumps only on triangles  $\triangle \in \mathcal{C}_{n}^{k}(v_n)$ and we always have $S(\hat{v}^{k}_n) \cap \overline{\triangle} \le  \frac{9}{2} \eps_n$, we conclude the proof of  \eqref{eq: jump bound}.  The  bound  \eqref{eq: not so uniform, but neu} along with Lemma \ref{lem:help2}(ii) for $r = 2R $ and \eqref{eq: ela en}  shows
$$\Vert \nabla  \hat{v}^{k}_{n} \Vert_{L^{\MMM 2 \EEE}(\Omega_n)}^{\MMM 2 \EEE} = \int_{\Omega_n} |\nabla  \hat{v}^{k}_{n}(x)|^{\MMM 2 \EEE}  \, {\rm d}x \leq C  \int_{\Omega_n} \Psi_{n}^{\rm cell}(\nabla \BBB \hat{v}^{k}_{n} \EEE (x)) \,{\rm d}x = C \,\mathcal{E}_n^{k, {\rm{ela}}}(v_{n}) \le CM. $$
% $$\Vert \nabla  \hat{v}^{k, i}_{n} \Vert_{L^4(\Omega_n)}^4 = \eps_n^{-1} \Vert  \eps_n^{1/4}  \nabla \hat{v}^{k, i}_{n} \Vert_{L^4(\Omega_n)}^4 \le  C\eps_n^{-1} \int_{\Omega_n} \Psi_\triangle(\eps_n^{1/4} \nabla  \hat{v}^{k, i}_{n}) \, {\rm d}x   \le  C  \mathcal{E}_n^{k, {\rm{ela}}}(v_{n}) \le CM. $$
Thus, by assumption on $v_n$ and \eqref{eq: jump bound} we  derive 
$$  \Vert \nabla  \hat{v}^{k}_{n} \Vert_{L^{\MMM 2 \EEE}(\Omega_n)} +      \mathcal{H}^{1}(S(\hat{v}^{k}_{n}))   +  \|  \hat{v}^{k}_{n}  \|_{L^{\infty}(\Omega_n)} \le C \quad \text{ for all $n \in \N$}. $$ 
Moreover, \BBB as $\Omega$ is a Lipschitz set, \EEE one can check that $\sup_n\mathcal{H}^1(\partial \Omega_n) < +\infty$.  Thus, we can apply Ambrosio's compactness theorem for $SBV$ functions  \cite[Theorem 4.8]{Ambrosio-Fusco-Pallara:2000} to obtain a function $v\in SBV^{\MMM 2 \EEE}(\Omega)$ such that  $\chi_{\Omega_{n}} \hat{v}^{k}_{n}\to v$ in $L^1(\Omega)$ and  $\chi_{\Omega_{n}}\nabla \hat{v}_n^{k}  \rightharpoonup\nabla v$ weakly in $L^{\MMM 2 \EEE}(\Omega;\R^2)$. 

By $\mathcal{L}^2(\bigcup_{\mathcal{C}^{k}_{n}({v_n})}  \triangle) \le \frac{\sqrt{3}}{4}\eps_n^2  \# \mathcal{C}^{k}_{n}({v_n})  $, \eqref{eq: jump bound}, \BBB the bound on  $\|\tilde{v}_n\|_{L^{\infty}(\Omega_n)}$, \EEE and the fact that $\tilde{v}_n = \hat{v}_n^{k}$ on $\Omega_n \setminus \bigcup_{\mathcal{C}^{k}_{n}({v_n})} \triangle$,   we eventually also  get $\chi_{\Omega_{n}}\tilde{v}_n \to v$ in $L^1(\Omega)$, i.e., $v_n$ AC-converges to $v$. Since $g_n$ AC-converges to $g$ and $v_n \in \mathcal{A}(g_n)$, we  conclude the proof by observing   $v=g$ on $\Omega \setminus \overline{U}$.
\end{proof}

\BBB Now, \EEE we   prove a lower semicontinuity result for  the elastic energy. 
\begin{lemma}[Lower semicontinuity]\label{elastic-lowersim}
Let $(v_n)_n$ be a sequence of displacements with displacement history  $(u^j_{n})_{j < k_n}$ for $k_n \in \N$, satisfying  $\mathcal{E}_{n}(v_n; (u^j_{n})_{j < k_n}) + \|\tilde{v}_n\|_{L^{\infty}(\Omega_n)} \le M$ for some $M>0$. Denote by   $v\in SBV^{\MMM 2 \EEE}(\Omega)$ the AC-limit given in Proposition \ref{bounded-jump}. Then,  we have
    \begin{equation*}%\label{el-lowersim-ineq}
    \liminf_{n \to \infty} \mathcal{E}_n^{k_n, {\rm{ela}}}(v_{n})  \geq \int_{\Omega} \Phi(\nabla v)\, {\rm d}x\,.\end{equation*} 
\end{lemma}
\begin{proof}
For convenience, we again simply write $k$ in place of $k_n$.  Let $\omega$ be the function from Lemma \ref{lem:help2}(i) satisfying  $\frac{\omega(y)}{|y|^{\MMM 2 \EEE}}\to 0$ for $|y|\to 0$. Then, by \eqref{eq: ela en},   by Lemma \ref{lem:help2}(i), and the fact that $\Phi$ is homogeneous of degree ${\MMM 2 \EEE}$, we   have
    \begin{equation*}%\label{elasticsplit}
\mathcal{E}_n^{k, {\rm{ela}}}(v_{n})  =  \int_{\Omega_{n}} \Psi_{n}^{\rm cell}( \nabla \hat{v}^{k}_n(x)) \,{\rm d}x   =\int_{\Omega_{n}}\Bigl( \Phi(\nabla \hat{v}^{k}_{n}) +\frac{1}{\varepsilon_n}\omega(\varepsilon_n^{{\MMM 1/2 \EEE}} \nabla \hat{v}^{k}_{n})\Bigr)\, {\rm d}x.
\end{equation*}
% \RRR From here on the proof needs to be changed a bit because it only works if you cut out a set where the gradients are too high. Look at other papers where the linearization is done. \EEE    
Consider the function $\chi_{n}(x)\defas \chi_{[0,\varepsilon_n^{-1/{\MMM 4\EEE}})}(|\nabla \hat{v}^{k}_n(x)|)$ and note by \eqref{eq: jump bound2} that for $n\to \infty$ we have $\chi_{n}\to 1$ in measure on $\Omega$. In particular, we can estimate
\[\mathcal{E}_n^{k, {\rm{ela}}}(v_{n})\geq \int_{\Omega_{n}}\chi_{n}(x)\Bigl( \Phi(\nabla \hat{v}^{k}_{n}(x)) +\frac{1}{\varepsilon_n}\omega(\varepsilon_n^{{\MMM 1/2 \EEE}} \nabla \hat{v}^{k}_{n}(x))\Bigr)\, {\rm d}x\,.\]
Note that we can express the second part of the integrand as 
    \[\chi_n(x)\frac{1}{\varepsilon_n}\omega(\varepsilon_n^{{\MMM 1/2 \EEE}} \nabla \hat{v}^{k}_{n})=\chi_n(x)|\nabla \hat{v}^{k}_{n}|^{\MMM 2 \EEE}\frac{\omega(\varepsilon_n^{\MMM 1/2 \EEE} \nabla \hat{v}^{k}_{n})}{|\varepsilon_n^{\MMM 1/2 \EEE} \nabla \hat{v}^{k}_{n}|^{\MMM 2 \EEE}}\,.\]
% By definition of $\omega$ we now deduce that $\chi_n(x)|\nabla \hat{v}^{k,i}_{n}|^4\frac{\omega(\varepsilon_n^{1/4} \nabla \hat{v}^{k,i}_{n})}{|\varepsilon_n^{1/4} \nabla \hat{v}^{k,i}_{n}|^4}$ converges uniformly to $0$ as 
% \[\sup_{n\in \N}\Big
% \{\frac{\omega(\varepsilon_n^{1/4}\nabla \hat{v}^{k,i}_{n})}{|\varepsilon_n^{1/4}\nabla \hat{v}^{k,i}_{n}|^4}\colon\,|\nabla \hat{v}^{k,i}_{n}|\leq  \varepsilon_n^{-1/8} \Big\}=\sup_{n\in \N}\Big
% \{\frac{\omega(\varepsilon_n^{1/4}\nabla \hat{v}^{k,i}_{n})}{|\varepsilon_n^{1/4}\nabla \hat{v}^{k,i}_{n}|^4}\colon\,\varepsilon_n^{1/4}|\nabla \hat{v}^{k,i}_{n}|\leq  \varepsilon_n^{1/4} \Big\} \to 0\,. \]  
By definition of $\omega$ and \BBB $\chi_n$ \EEE we deduce that $\chi_n(x)\frac{\omega(\varepsilon_n^{\MMM 1/2 \EEE} \nabla \hat{v}^{k}_{n})}{|\varepsilon_n^{\MMM 1/2 \EEE} \nabla \hat{v}^{k}_{n}|^{\MMM 2 \EEE}}$ converges uniformly to $0$. Since by \eqref{eq: jump bound2} we have $\|\nabla \hat{v}^{k}_{n}\|_{L^{\MMM 2 \EEE}(\Omega_n)}\leq C$, this yields
    \begin{equation*}
    \int_{\Omega_{n}}\chi_n(x)\frac{1}{\varepsilon_n}\omega(\varepsilon_n^{\MMM 1/2 \EEE} \nabla \hat{v}^{k}_{n})  \, \BBB {\rm d}x \EEE\to 0 \quad \BBB \text{as $n \to \infty$}\,. \EEE \end{equation*}
Since $\chi_{\Omega_{n}}\nabla \hat{v}^{k}_{n} \wto \nabla v$ weakly in $L^{\MMM 2 \EEE}(\Omega;\R^2)$, we have $\chi_n\chi_{\Omega_{n}}\nabla \hat{v}^{k}_{n} \wto \nabla v$ weakly in $L^{\MMM 2 \EEE}(\Omega;\R^2)$. Because $\Phi$  defined in \eqref{eq: lineraition} is convex with $\Phi(0)=0$, we get 
\begin{align}\label{getsomenumber}
\liminf_{n \to \infty} \mathcal{E}_n^{k, {\rm{ela}}}(v_{n})= \liminf_{n \to \infty} \int_{\Omega_{n}} \chi_n\Phi(\nabla \hat{v}^{k}_{n})\, {\rm d}x \geq \int_{\Omega} \Phi(\nabla v)\, {\rm d}x\,.
\end{align}
    This concludes the proof. 
\end{proof}

Now, we will show that the lower semicontinuity result along with energy convergence implies strong convergence. For this we follow a standard argument, see e.g.\ \cite{DalMasoNegriPercivale:02, Schmidt:2009, Schmidt:08} or \cite{almi, FS153} for results in the context of fracture. More precisely, we have  the following. 
\begin{lemma}[Strong convergence]\label{rem: the new one}
Assume the setting of Lemma \ref{elastic-lowersim}. If we have
\begin{align}\label{energy convi}
\lim_{n \to \infty} \int_{\Omega_n} \Psi_n^{\rm cell}(\nabla \hat{v}^{k_n}_{n} ) \, {\rm d}x =    \lim_{n \to \infty} \mathcal{E}_n^{k_n, {\rm{ela}}}(v_{n})  = \int_{\Omega} \Phi(\nabla v)\, {\rm d}x, 
   \end{align}
then we deduce 
\begin{align}\label{energy convi2}
\chi_{\Omega_{n}}\nabla \hat{v}^{k_n}_{n} \to \nabla v \; \text{ strongly in } \;\; L^{\MMM 2 \EEE}(\Omega;\R^2)\,.
\end{align}   
\end{lemma}
\begin{proof}
\BBB As before, we   simply write $k$ in place of $k_n$. \EEE We first check that $h_n := |\chi_{\Omega_{n}}\nabla \hat{v}^{k}_{n}|^{\MMM 2 \EEE}$ is equi-integrable on $\Omega$. If not, we \BBB would \EEE have
$$\lim_{M \to \infty} \limsup_{n \to \infty} \int_{ \lbrace h_n > M \rbrace} h_n \, {\rm d}x \ge \eta $$
for some $\eta>0$. By a diagonal argument we then find a sequence $M_n \to \infty$ with   $  M_n \le  \varepsilon_n^{\MMM -1 \EEE} $  such that
$$ \limsup_{n \to \infty} \int_{ \lbrace h_n > M_n \rbrace} h_n \, {\rm d}x \ge \eta \,.$$
 With $\bar{\chi}_n := \chi_{ \lbrace h_n \le M_n \rbrace}$, we have $\bar{\chi}_n \to 1$ a.e.\ on $\Omega$ by \eqref{eq: jump bound2}. Then, repeating the lower semicontuinuity proof \BBB in \eqref{getsomenumber} \EEE and using Lemma \ref{lem:help2}(ii)  together with \eqref{eq: not so uniform, but neu}    we get 
\begin{align*}
\BBB \limsup_{n \to \infty} \EEE \int_{\Omega_n} \Psi_n^{\rm cell}(\nabla \hat{v}^{k}_{n}) \, {\rm d}x & \ge  \liminf_{n \to \infty} \int_{\Omega_n} \bar{\chi}_n \Psi_n^{\rm cell}(\nabla \hat{v}^{k}_{n}) \, {\rm d}x +  \BBB \limsup_{n \to \infty} \EEE \int_{\Omega_n} (1-\bar{\chi}_n) c| \nabla \hat{v}^{k}_{n}|^{\MMM 2 \EEE} \, {\rm d}x \\
& \ge  \int_{\Omega} \Phi(\nabla v)\, {\rm d}x +  c\eta.
\end{align*}
 This, however, contradicts \eqref{energy convi} and shows that $(h_n)_n$ is equi-integrable. \BBB As in the previous proof, we define \EEE $\chi_{n}(x) = \chi_{[0,\varepsilon_n^{\MMM -1/4 \EEE})}(|\nabla \hat{v}^{k}_n(x)|)$. Then, the equi-integrability   and Lemma  \ref{lem:help2}(iii)  yield
 \[\lim_{n\to \infty} \int_{\Omega_n} (1-\chi_n(x))   \Psi^{\rm cell}_n(\nabla \hat{v}^{k}_n)   \, {\rm d}x =0, \quad \quad \lim_{n\to \infty} \int_{\Omega_n} (1-\chi_n(x))  \Phi(\nabla \hat{v}^{k}_{n})   \, {\rm d}x =0 \,. \]
Along with Lemma \ref{lem:help2}(i) this shows 
 $$\Big| \int_{\Omega_n} \Psi_n^{\rm cell}(\nabla \hat{v}^{k}_{n}) \, {\rm d}x - \int_{\Omega_n} \Phi(\nabla \hat{v}^{k}_{n})\, {\rm d}x  \Big|  \to 0 $$
 and thus, again by \eqref{energy convi},
 \begin{equation}\label{phi-conv}\lim_{n \to \infty}  \int_{\Omega} \chi_{\Omega_n} \Phi(\nabla \hat{v}^{k}_{n})\, {\rm d}x = \int_{\Omega} \Phi(\nabla v)\, {\rm d}x  \,.\end{equation}
Now, since $\Phi$ is  strictly convex and homogeneous of degree \MMM 2\EEE, see \eqref{eq: lineraition}, there exists a constant $c>0$ such that $y\mapsto \Phi(y)-  c |y|^{\MMM 2 \EEE}$ is convex. Hence, we can estimate by \eqref{phi-conv}
\begin{align*}
\liminf_{n\to \infty}\int_{\Omega_n} - c|\nabla \hat{v}^{k}_{n}|^{\MMM 2 \EEE} \,  {\rm d}x & = \liminf_{n\to \infty} \int_{\Omega_n} \Big( \Phi (\BBB \nabla \EEE \hat{v}^{k}_{n})- c|\nabla \hat{v}^{k}_{n}|^{\MMM 2 \EEE}-  \Phi(\BBB \nabla \EEE\hat{v}^{k}_{n}) \Big) \, {\rm d}x  \\
& \ge \int_{\Omega} \big( \Phi(\nabla v)- c|\nabla v|^{\MMM 2 \EEE} \big) \, {\rm d}x\, - \int_{\Omega}\Phi(\nabla v)\, {\rm d}x\,=-\int_{\Omega} c|\nabla v|^{\MMM 2 \EEE} {\rm d}x\,.
\end{align*}
This implies $\limsup_{n\to \infty} \|\nabla \hat{v}^{k}_{n} \|_{L^{\MMM 2 \EEE}(\Omega_n)} \leq \|\nabla v \|_{L^{\MMM 2 \EEE}(\Omega)} $.  
Since $L^{\MMM 2 \EEE}(\Omega; \R^2)$ is a uniformly convex Banach space, we can apply \cite[Proposition 3.32]{brezis} to infer that \eqref{energy convi2} holds. 
\end{proof}

\subsection{Convergence of sets}\label{set conv-sec}

After we have discussed compactness of displacement fields and corresponding lower semicontinuity results for elastic energies in Subsection \ref{sec: comp}, this section is devoted to  results for crack sets. We start by recalling the notion of $\sigma$-convergence introduced in \cite{GiacPonsi}.

We denote by $\mathcal{A}(\Omega)$ all open subsets of $\Omega$ and define the family of sets with finite perimeter in $\Omega$ by 
\[P(\Omega)\defas \{v\in SBV(\Omega)\colon u(x)\in \{0,1\} \,\text{for a.e.\ } x\in \Omega\}\,.\] For a sequence $(K_n)_{n\in \N}$ of rectifiable sets in $\Omega$, we  consider the functionals $\mathcal{H}^{-}_n \colon P(\Omega)\times \mathcal{A}(\Omega)\to \R$ defined by 
\[\mathcal{H}_{n}^{-}(u, A)=\mathcal{H}^{1}\big((S(u)\setminus K_n) \cap A \big) \,.\]
Following \cite[Theorem 3.2]{AmbrosioBraides}, we assume that for every $A\in \mathcal{A}(\Omega)$, $ \mathcal{H}^{-}_n (\cdot, A)$  
$\Gamma$-converges with respect to the strong topology of $L^1(\Omega)$  to a functional $\mathcal{H}^{-}(\cdot ,A)$, which by   \cite[Theorem 3]{BFLM} \EEE is of the form
\[\mathcal{H}^{-}(u,A)\defas \int_{S(u)\cap A} h^{-}(x,\nu)\, {\rm d}x\]
for some function $h^{-}\colon\Omega\times \mathbb{S}^1 \to [0,\infty)$. 

\begin{definition}[$\sigma$-convergence]\label{def: sigma conv} 
Let $(K_n)_{n\in \N}$ be a sequence of rectifiable
sets in $\Omega$. We say that $K_n$ $\sigma$-converges to $K$ if the functionals $\mathcal{H}_{n}^{-}$ $\Gamma$-converge in the strong topology of $L^1(\Omega)$ to the functional $\mathcal{H}^{-}$, and $K$ is the (unique) rectifiable set $K$ in $\Omega$  such that 
\begin{equation*}
    h^{-}(x, \nu_K(x)) = 0 \quad  \text{for  $\mathcal{H}^{1}$-a.e.\ $x\in K$}, 
    \end{equation*}
and such that for every rectifiable set $H\subset \Omega$ we have
 \[h^{-}(x, \nu_H (x)) = 0
\quad \text{for   $\mathcal{H}^{1}$-a.e.\  \ $x\in H \Rightarrow H\, \tilde{\subset} \,  K$}.\]  
\end{definition}

% \RRR Adjust the following. Maybe we want to keep the main compactness and lsc results \EEE 

Unfortunately, we cannot employ this notion in our setting, but  in view of the definitions in \eqref{eq: hat sets} and \eqref{eq: small jumpi}, we have to deal with two parts of the jump set
\begin{equation}\label{jump-sets-SL}
{K}^{k,{\rm L}}_n(v_n) = \bigcup_{i=1,2,3} \Kl{k,i}{v_n}\quad \text{ and } \quad {K}^{k,{\rm S}}_n(v_n)=\bigcup_{\Cn{k}{v_n}}  \partial \triangle,\end{equation}
where intuitively only ${K}^{k,{\rm L}}_n(v_n)$ should contribute to the crack in the continuum limit whereas ${K}^{k,{\rm S}}_n(v_n)$ should not affect it, see the comment in Remark \ref{good remark}(iii). To solve this issue, we introduce another jump interpolation which only exhibits discontinuities on  ${K}^{k,{\rm L}}_n(v_n) $ and invoke the notion of \emph{$\sigma^p$-convergence}, introduced in \cite{dMasoFranToad}, where $p$ stands for the integrability of the absolutely continuous parts of derivatives. \\ In the following, we say that $v_j \to v$ \emph{weakly in $SBV^p(\Omega)$} if $v_j,v \in SBV^p(\Omega) \cap L^\infty(\Omega)$, $v_j \to v$ a.e.\ in $\Omega$, $\nabla v_j \rightharpoonup \nabla v$ weakly in $L^p(\Omega;\R^2)$, and $\Vert v_j \Vert_\infty + \mathcal{H}^1(S(v_j))$ is uniformly bounded in $j$.

\begin{definition}[$\sigma^p$-convergence]\label{def: sigma-p} 
Let $(K_n)_{n\in \N}$ be a sequence of rectifiable sets in $\Omega$. We say that $K_n$ $\sigma^{p}$-converges to $K \subset \Omega$ in $\Omega$ if $\mathcal{H}^1({K_n})$ is uniformly bounded in $n$ and if the following is satisfied:
\begin{enumerate}[label=(\roman*)]
    \item If a sequence $(v_j)_j \subset SBV^p(\Omega)$ with $v_j\to v$ weakly in $SBV^{p}(\Omega)$ fulfills $S(v_j) \, \tilde{\subset} \,  K_{n_j}$ for some sequence $n_j \to \infty$, then $S(v) \, \tilde{\subset} \,  K$. 
    \item There exists a function $v\in SBV^p(\Omega)$ and a sequence $(v_n)_n$ that converges weakly to $v$ in $SBV^{p}(\Omega)$ such that $S(v_n)\, \tilde{\subset} \,  K_n$ for each $n\in \N$ and $K \, \tilde{=} \,  S(v)$, where $\, \tilde{=} \, $ stands for equality up to an $\mathcal{H}^1$-negligible set.          
\end{enumerate}     
\end{definition}

\begin{remark}\label{sigma-p-remark}
    (a) The $\sigma^{p}$-limit of any sequence is always rectifiable since $S(u)$ is rectifiable for each $u\in SBV(\Omega)$.  \\
    (b) Note that $\sigma$-convergence is invariant under perturbations with vanishing $\mathcal{H}^1$-measure, i.e, given $(\Gamma_n)_n$ and $(\Gamma^*_n)_n$  with $\mathcal{H}^1(\Gamma_n^*) \to 0$, the $\sigma$-limits of $(\Gamma_n)_n$ and $(\Gamma_n \cup \Gamma^*_n)_n$, if existent, coincide, \BBB see \cite[Remark~5.2]{GiacPonsi}. \EEE The same property for $\sigma^p$ convergence does not hold in general, \BBB see \cite[Remark~5.10]{GiacPonsi}. \EEE \\
    (c) If $\Gamma_n$ and $\Gamma'_n$ $\sigma^p$-converge to $\Gamma$ and $\Gamma'$ and $\Gamma_n\subset \Gamma'_n$ for all $n\in \N$, then it follows $\Gamma\subset \Gamma'$, \BBB see \cite[Remark~4.2]{dMasoFranToad}. \EEE \\
    (d)  Note that the $\sigma$-limit of a sequence always contains its $\sigma^p$-limit, see \cite[Corollary 5.9]{GiacPonsi}. The opposite inclusion is false. See \cite[Remark 5.10]{GiacPonsi} for an example, where the $\sigma^p$-limit is strictly contained in the $\sigma$-limit. \EEE 
\end{remark}

In the following, we  need  the  compactness and lower semicontinuity properties of $\sigma^{p}$-convergence proved in  \cite[Theorems 4.3, 4.7, and 4.8]{dMasoFranToad}. 
 
\begin{theorem}\label{helly-sigma-p}
     Let $p \in (1,\infty)$.    Let $t\mapsto \Gamma_n(t)$ be a sequence of increasing set functions defined on an interval $I\subset \R$ with values contained in $\Omega$, i.e., $\Gamma(s)\subset \Gamma(t)\subset \Omega$ for every $s,t\in I$ with $s<t$. Assume that $\mathcal{H}^{1}(\Gamma_n(t))$ is bounded uniformly with respect to $n$ and $t$. Then, there exist a subsequence $(\Gamma_{n_k})_k$ and an increasing set function $t\mapsto \Gamma(t)$ on $I$ such that for every $t\in I$ we have:
    \begin{enumerate}[label=(\alph*)]
        \item $\Gamma_{n_k}(t)\to \Gamma(t)$ in the sense of $\sigma^p$-convergence.
        \item Let $\varphi\colon \R^{2}\to [0,\infty) $ be a norm on $\R^2$. Then, we have 
    \begin{align*}
        \int_{\Gamma(t)} \varphi(\nu_\Gamma) \, {\rm d}\mathcal{H}^1(s) \leq \liminf_{k\to \infty} \int_{\Gamma_{n_k}(t)} \varphi(\nu_{\Gamma_{n_k}}) \, {\rm d}\mathcal{H}^1\,,
    \end{align*}
    where $\nu_\Gamma$ and  $\nu_{\Gamma_{n_k}}$ denote measure-theoretic unit normals to $\Gamma(t)$ and $\Gamma_{n_k}(t)$, respectively.
\end{enumerate}
\end{theorem}

\begin{remark}\label{wurschel remark}
{\normalfont
The same two properties   (a) and (b)  also hold for $\sigma$-convergence in place of $\sigma^p$-convergence, see \cite[Propositions 5.3 and 5.7]{GiacPonsi}.}
\end{remark} 
\noindent
Recall \eqref{jump-sets-SL}. The essential point for us will be the following property. 
 
\begin{proposition}[\BBB$\sigma^{\OOO 3/2 \EEE }$-limit of crack set\EEE]\label{crucial-inclusion}
\EEE Let $(v_n)_n$ be a sequence of atomistic displacements and, for each $n \in \N$, let $(u^j_{n})_{j < k_n}$ be a displacement history for $k_n \in \N$ depending on $n$. Suppose that $\mathcal{E}_{n}(v_n; (u^j_{n})_{j < k_n}) + \|\tilde{v}_n\|_{L^{\infty}(\Omega_n)} \le M$ for some $M>0$. Assume that $K_n^{k_n,{\rm L}}(v_n)$ $\sigma^{\OOO 3/2 \EEE}$-converges to $K$ and that $v_n$  AC-converges to some $v \in SBV^{\MMM  2 \EEE}(\Omega)$. Then, we have $S(v) \, \tilde{\subset} \,  K$.  
% For a $t\in [0,T]$ let $K_n^{\rm L}(t)$ as in \eqref{def: KL} and assume that $ K_n^{\rm L}(t)$ $\sigma^2$-converges to $K$. Suppose that $u_n(t)$ given in \eqref{definitionofun} AC-converges to some $u \in SBV(\Omega)$. Then, we have $S(u) \subset K$.  
\end{proposition} 
\begin{remark}
{\normalfont
(i) The \BBB gist \EEE is that the limit is taken with respect to $K_n^{k_n,{\rm L}}(v_n)$, i.e., only with respect to the part of the crack set which should be relevant for the continuum limit.  Note that the result is not an immediate consequence of (i) in the definition of $\sigma^p$-convergence as $S(\hat{v}^{k_n}_n)$ is \emph{not contained} in $K_n^{k_n,{\rm L}}(v_n)$, cf.\ Remark \ref{good remark}(ii).

(ii)  Note that the choice \JJJ $ p = 3/2$ \EEE is for definiteness only and could be replaced by any $1 < p < 2$. Later we see that along the evolution   the $\sigma$-limit and the $\sigma^{p}$-limit coincide, see  Step 5 \BBB in \EEE the proof of Theorem~\ref{maintheorem}. 
}
\end{remark}

\begin{proof}
% Recall the definition $u_n(t)=u_n^{k(t)}$ for $t\in [t_{n}^{k(t)},t_{n}^{k(t)+1})$. Since $t\in [0,T]$ is fixed, we will drop the dependency of $k$ on $t$ and just denote $k\defas k(t)$. 
For convenience, we drop the dependency of $k$ on $n$ in the proof and simply write $k$ in place of $k_n$.
The idea is  to construct an interpolation $\bar{v}_n\colon\Omega_n\to \R$ of $v_n\colon\mathcal{L}_{n}(\Omega_n)\to \R$ such that 
\begin{enumerate}
    \item $S(\bar{v}_n)\subset K_n^{k,{\rm L}}(v_n)$,
    \item  $\chi_{\Omega_n}\bar{v}_n\to v$ weakly in $SBV^{\OOO 3/2 \EEE}(\Omega)$. 
\end{enumerate}
Then, from  item (i) of the definition of $\sigma^{\OOO 3/2 \EEE}$-convergence we obtain  $S(v) \, \tilde{\subset} \,  K$. 

For triangles $\triangle \in \mathcal{T}_{n}\setminus\Cl{k}{v_n}$ we set $\bar{v}_n\defas \tilde{v}_n$. Recall, that by definition (see \eqref{big jump tri}) we then have
 \begin{equation}\label{grad: outside-jump}
 |\nabla \bar{v}_n|= |(v_n)_{\triangle}|\leq C\varepsilon_n^{\MMM -1/2 \EEE} R_n\quad \text{on triangles $\triangle \in \mathcal{T}_{n}\setminus \Cl{k}{v_n}$}\,.
 \end{equation} 
Now we consider triangles $\triangle \in \Cl{k}{v_n}$ and recall the division into sub-triangles introduced in Figure \ref{different-triangles}.

\BBB First, consider $\triangle \in \mathcal{C}^{k,3}_n(v_n)$. On each subtriangle $\triangle_{\alpha}$, $\alpha = 1,2,3$, we   define  $\bar{v}_n$ as an affine interpolation  such that $\bar{v}_n=v_n$ on the vertex and 
\begin{align}\label{affiness}
\nabla \bar{v}_n \cdot \mathbf{v} = \begin{cases}  \nabla \tilde{v}_n \cdot \mathbf{v}   & \text{ if }    |\nabla \tilde{v}_n \cdot \mathbf{v}|  \le 2R_n \varepsilon_n^{\MMM -1/2 \EEE}, \\    0   & \text{ if }    |\nabla \tilde{v}_n \cdot \mathbf{v}| >2 R_n \varepsilon_n^{\MMM -1/2 \EEE}   \end{cases} \quad \text{ for $\mathbf{v} \in \lbrace \mathbf{v}_\beta , \mathbf{v}_\gamma  \rbrace$ }
\end{align}
where $\alpha,\beta,\gamma$ are pairwise distinct. In $\triangle_{\rm mid}$, we define $\bar{v}_n$ as a constant. 

If $\triangle \in \mathcal{C}^{k,2}_n(v_n)$, assuming $M^{k,\mathbf{v}_\alpha}_{n, \triangle}(v_n)\leq R_n$, we  define  $\bar{v}_n$ on $\triangle_{\alpha}$,   $\triangle_{\beta}$,  $\triangle_{\gamma}$ as in \eqref{affiness},  and in   $\triangle_{\rm mid}$ \EEE we   choose an affine interpolation between the \BBB affine functions attained on \EEE $h^{\beta}_{\triangle}$ and $h^{\gamma}_{\triangle}$. 

Note that by our construction the only discontinuities in \BBB the closure of a \EEE triangle $\triangle$ appear on $\bigcup_{i=1}^{3}\, \BBB h_{\scalebox{1.0}{$\scriptscriptstyle \triangle$}}^{i} \EEE $ for $ \triangle \in \mathcal{C}^{k,3}_n(v_n) $ and on $h_{\scalebox{1.0}{$\scriptscriptstyle \triangle$}}^m$  for $\triangle \in \mathcal{C}^{k,2}_n(v_n)$, where $m$ is the unique index with $M^{k, \mathbf{v}_m}_{n, \triangle}(v_n) \le  R_n$. \BBB In fact, given $(x,x') \in \mathcal{\rm{NN}}_n(\Omega)$ with $x-x'$ parallel to $\mathbf{v}_\alpha$, denote by $\triangle$ and $\triangle'$ the triangles adjacent to the spring.  Then, we have $\bar{v}_n(x) = v_n(x)$, $\bar{v}_n(x') = v_n(x')$, and on the corresponding sub-triangles $\triangle_\beta\cup\triangle_\gamma\cup \triangle'_\beta \cup\triangle'_\gamma$  it either holds $\nabla \bar{v} \cdot \mathbf{v}_\alpha = \nabla \tilde{v}_n \cdot \mathbf{v}_\alpha $ or $\nabla \bar{v} \cdot \mathbf{v}_\alpha  = 0$, where we also use \eqref{big jump tri} to cover the case where one of triangles $\triangle$ and $\triangle'$ lie in $\Cl{k}{v_n}$  and the other one does not.  

  Thus, \EEE  in view of \eqref{eq: hat sets}, we get that $S(\bar{v}_n)\subset K_n^{k,{\rm L}}(v_n)$. Moreover, the construction leads to $|\nabla \bar{v}_n(x)|\leq C  R_n \varepsilon_{n}^{\MMM -1/2 \EEE} $ for $x\in \triangle$ with $ \EEE \triangle \in \BBB \mathcal{C}^{k,{\rm L}}_n(v_n) $. Putting this together with \eqref{grad: outside-jump} yields
 \[|\nabla \bar{v}_n|\leq C R_n \varepsilon_n^{\MMM -1/2 \EEE}\,\quad \text{a.e.\ on}\;\,\Omega_n\,.\]
Using $\| \tilde{v}_n   \|_{L^{\infty}(\Omega_n)} \le C$, by construction we also have  $\| \bar{v}_n   \|_{L^{\infty}(\Omega_n)} \le C$. 
 
Now, recall that by Proposition \ref{bounded-jump} we  have 
\begin{equation}\label{evolution: bounds}
    \# \mathcal{C}^{k}_{n}(v_n)\leq \frac{C}{\varepsilon_n}\quad \text{and} \quad \|\nabla \hat{v}^k_n\|_{L^{\MMM 2 \EEE}(\Omega_n)} ^{\MMM 2 \EEE}\le C. 
\end{equation}
Invoking \eqref{rnepsn} this leads to 
\[\sum_{\triangle \in {\mathcal{C}^{k}_n(v_n)} } \int_{\triangle} |\nabla \bar{v}_n(x)|^{\OOO 3/2 \EEE }\, {\rm d}x \JJJ  \leq C \varepsilon_{n}^2 \, \# \mathcal{C}^{k}_n(v_n)  \, R_n^{3/2}\varepsilon_n^{ -3/4 } \le C R_n^{3/2 }\varepsilon_n^{1/4} \to 0\,. \]
\EEE
As outside of the  triangles   $\mathcal{C}^{k}_n(v_n)$ we have $\bar{v}_n= \hat{v}^k_n$, we can estimate by \eqref{evolution: bounds}  \BBB and H\"older's inequality \EEE 
\[ \limsup_{n \to \infty }\int_{\Omega_n}|\nabla \bar{v}_n(x)|^{\OOO 3/2 \EEE }\, {\rm d}x \leq  \limsup_{n \to \infty }  C\, \Bigl(\int_{\Omega_n}|\nabla \bar{v}_n(x)|^{\OOO 2}\, {\rm d}x\Bigr)^{\OOO 3/4}  \leq  \limsup_{n \to \infty }  \BBB  C\|\nabla \hat{v}^k_n\|_{L^{2}(\Omega_n)}^{\OOO 3/2 \EEE } \EEE \le C\,. \] 
Eventually, \eqref{evolution: bounds} also yields 
\[\mathcal{H}^{1}(S(\bar{v}_n))\leq \mathcal{H}^{1}(K^{k,{\rm L}}_n(v_n))\leq C \# \mathcal{C}^{k,L}_{n}(v_n)\, \varepsilon_n \BBB \leq C \# \mathcal{C}^{k}_{n}(v_n)\, \varepsilon_n \EEE  \leq C\,. \]
Altogether, we have $\|\nabla \bar{v}_n\|_{L^{\OOO 3/2 \EEE }(\Omega_n)} + \mathcal{H}^1(S(\bar{v}_n))+ \|\bar{v}_n\|_{\infty}<\infty $, and thus by Ambrosio's compactness theorem and the assumption that $v_n$  AC-converges to $v$ we conclude   $\chi_{\Omega_n}\bar{v}_n\to  v $ \BBB weakly \EEE in $SBV^{\OOO 3/2 \EEE }(\Omega)$. This proves the assertion.   
\end{proof}     
\section{Convergence of quasi-static crack growth: Proof of Theorem \ref{maintheorem}}\label{se:c main}

This section is devoted to the proof of  Theorem \ref{maintheorem}. In Subsection \ref{sec: 5.1}, we start by considering the time-discrete atomistic problems. Afterwards, in Subsection \ref{sec: 5.2} we establish  compactness results for displacements and crack sets. Eventually, Subsection \ref{sec: 5.3} is subject to the proof of the main result, namely the convergence to an irreversible quasi-static crack evolution $t\mapsto (u(t),K(t))$.

\subsection{Time-discrete atomistic evolution}\label{sec: 5.1}

We recall the definition of the time-discrete evolution in \eqref{definitionofun}. In particular, we denote by  $\{0=t^{0}_n< t^{1}_n<\dots< t^{T /\delta_n}_n =T\}$ a  time discretization of the interval $[0,T]$, and let $(u^k_{n})_{k}$ be the   corresponding  displacements. Recalling the boundary values $(g(t^k_n))_k$, we define $g_n^k \colon \L_n(\Omega) \to \R$ by $g_n^k(x) = g(t^k_n,x)$.

\begin{lemma}\label{lemma: dm}
The minimization problems \eqref{minimizing-scheme0} and \eqref{minimizing-scheme} admit solutions. The solution \BBB $t \mapsto u_n(t)$ \EEE given in \eqref{definitionofun} satisfies $\sup_{t\in [0,T]}\|\tilde{u}_n(t)\|_{L^\infty(\Omega_n)} < +\infty$.
\end{lemma}

\begin{proof}
We observe that we can restrict the minimization problem to functions $v \in \mathcal{A}_n (g(t^{k}_n)) $ which satisfy $|v(x)| \le \max_{x \in \L_n(\Omega)} |g^k_n(x)| \BBB \le C \EEE $ for all $x \in \L_n(\Omega)$ by the monotonicity of the function $\Psi$ \BBB and the regularity of $g$. \EEE As  the space of displacements from $\L_n(\Omega)$ to $\R$  is finite dimensional, this along with the continuity of $\Psi$ guarantees existence by the Direct Method.  
\end{proof}
\noindent
\BBB Similarly as in \eqref{definitionofun}  we write \EEE
\begin{align*}
\hat{u}_n(t):= \hat{u}^{k}_n \; \quad  \text{for}\; t\in[t^{k}_n,t^{k+1}_n)\,,
\end{align*}
%\label{definitionofun-jump}
 where  the interpolation $\hat{u}^k_n$ is defined in Subsection \ref{sec: jump interp} (for general $v_n$ in place of $u_n^k$). Recalling \eqref{basicenergy-neu} and \eqref{energy-triangle}, given an arbitrary $v_n \colon \L_n(\Omega) \to \R$, we also introduce the shorthand notation 
 \begin{align}\label{energy along evo}
\mathcal{E}_n^{0}(v_n) \defas     \mathscr{E}_n(v_n), \quad \quad   \mathcal{E}_n^{k}(v_n) \defas    \mathcal{E}_n(v_n ;(u_n^j)_{j<k}) \quad \text{for $k \ge $1}.
\end{align}

\begin{lemma}[Bound on derivatives]\label{lemma control on ela}
 Let \BBB $t \mapsto u_n(t)$ \EEE be the discrete evolution defined in \eqref{definitionofun}. We have  $\mathcal{E}_n^{0}(u^0_n) \le C$. Moreover,   it holds that 
\begin{align}\label{eq: hat finally}
 \int_{\Omega_n} |\nabla \hat{u}_n(t)|^{\MMM 2 \EEE}  \, {\rm d} x \le C \quad \text{ for all $t \in [0,T]$ and for all $n \in \N$}
 \end{align}
    for some  $C>0$ depending on $W$ and $g$. 
\end{lemma}

\begin{proof}
By \eqref{basicenergy-neu} and by writing the energy as an integral using \eqref{eq: ela en}--\eqref{psitriangle}, we get  
\begin{align}\label{E-g}
\mathscr{E}_n(g^0_n)  = \frac{\eps_n}{2}\, \sum_{(x,x') \in \mathcal{\rm{NN}_{\varepsilon}}(\Omega)} \Psi\Bigl(\frac{ \eps_n^{\MMM 1/2 \EEE} |g^0_n(x)-g^0_n(x')|}{\varepsilon_n}\Bigr) =  \sum_{\triangle \in \mathcal{T}_n}   \int_{\triangle}  \Psi_n^{\rm cell}  (  (g_n^0)_\triangle) + \BBB \mathcal{E}^{\rm bdy}_n(g_n^0), \EEE    
\end{align}
where $\mathcal{E}^{\rm bdy}(g_n^0)$ accounts for springs at the boundary, similarly to \eqref{bdy energy}. \BBB As $\mathcal{E}^{\rm bdy}_n(g_n^0)$  \EEE  is bounded and as $(g_n^0)_\triangle$ is uniformly bounded by $g \in W^{1,1}([0,T]; W^{2,\infty}(\Omega))$, Lemma \ref{lem:help2}\BBB(iii) \EEE yields $\mathscr{E}_n(g^0_n) \le C$ independently of $n$ for some $C>0$ large enough, which depends on $\Phi$ (and thus on $W$) and $g$. In particular, by \eqref{minimizing-scheme0} this implies $\mathscr{E}_n(u^0_n) \le C$.

In a similar fashion, we obtain a control on  $u_n^k$ for $k \ge 1$. As $\eps_n^{\MMM 1/2 \EEE} \Vert \nabla g \Vert_{L^\infty([0,T]; L^\infty(\Omega))} \le R$ for $n$ large enough, \eqref {energy-in-u} implies
\begin{align*}
\mathcal{E}^k_n(g_n^k)=\frac{\varepsilon_n}{2}\, \sum_{M_{x,x'}^{\eps_n}((u_n^j)_{j<k})\leq R}  \Psi\big(\eps_n^{\MMM -1/2 \EEE} |g^k_n(x) - g^k_n(x')  | \big)  +  \frac{\varepsilon_n}{2}\,   \sum_{M_{x,x'}^{\eps_n}((u_n^j)_{j<k})> R} \Psi (M_{x,x'}^{\eps_n}((u_n^j)_{j<k}) ) , 
\end{align*}
where the sums run over all $(x,x')\in \mathcal{\rm{NN}}_n(\Omega)$. Note that by \eqref{minimizing-scheme} we have $\mathcal{E}^k_n(u_n^k) \le \mathcal{E}^k_n(g_n^k)$. \EEE Using the monotonicity of the memory variable in $k$ and the monotonicity of $\Psi$, this implies 
$$\frac{\varepsilon_n}{2}\, \sum_{M_{x,x'}^{\eps_n}((u_n^j)_{j\le k}) \EEE \leq R}  \Psi\big(\eps_n^{\MMM -1/2 \EEE} |u^k_n(x) - u^k_n(x')|  \big)  \le \frac{\varepsilon_n}{2}\, \sum_{ M_{x,x'}^{\eps_n}((u_n^j)_{j\le k}) \EEE\leq R}  \Psi\big(\eps_n^{\MMM -1/2 \EEE} |g^k_n(x) - g^k_n(x') | \big)  \le C  $$
for $C$ depending on $\Phi$ and $g$, where the last step follows by repeating the argument in \eqref{E-g}.  Recalling the representation of the elastic energy in \eqref{eq: ela en}--\eqref{psitriangle}, we thus find $ \int_{\Omega_n} \Psi_n^{\rm cell}(\nabla \hat{u}^k_n)  \, {\rm d} x  \le C$. By \eqref{eq: not so uniform, but neu} and  the property of $\Psi_n^{\rm cell}$ stated in Lemma~\ref{lem:help2}(ii) \BBB we get \EEE the bound 
$ \int_{\Omega_n} |\nabla \hat{u}^k_n|^{\MMM 2 \EEE} \, {\rm d} x \le C$. This shows \eqref{eq: hat finally} and concludes the proof. 
\end{proof}

We continue with an  energy estimate which will be vital for the proof of  \eqref{energybalance} and also delivers a priori bounds for the discrete evolutions. To this end, we   discretize the boundary values $g\in W^{1,1}([0,T];W^{2,\infty}(\Omega))$  in  space. Given $g(\tau)$ for $\tau \in [0,T]$, we define the function  $g_n(\tau) \colon \mathcal{L}_n(\Omega) \to \R$ by $g_n(\tau)(x) = g(\tau, x)$ for $x \in \mathcal{L}_n(\Omega)$. As before, by $\tilde{g}_n(\tau) \colon \Omega_n \to \R$ we denote the interpolation which is piecewise affine on the triangles of $\mathcal{T}_n$. In the same way, we define the time derivative $\partial_t \tilde{g}_n(\tau)$.

\begin{lemma}\label{e-balance-lemma1}
 Let \BBB $t \mapsto u_n(t)$ \EEE be the discrete evolution defined in \eqref{definitionofun}. \BBB Let \EEE $k=0,\ldots,(T /\delta_n) - 1$. Then, there exists $(\eta_n)_n$ \BBB independent of $k$ \EEE with $\eta_n \to 0$ such that
    \begin{equation}\label{eq: for en bal}
        \mathcal{E}_n^k( \BBB u_n^k \EEE ) - \mathcal{E}_{n}^{0}(u^0_n) \leq \int_{0}^{t_{n}^{k}} \int_{\Omega_n} D \Psi_n^{\rm cell }(\nabla \hat{u}_n(\tau)) \cdot \nabla \partial_{t}\tilde{g}_n(\tau)\,{\rm d} x\,{\rm d}\tau  +  ( 1  + \eps_n \# \mathcal{C}^k_n)   \eta_n, 
    \end{equation}
where $\mathcal{C}^k_n$ is defined in \eqref{nom jump tri}. 
\end{lemma}

\begin{proof}
The proof follows a rather standard strategy, see e.g.\ \cite[Lemma 6.1]{dMasoFranToad}. Here, we include full details and in particular describe how the estimates are adapted to our discrete setting. For each step $l$, where $0<l\leq k$, we consider the test function $ \BBB \xi_{n}^l \EEE :=u_{n}^{l-1}+g_{n}^{l}-g_{n}^{l-1}$ and obtain by \eqref{minimizing-scheme}   
    \begin{equation}\label{discrete-global-stability}
\mathcal{E}^l_{n}(u^l_{n})     \leq \mathcal{E}^l_{n}(\xi_{n}^l)\,=\mathcal{E}^l_{n}(u_{n}^{l-1}+g_{n}^{l}-g_{n}^{l-1}).
\end{equation}
    We aim at estimating the right-hand side of this inequality. Our goal is to prove that there exists a bounded sequence  $(\vartheta_n)_n$ in $L^\infty([0,T] \times \Omega )$ with $ \Vert \vartheta_n(\tau) \Vert_{L^{\MMM 2 \EEE}(\Omega)} \to 0 $  uniformly in $\tau$ such that, for each $l$, we have 
    \begin{align}\label{wanttoprove}
          \mathcal{E}^l_{n}(u_{n}^{l-1}+g_{n}^{l}-g_{n}^{l-1})- \mathcal{E}^l_{n}(u_{n}^{l-1})  & \leq \int_{t_{n}^{l-1}}^{t_{n}^{l}}  \int_{\Omega_{n}}   D \Psi^{\rm cell}_n  \big(\nabla \hat{u}_{n}^{l-1}+  \vartheta_n(\tau) \big)   \cdot \nabla \partial_{t} \tilde{g}_n(\tau)\,{\rm d}x \, {\rm d} \tau \notag   \\ & \ \ \ +C( 1  + \eps_n \# \mathcal{C}^l_n)  \, \eps_n^{\MMM 1/2 \EEE} \int_{t_{n}^{l-1}}^{t_n^l} \| \nabla  \partial _t  \tilde{g}_{n}(\tau)  \|_{L^{\infty}(\Omega_n)}  \, {\rm d}\tau.
    \end{align}
We defer the proof of \eqref{wanttoprove}  to Step 2 below and first show how to derive the estimate   \eqref{eq: for en bal}  from this.

\emph{Step 1: Proof of \eqref{eq: for en bal}.} Assume \eqref{wanttoprove} holds. Since $\mathcal{E}^{l-1}_{n}(u_{n}^{l-1}) = \mathcal{E}^l_{n}(u_{n}^{l-1})$ for each step $l$, summing up over all time steps $0 \le l \le k$ and using \eqref{discrete-global-stability} yields a telescopic sum on the left-hand side \BBB of \eqref{wanttoprove} \EEE such that we get 
\begin{equation}\label{eq: a number}
    \begin{aligned}
        \mathcal{E}^k_{n}(u_n^k)- \mathcal{E}_{n}^{0}(u_n^0)
 \leq  \int_{0}^{t_{n}^{k}} \int_{\Omega_{n}}   D \Psi^{\rm cell}_n  \big(\nabla \hat{u}_{n}(\tau)+  \vartheta_n(\tau) \big)   \cdot \nabla \partial_{t} \tilde{g}_n(\tau)\,{\rm d}x \, {\rm d} \tau + C( 1  + \eps_n \# \mathcal{C}^k_n) \, \eps_n^{\MMM 1/2 \EEE} 
    \end{aligned}
\end{equation}
 \BBB for a  constant $C$ depending \EEE on $g \in W^{1,1}([0,T]; W^{2,\infty}(\Omega))$, where we used that the sets $(\mathcal{C}^l_n)_l$ are increasing in $l$. By \MMM $|D\Psi^{\rm cell}_n(z)| \le  C|z| $ \EEE for $z \in \R^2$, see Lemma \ref{lem:help2}(iii), \EEE and  the fact that  $(D\Psi^{\rm cell}_n)_n$ satisfies the continuity property stated in Lemma \ref{lem:help2}(iv), we can invoke \cite[Lemma 2.4]{GiacPonsi} to find  
$$ \Big|  \int_{\Omega_{n}}   D \Psi^{\rm cell}_n  \big(\nabla \hat{u}_{n}(\tau)+  \vartheta_n(\tau) \big)   \cdot \nabla \partial_{t} g(\tau)\,{\rm d}x -   \int_{\Omega_{n}}   D \Psi^{\rm cell}_n  \big(\nabla \hat{u}_{n}(\tau)  \big)   \cdot \nabla \partial_{t} g(\tau)\,{\rm d}x \Big| \to 0 $$
for a.e.\ $\tau\in [0,T]$. \BBB Here, we also employed $\|\vartheta_n(\tau)\|_{L^{\MMM 2 \EEE}(\Omega)}\to 0$ uniformly in $\tau$. \EEE  Again using $|D\Psi^{\rm cell}_n(z)| \le  C|z| $  we get  
\begin{align}\label{eq: nochmal}
\Vert  D \Psi^{\rm cell}_n  \big(\nabla \hat{u}_{n}(\tau) \big) \Vert_{L^{\MMM 2 \EEE}(\Omega_n)} \le C, \quad \quad \Vert  D \Psi^{\rm cell}_n  \big(\nabla \hat{u}_{n}(\tau) +  \vartheta_n(\tau) \big) \Vert_{L^{\MMM 2 \EEE}(\Omega_n)} \le C \quad \BBB \text{for a.e.\ $\tau \in [0,T]$}\,, \EEE
\end{align}
where we used \eqref{eq: hat finally} and for the second estimate also that $\|\vartheta_n(\tau)\|_{L^{\MMM 2 \EEE}(\Omega)}\to 0$ uniformly in $\tau$. This along with the fact that $\chi_{\Omega_n}\nabla \partial_t \tilde{g}_n$ converges strongly to $ \nabla \partial_t  g$ in $L^1([0,T];L^{\infty}(\Omega;\R^2))$ yields for a.e.\ $\tau\in [0,T]$ 
\[ h_n(\tau):= \Big|  \int_{\Omega_{n}}   D \Psi^{\rm cell}_n  \big(\nabla \hat{u}_{n}(\tau)+  \vartheta_n(\tau) \big)   \cdot \nabla \partial_{t} \tilde{g}_n(\tau)\,{\rm d}x -   \int_{\Omega_{n}}   D \Psi^{\rm cell}_n  \big(\nabla \hat{u}_{n}(\tau)  \big)   \cdot \nabla \partial_{t} \tilde{g}_n(\tau)\,{\rm d}x \Big| \to 0 \,.\]
By \eqref{eq: nochmal} we then also get  by Hölder's inequality  
\begin{align*}
|h_n(\tau)| \le     C  \|\nabla \partial_t \tilde{g}_n(\tau)\|_{L^{\infty}(\Omega_n)}.
\end{align*}
\BBB Using \EEE again the convergence of   $\chi_{\Omega_n}\nabla \partial_t \tilde{g}_n$  to $ \nabla \partial_t  g$, we obtain   $\int_0^T h_n(\tau)\, {\rm d}\tau \to 0$ as $n \to \infty$ by   dominated convergence. This along with \eqref{eq: a number}   and setting $\eta_n:= C\varepsilon_n^{\MMM 1/2 \EEE} \BBB + \int_0^T h_n(\tau)\, {\rm d}\tau \EEE$   shows the assertion. 
% \RRR Not directly. One has to add an dominated convergence argument exactly like in \cite{dMasoFranToad} equation (6.34)
% \EEE

\emph{Step 2: Proof of \eqref{wanttoprove}.}  We recall the notation of $\mathcal{C}_n^l$ and  $M_{n,\triangle}^{l, \mathbf{v}} $ for all $\triangle \in \mathcal{T}_n$ and $\mathbf{v} \in \mathcal{V}$, see \eqref{eq: noacons} and  \eqref{nom jump tri}. 
We can write the memory variable of the test function in terms of the piecewise affine interpolation for each $\triangle\in \mathcal{C}_{n}^{l}(\xi_{n}^l)$ and $\mathbf{v} \in \mathcal{V}$, namely
 \begin{align}\label{more detials}
 M_{n, \triangle}^{l, \mathbf{v}}(\xi_{n}^l)=\OOO \varepsilon_n^{1/2} \EEE \sup_{m<l}|\nabla \tilde{u}_{n}^{m} \cdot \mathbf{v}|\vee \OOO \varepsilon_n^{1/2} \EEE |\nabla \BBB \tilde{\xi}_{n}^l \EEE \cdot \mathbf{v}|\,.
 \end{align}
  We  estimate the last term by
    \begin{align}\label{compare-memories0}
     \varepsilon_n^{\MMM 1/2 \EEE} |\nabla \tilde{\xi}_{n}^l \cdot \mathbf{v}|&=\varepsilon_n^{\MMM 1/2 \EEE} |\nabla (\tilde{u}_{n}^{l-1} + \tilde{g}_{n}^{l}- \tilde{g}_{n}^{l-1}) \cdot \mathbf{v}|\leq \varepsilon_n^{\MMM 1/2 \EEE} |\nabla \tilde{u}_{n}^{l-1} \cdot \mathbf{v}| + \varepsilon_n^{\MMM 1/2 \EEE} \|\nabla \tilde{g}_{n}^{l}- \nabla \tilde{g}_{n}^{l-1}\|_{L^{\infty}(\Omega_n)} \notag    \\
& \le  \varepsilon_n^{\MMM 1/2 \EEE} |\nabla \tilde{u}_{n}^{l-1} \cdot \mathbf{v}| + \varepsilon_n^{\MMM 1/2 \EEE}    \int_{t_{n}^{l-1}}^{t_n^l} \| \partial _t \nabla \tilde{g}_{n}(\tau)  \|_{L^{\infty}(\Omega_n)}  \, {\rm d}\tau \,,
      \end{align}
where  $\tilde{g}_{n}^{l}$ denote the piecewise affine interpolation of the atomistic boundary data ${g}_{n}^{l}$. Hence, we obtain  
\begin{equation}\label{compare-memories}
M_{n,\triangle}^{l, \mathbf{v}}(\xi_{n}^l) = M_{n,\triangle}^{l, \mathbf{v}} (u_{n}^{l-1}+g_{n}^{l}-g_{n}^{l-1})\leq M_{n,\triangle}^{l, \mathbf{v}}(u_n^{l-1}) + \varepsilon_n^{\MMM 1/2 \EEE} \int_{t_{n}^{l-1}}^{t_n^l} \| \partial _t \nabla \tilde{g}_{n}(\tau)  \|_{L^{\infty}(\Omega_n)}  \, {\rm d}\tau\,.
\end{equation}
 Recalling \eqref{energy-triangle}, we can write 
\begin{align}\label{eq: nother spli}
 \mathcal{E}^l_{n}(\xi_{n}^l) -  \mathcal{E}^l_{n}(u^{l-1}_{n})    & =  \frac{\varepsilon_n}{2}\sum_{\T_{n}\setminus \mathcal{C}^l_n} \sum_{\mathbf{v}\in\mathcal{V}} \Big( \Psi(\varepsilon_n^{\MMM 1/2 \EEE}|(\xi_{n}^l)_{\triangle}  \cdot \mathbf{v}|)  -  \Psi(\varepsilon_n^{\MMM 1/2 \EEE}|({u}^{l-1}_n)_{\triangle}  \cdot \mathbf{v}|) \Big) + \mathcal{G}^l_n(\xi_{n}^l) \notag \\ & \  \ \    + \mathcal{E}_n^{l, {\rm{bdy}}}(\xi_{n}^l) - \mathcal{E}_n^{l, {\rm{bdy}}}(u^{l-1}_{n}), 
   \end{align}
   where the boundary energy is defined in \eqref{bdy energy} and where we have set 
    $$ \mathcal{G}^l_n(\xi_{n}^l) :=  \frac{\varepsilon_n}{2}\sum_{\mathcal{C}_{n}^{l}} \Big( \sum_{M_{n,\triangle}^{l, \mathbf{v}}>R} \big(\Psi(M_{n,\triangle}^{l, \mathbf{v}} (\xi_{n}^l))   -\Psi( M_{n,\triangle}^{l, \mathbf{v}}   )  \big)  +  \sum_{M_{n,\triangle}^{l, \mathbf{v}} \le R} \big(\Psi(\varepsilon_n^{\MMM 1/2 \EEE}\,|(\xi_{n}^l)_{\triangle}  \cdot \mathbf{v}|) - \Psi(\varepsilon_n^{\MMM 1/2 \EEE}\,|({u}^{l-1}_n)_{\triangle}  \cdot \mathbf{v}|)\big)\Big).$$
Note that here we use that $M_{n,\triangle}^{l, \mathbf{v}}(\xi_{n}^l) = \varepsilon_n^{\MMM 1/2 \EEE}| (\xi_{n}^l)_{\triangle}  \cdot \mathbf{v}|$ whenever $M_{n,\triangle}^{l, \mathbf{v}}(\xi_{n}^l) >R$ and $M_{n,\triangle}^{l, \mathbf{v}} \le R$, cf.\ \eqref{more detials}.  In particular, the contributions of triangles $\triangle \in \mathcal{C}_{n}^{l}(\xi_{n}^l)\setminus \mathcal{C}_{n}^{l} $ are captured \BBB in  \EEE the first term of \eqref{eq: nother spli}. 

\BBB We start by estimating \EEE the term $\mathcal{G}^l_n(\xi_{n}^l)$. In view of \eqref{compare-memories0}--\eqref{compare-memories}, \eqref{eq: not so uniform}, and the monotonicity and   Lipschitz continuity of  $\Psi$, we deduce for each $\triangle \in \mathcal{C}^l_n$ that 
\begin{align*}
\Psi(\varepsilon_n^{\MMM 1/2 \EEE}\,|(\xi_{n}^l)_{\triangle}  \cdot \mathbf{v}|) - \Psi(\varepsilon_n^{\MMM 1/2 \EEE}\,|({u}^{l-1}_n)_{\triangle}  \cdot \mathbf{v}|)& \le L\varepsilon_n^{\MMM 1/2 \EEE}\int_{t_{n}^{l-1}}^{t_n^l} \| \partial _t \nabla \tilde{g}_{n}(\tau)  \|_{L^{\infty}(\Omega_n)}  \, {\rm d}\tau , \\
 \Psi( M_{n,\triangle}^{l, \mathbf{v}} (\xi_{n}^l))-\Psi(M_{n,\triangle}^{l, \mathbf{v}}) & \leq L\varepsilon_n^{\MMM 1/2 \EEE}\int_{t_{n}^{l-1}}^{t_n^l} \| \partial _t \nabla \tilde{g}_{n}(\tau)  \|_{L^{\infty}(\Omega_n)}  \, {\rm d}\tau \,,
 \end{align*}
     where $L$ denotes the Lipschitz constant \BBB of $\Psi$. \EEE Here we distinguished the directions $\mathbf{v} \in \mathcal{V}$ with $M_{n,\triangle}^{l, \mathbf{v}}\leq R$ and $M_{n,\triangle}^{l, \mathbf{v}}>R$, respectively.  %exceeding the threshold $R$ or not. 
     This induces 
    \begin{equation}\label{estimate-crackpart}
    \mathcal{G}^l_n(\xi_{n}^l) \le C\eps_n \# \mathcal{C}^l_n \varepsilon_n^{\MMM 1/2 \EEE} \int_{t_{n}^{l-1}}^{t_n^l} \| \partial _t \nabla \tilde{g}_{n}(\tau)  \|_{L^{\infty}(\Omega_n)}  \, {\rm d}\tau.    
     \end{equation}     
 By a similar argument, noting that $\# \lbrace \triangle \in \mathcal{T}_n \colon \partial \triangle \cap \partial \Omega_n \neq \emptyset  \rbrace \le C\eps_n^{-1}$, we get    
        \begin{equation}\label{estimate-crackpart2}
\mathcal{E}_n^{l, {\rm{bdy}}}(\xi_{n}^l) - \mathcal{E}_n^{l, {\rm{bdy}}}(u^{l-1}_{n}) \le C  \varepsilon_n^{\MMM 1/2 \EEE} \int_{t_{n}^{l-1}}^{t_n^l} \| \partial _t \nabla \tilde{g}_{n}(\tau)  \|_{L^{\infty}(\Omega_n)}  \, {\rm d}\tau.    
     \end{equation}  
 We now address the elastic energy. To this end,  recall the definition of the jump interpolation $\hat{u}_{n}^{l-1}$, which outside of triangles $\triangle \in \mathcal{C}_n^l  $ coincides with the piecewise affine interpolation $\tilde{u}_{n}^{l-1}$. We define a  further  jump interpolation of ${\xi}_{n}$ by
    \[ \BBB \check{\xi}_{n}^l \EEE (x)=\begin{cases}
        \tilde{u}_{n}^{l-1}+(\tilde{g}_{n}^{l}-\tilde{g}_{n}^{l-1})& \text{ if }  x\in \triangle \in \mathcal{T}_{n}\setminus \mathcal{C}_n^l\\
        \hat{u}_{n}^{l-1} ,
        &  \text{ if }  x\in \triangle \in \mathcal{C}_n^l\,.
    \end{cases}\] 
 Note carefully that for technical reasons this definition slightly differs from the interpolation $\hat{\xi}_n$ introduced in Subsection \ref{sec: jump interp},  because on triangles $\triangle \in C_{n}^{l}(\xi_{n}^l)\setminus \mathcal{C}_{n}^{l}$ we use the piecewise affine interpolation $\tilde{u}_{n}^{l-1}$ instead of introducing a jump, and on $\mathcal{C}_{n}^{l}$ we simply use the interpolation of $ \hat{u}_{n}^{l-1}$. Recalling the representation of the elastic energy in \eqref{eq: ela en}, in particular using the fact that $\nabla \check{\xi}_{n}^l =  \nabla  \hat{u}_{n}^{l-1} = 0$ on triangles in $\mathcal{C}_{n}^{l}$, we obtain 
 \[\frac{\varepsilon_n}{2}\sum_{\T_{n}\setminus \mathcal{C}^l_n} \sum_{\mathbf{v}\in\mathcal{V}} \Big( \Psi(\varepsilon_n^{\MMM 1/2 \EEE}|(\xi_{n}^l)_{\triangle}  \cdot \mathbf{v}|)  -  \Psi(\varepsilon_n^{\MMM 1/2 \EEE}|({u}^{l-1}_n)_{\triangle}  \cdot \mathbf{v}|) \Big)= \int_{\Omega_{n}} \big( \Psi_{n}^{\rm cell}( \nabla \check{\xi}^l_n(x)) -  \Psi^{\rm cell}_n(  \nabla \hat{u}^{l-1}_n(x)) \big) \,{\rm d}x\,.\]  
 By the estimates \eqref{estimate-crackpart}--\eqref{estimate-crackpart2}, in view of \eqref{eq: nother spli}, we hence can write 
    \begin{equation}\label{elastic-brings-dissipation000}
     \mathcal{E}^k_{n}(\xi_{n}^l) -  \mathcal{E}^k_{n}(u^{l-1}_{n})    \le     \int_{\Omega_{n}} \hspace{-0.2cm} \big( \Psi_{n}^{\rm cell}( \nabla \check{\xi}^l_n(x)) -  \Psi^{\rm cell}_n(  \nabla \hat{u}^{l-1}_n(x)) \big) \,{\rm d}x \,+ \, C( 1  + \eps_n \# \mathcal{C}^l_n) \,  \varepsilon_n^{\MMM 1/2 \EEE} \int_{t_{n}^{l-1}}^{t_n^l} \| \nabla \partial _t \tilde{g}_{n}(\tau)  \|_{L^{\infty}(\Omega_n)}  \, {\rm d}\tau.
\end{equation}
Now, we want to prove the existence of a function $\vartheta_n\in L^\infty([0,T] \times \Omega)$ such that $\|\vartheta_n(s)\|_{L^{\MMM 2 \EEE }\BBB (\Omega) \EEE }\to 0$ uniformly for $s \in [0,T]$, and 
    \begin{equation}\label{elastic-brings-dissipation}
        \int_{\Omega_{n}} \big( \Psi^{\rm cell}_n(\nabla \check{\xi}^l_{n})- \Psi^{\rm cell}_n(\nabla \hat{u}_{n}^{l-1}) \big) \, {\rm d}x\BBB = \EEE \int_{t_{n}^{l-1}}^{t_{n}^{l}} \int_{\Omega} D\Psi^{\rm cell}_n(\nabla \hat{u}_n(\tau)+\vartheta_n(\tau)) \cdot \nabla \partial_{t}\tilde{g}_n(\tau)\,{\rm d} x\,{\rm d}\tau\,.
    \end{equation}
\BBB Once this is achieved, \eqref{elastic-brings-dissipation000}--\eqref{elastic-brings-dissipation} imply \eqref{wanttoprove}. 

Let us now show \eqref{elastic-brings-dissipation}. \EEE  By the mean value theorem there exists a $\rho^{l-1}_{n}\in [0,1]$ such that
    \begin{equation*}
    \begin{aligned}
    \int_{\Omega_{n}} \big( \Psi^{\rm cell}_n(\nabla \check{\xi}^l_{n})- \Psi^{\rm cell}_n(\nabla \hat{u}_{n}^{l-1}) \big)\, {\rm d}x &    =   \int_{\Omega_{n}}  \big(\Psi^{\rm cell}_n (\nabla \hat{u}_{n}^{l-1}+\nabla \tilde{g}_{n}^{l}-\nabla \tilde{g}_{n}^{l-1})- \Psi^{\rm cell}_n (\nabla \hat{u}^{l-1}_{n}) \big) \, {\rm d}x \\ &=  \int_{\Omega_{n}} D \Psi^{\rm cell}_n  \big(\nabla \hat{u}_{n}^{l-1}+ \rho^{l-1}_{n}(\nabla \tilde{g}_{n}^{l}-\nabla \tilde{g}_{n}^{l-1})\big)\cdot \big(\nabla \tilde{g}_{n}^{l}-\nabla \tilde{g}^{l-1}_{n} \big)\, {\rm d}x \,. 
    \end{aligned}
    \end{equation*}
    In the first   equality we used the fact that, by definition, $\check{\xi}^l_{n}= \hat{u}_{n}^{l-1}+ \tilde{g}_{n}^{l}- \tilde{g}_{n}^{l-1}$, whenever $\check{\xi}^l_{n}$ and $\hat{u}_{n}^{l-1}$ differ. We define the piecewise constant function $\vartheta_{n}: [0,T]\to L^{\MMM 2 \EEE }(\Omega)$  for $s\in [t_{n}^{l-1},t_{n}^{l})$ by
  \[\vartheta_{n}(s):=\rho^{l-1}_{n}(\nabla \tilde{g}_n^{l}-\nabla \tilde{g}_n^{l-1})=\rho^{l-1}_{n} \int_{t_{n}^{l-1}}^{t_{n}^{l}} \nabla \partial_{t} \tilde{g}_{n}(\tau)\, {\rm d} \tau \quad  \text{ on $\Omega_n$}\,, \]
and $0$ outside of $\Omega_n$. Clearly, $(\vartheta_{n})_n$ is bounded in $L^\infty([0,T] \times \Omega)$ since the function $\tau\mapsto \nabla \partial_{t} \tilde{g}_{n}(\tau)$ belongs to $L^1([0,T];L^{\infty}(\Omega_n))$. Moreover, we can use the absolute continuity of the integral and $|t_{n}^{l}-t_{n}^{l-1}| =\delta_n \to 0$ to conclude that $\|\vartheta_n(s)\|_{L^{\MMM 2\EEE}(\Omega)}\to 0$ uniformly in $s$. 
    \BBB Due \EEE to the fact that $u_n$ and $\vartheta_n$ are constant in time on the interval $[t_{n}^{l-1},  {t_{n}^{l}} )$,  we then obtain \eqref{elastic-brings-dissipation} by
    \begin{align*}
   \int_{\Omega_{n}} \big( \Psi^{\rm cell}_n(\nabla \check{\xi}^l_{n})- \Psi^{\rm cell}_n(\nabla \hat{u}_{n}^{l-1}) \big)\, {\rm d}x & \BBB = \EEE
\int_{\Omega_{n}} D \Psi^{\rm cell}_n  \big(\nabla \hat{u}_{n}^{l-1}+  \vartheta_n(\tau) \big)  \cdot \Bigl(\int_{t_{n}^{l-1}}^{t_{n}^{l}}\nabla \partial_{t} \tilde{g}_n(\tau)\,{\rm d} \tau\Bigr) {\rm d}x \, \notag \\ & =\int_{t_{n}^{l-1}}^{t_{n}^{l}}  \int_{\Omega_{n}}   D \Psi^{\rm cell}_n  \big(\nabla \hat{u}_{n}^{l-1}+  \vartheta_n(\tau) \big)   \cdot \nabla \partial_{t} \tilde{g}_n(\tau)\,{\rm d}x \, {\rm d} \tau. 
\end{align*}
\BBB This shows \eqref{elastic-brings-dissipation} and \EEE concludes the proof. 
\end{proof}

%In this section we would like to verify the energy balance \eqref{energybalance}. Note that in the following the deformation history $(u^{m}_{\varepsilon,\delta})_{m<k}$ will remain fixed. For the convenience of the reader we hence drop the history-dependence of the alternative interpolation and use for a discrete deformation $\varphi_{\varepsilon}:\mathcal{L}_{\varepsilon}(\Omega)\to \R$ the notation $\hat{\varphi}_{\varepsilon}:=\hat{\varphi}^{k-1}_{\varepsilon}$. In particular, we mean by $\hat{u}^{m}_{\varepsilon}:=\hat{u}^{m,k-1}_{\varepsilon}$ the interpolation of $u^{m}_{\varepsilon}$ with respect to $C_{\varepsilon}^{k-1}$ . \RRR Note carefully, that this leads to a subtle change of notation, since before we have used the superscript to denote the corresponding time step to which respect we form the alternative interpolation! \EEE 

As a direct consequence, we obtain the following bound on the energy. 

\begin{corollary}[Energy bound]\label{cor: energy bound}
 Let \BBB $t \mapsto u_n(t)$ \EEE be the discrete evolution defined in \eqref{definitionofun}.
Then, there exists a constant $C>0$ only depending on $W$ and $g$ such that 
$${\mathcal{E}^k_n( \BBB u^k_n \EEE ) \le C \quad \text{for all $k = 0, \ldots, T/\delta_n \BBB -1 \EEE $ and $n \in \N$.} }$$
\end{corollary}

\begin{proof}
We show that there is a constant $\hat{C}$ depending on $W$ and $g$, but independent of $k$ and $n$, such that 
\begin{align}\label{vio}
\int_{0}^{t_{n}^{k}} \int_{\Omega_n} D \Psi_n^{\rm cell}(\nabla \hat{u}_n(\tau)) \cdot \nabla \partial_{t}\tilde{g}_n(\tau)\,{\rm d} x\,{\rm d}\tau \le   \hat{C}.
\end{align} 
\BBB Once this is achieved, we conclude as follows. In \EEE view of \eqref{energy-triangle} and the fact that $\Psi(R)>0$, we can find a constant $\tilde{C}$ such that $  \varepsilon_n \#  \mathcal{C}_{n}^{k} \leq  \tilde{C}\, \mathcal{E}^k_n(u_n(t))$.   Therefore, we obtain from Lemma \ref{e-balance-lemma1} and $\mathcal{E}^0_n(u^0_n) \le C$ (see Lemma \ref{lemma control on ela})  
 \[ (1- \tilde{C}\eta_n)\mathcal{E}^k_n( \BBB u_n^k \EEE )\leq \hat{C}+  C + \eta_n \leq C,\]
which implies the desired bound.   

To see \eqref{vio}, we argue by contradiction. If \eqref{vio} is violated, by H\"older's inequality we find for each $\hat{C}>0$ some $n\in \N$ large enough  (depending on 
$\hat{C}$) such that
$$ {\rm ess \, sup}_{\tau \in (0,t_n^k)} \int_{\Omega_n} |D \Psi_n^{\rm cell}(\nabla \hat{u}_n(\tau))|  \,{\rm d} x >\hat{C} \Vert \partial_t \nabla g\Vert_{L^1([0,T]; L^\infty(\Omega))}^{-1}  \,.$$
Using Lemma \ref{lem:help2}(iii) \BBB we \EEE find $c>0$ such that 
$${\rm ess \, sup}_{\tau \in (0,t_n^k)}  \int_{\Omega_n} |\nabla \hat{u}_n(\tau)|  \, {\rm d}x   \ge  c\hat{C}  \Vert \partial_t \nabla g\Vert_{L^1([0,T]; L^\infty(\Omega))}^{-1}. $$
For $\hat{C}$ large enough, using H\"older's inequality this contradicts the bound in  \eqref{eq: hat finally}. This shows \eqref{vio} and concludes the proof.  
\end{proof}

\subsection{Compactness for crack sets and displacement}\label{sec: 5.2}

Based on the energy bound in Corollary \ref{cor: energy bound}, we can pass to the limit in the crack sets and displacements by compactness arguments.  To this end, it will be  convenient to express some of the quantities considered in Section \ref{sec: prel} in terms of the time $t$ in place of the iteration step.  As before,  let $\{0=t^{0}_n< t^{1}_n<\dots<t^{T/\delta_n}_n=T\}$ be a  time discretization of the interval $[0,T]$, and let $(u^j_{n})_{j<k}$ be a  corresponding  displacement history. Recalling \eqref{energy-triangle} and \eqref{energy along evo}, for $t \in [t^k_n,t_n^{k+1})$ we define
\begin{align}\label{energy along evo2}
   \mathcal{E}_n(v_n;t)  := \mathcal{E}_n^{k}(v_n) =     \mathcal{E}_n(v_n ;(u_n^j)_{j<k})
\end{align}
for each $v_n \colon \L_n(\Omega) \to \R$. We use similar notation for the parts of the energy introduced in \eqref{energy:main}. Given the evolution \BBB $t \mapsto u_n(t)$ \EEE defined in \eqref{definitionofun} and recalling  \eqref{jump-sets-SL},  for $t \in [t^k_n,t_n^{k+1})$ we define the jump sets
\begin{equation}\label{def: sets-evol}
    {K}_n^{i}(t) \defas {K}_n^{k,i}(u_n(t)), \quad \quad K^{{\rm L}}_n(t)\defas \bigcup_{i=1}^{3} {K}_n^{i}(t),    \quad  \quad  {K}_n^{\rm S}(t) \defas  {K}_n^{k,{\rm S}}(u_n(t)). 
\end{equation}
One  drawback of the definition of the sets in \eqref{def: sets-evol} is the ambiguity that comes from the three different variants $i=1,2,3$. We will now prove that this ambiguity does not affect the limit, i.e., for a suitable subsequence, the $\sigma^{\OOO 3/2 \EEE}$-limits of all three variants and their union coincide.

\begin{proposition}\label{sigmaconvergence}
    Let \BBB $t \mapsto u_n(t)$ \EEE be the evolution defined in \eqref{definitionofun}. There exists an increasing set function $t\to K(t)$ on $[0,T]$ and a subsequence (not relabeled) such that for every $t\in [0,T]$  %$K_n(t)$ $\sigma$-converges to $K(t)$ %
    we have \[K^{i}_n(t) \quad  \text{$\sigma^{\OOO 3/2 \EEE}$-converges to $K(t)$ in $\Omega$ \quad \text{for all $i=1,2,3$}\,.}\]
    \[K^{{\rm L}}_n(t)  \quad  \text{$\sigma^{\OOO 3/2 \EEE}$-converges to $K(t)$ in $\Omega$\,.}\]
\end{proposition}
\begin{figure}
    \hspace{2.0cm}
\begin{tikzpicture}
    \draw [gray, thick] (0:0)  -- node[anchor=south east] {$\mathbf{v}_1$}++(60:4) --node[anchor=south west] {$\mathbf{v}_2$}++(-60:4) --node[anchor=north ] {$\mathbf{v}_3$} cycle;
    \draw[dashed,red](0:2)--node[anchor=  east] {$h_2$}++(120:2) ;
    \draw[dashed,red](0:2)--node[anchor=  west] {$h_1$}++(60:2) ;
    \draw[dashed,red](60:2)--node[anchor= south] {$h_3$}++(0:2)  ;
    --cycle; \draw [thick,-stealth](4,2) -- (6,2);
\draw [gray, thick] (6,0) --node[anchor=south east] {$\mathbf{v}_1$} ++(60:4) --node[anchor=south west] {$\mathbf{v}_2$}++(-60:4) --node[anchor=north ] {$\mathbf{v}_3$}cycle;     
\draw[dashed,red] (6,0)++(0:2)--node[anchor=  west] {$h_2$}++(120:2);
\draw[dashed,red](6,0)++(0:2)++(60:2)--node[anchor= south] {$h_3$}++(180:2);
\end{tikzpicture}
\caption{The modification in the proof of Proposition \ref{sigmaconvergence}.} \label{fig: extension}
\end{figure}
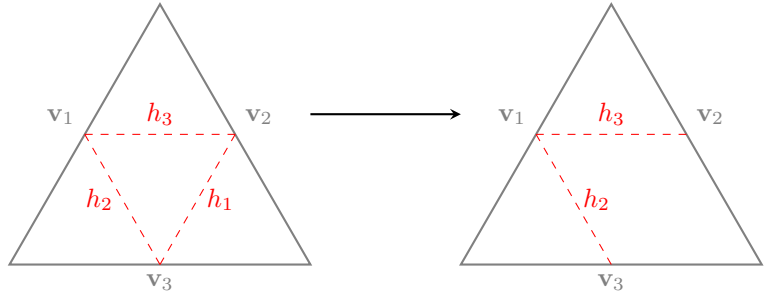

\begin{proof}
\OOO For notational convenience, we write $p=3/2$ in the proof. \EEE 
    Note that we can deduce by \BBB Proposition \EEE \ref{bounded-jump} and Corollary \ref{cor: energy bound} that $\mathcal{H}^{1}(K^{\rm L}_n(t))$ is uniformly bounded with respect to $n$ and $t$. Hence, we can apply Theorem \ref{helly-sigma-p} and thereby get a subsequence 
    $(K^{\rm L}_{n})_n$ (not relabeled) and an increasing set function $t\mapsto K(t)$ such that 
    \[K^{\rm L}_{n}(t) \;\text{$\sigma^{\MMM p \EEE}$-converges to $K(t)$ in $\Omega$}\,,\] 
    for all $t\in [0,T]$. We can repeat this argument for the variants $K^i_n(t)$ and, since $K^i_n(t) \BBB \, \tilde{\subset} \, \EEE K^{\rm L}_n(t)$ for $i=1,2,3$, the $\sigma^{\MMM p \EEE}$-limits $K^i(t)$ exist, up to passing to a further subsequence (not relabeled), and we have $K^i(t) \, \tilde{\subset} \, K(t)$ by Remark \ref{sigma-p-remark}(c). It suffices to show that $K^i(t) \, \tilde{\supset} \, K(t)$ for each $t \in [0,T]$.

Without restriction we show the argument for $i=1$. We fix a time $t\in [0,T]$, and for given $n$ we choose $k$ such that $t \in [t_n^k, t_n^{k+1})$. By the definition of $\sigma^{\MMM p \EEE}$-convergence we find a sequence $u_n\to u$ weakly in $SBV^{\MMM p \EEE}(\Omega)$ such that $S(u_n) \, \tilde{\subset} \,  K_n^{\rm L}(t)$  and $K(t) \, \tilde{=} \,  S(u)$. We now construct  a modification $v_n$ with $S(v_n)\, \tilde{\subset} \,  K^1_{n}(t)$. By construction  in \eqref{eq: hat sets}, see also \eqref{2,3},  we have $K^{\rm L}_n(t) \subset \bigcup_{\triangle\in  \Cl{k}{u^k_n}   }\triangle$. We define $v_n = u_n$ outside of these triangles and now consider each triangle in $\triangle\in \Cl{k}{u^k_n}$ separately. 

If $\triangle\in \mathcal{C}^{k,2}_n(u_n)$, then we set $v_n=u_n$ in $\triangle$. If $\triangle\in \mathcal{C}^{k,3}_n(u_n)$, we know that $S(u_n)$ lies in $ h^1_{\scalebox{1.0}{$\scriptscriptstyle \triangle$}} \cup  h^2_{\scalebox{1.0}{$\scriptscriptstyle \triangle$}} \cup  h^3_{\scalebox{1.0}{$\scriptscriptstyle \triangle$}}$. We now explain that the function can be modified such that it jumps only on $ h^2_{\scalebox{1.0}{$\scriptscriptstyle \triangle$}}\cup  h^3_{\scalebox{1.0}{$\scriptscriptstyle \triangle$}}$. In fact, as $u_n\in W^{1,\MMM p \EEE} (\Omega\setminus \overline{S(u_n)})$, we can extend $u_n$  from \BBB the sub-triangle with edge \EEE $ h^1_{\scalebox{1.0}{$\scriptscriptstyle \triangle$}}$ to the middle triangle, see Figure~\ref{fig: extension}, i.e., we replace $u_n$ by $v_n$ inside the middle triangle such that $v_n$ is continuous on $ h^1_{\scalebox{1.0}{$\scriptscriptstyle \triangle$}}$.  This construction leads to $S(v_n) \cap \triangle \subset K^{1}_n(t)$, with $\|v_n\|_{W^{1,\MMM p \EEE}(\triangle)}\leq  C\|u_n\|_{W^{1,\MMM p \EEE}(\triangle)}$. \BBB In turn, we obtain $\|\nabla v_n\|_{L^{\MMM p \EEE}(\Omega_n)} \le C$. \EEE

Now, $u_n \to u$ in $L^1(\Omega)$ along with the fact that $\mathcal{L}^2(\lbrace u_n \neq v_n \rbrace) \le \eps^2_n \# \mathcal{C}^k_n(u^k_n) \to 0$ also shows that $v_n \to u$ in $L^1(\Omega)$, i.e., $v_n \to u$ weakly in $SBV^{\MMM p \EEE}(\Omega)$. Using part (i) of Definition \ref{def: sigma-p}  
\BBB and  $S(v_n)  \, \tilde{\subset} \,  K^{1}_n(t)$ we conclude \EEE $K(t) \, \tilde{=} \,  S(u) \, \tilde{\subset} \,   K^1(t)$.
\end{proof}

% We now want to start with a convergence result for this sets in the $\sigma^2$-sense, see Definition \ref{sigma-p-definition}.  

% \begin{proposition}\label{sigmaconvergence}
%     Let $u_n(t)$ be the evolution defined in \eqref{definitionofun} and $K_n^{L}(t)$ as in \eqref{def: KL}. There exists an increasing set function $t\to K(t)$ on $[0,T]$ and a subsequence (not relabeled) such that for every $t\in [0,T]$  %$K_n(t)$ $\sigma^{2}$-converges to $K(t)$ %
%     we have \[K^{\rm L}_n(t) \quad  \text{$\sigma^2$- converges to $K(t)$ in $\Omega$\,.}\]
% \end{proposition}
% \begin{proof}
    
% \end{proof}

We now proceed with a compactness result for displacements and lower semicontinuity for the energies. 

\begin{proposition}[Compactness and lower semicontinuity]\label{liminfineq}
    For a fixed $t \in [0,T]$, let $(u_n(t))_{n\in \N}$ be the evolution defined in \eqref{definitionofun}. Let $K(t)$ be the $\sigma^{\OOO 3/2 \EEE}$-limit of ${K}^{{\rm L}}_n(t)$ as given in Proposition \ref{sigmaconvergence}. Then, there exists a subsequence $(n_l)_l$ depending on $t$ and some $u \in SBV^{\MMM 2 \EEE}(\Omega)$ such that  $u_{n_l}(t)$ AC-converges to $u$ as $l \to \infty$, $\liminf_{n \to \infty} \mathcal{E}_{n}(u_{n}(t);t) =  \liminf_{l \to \infty} \mathcal{E}_{n_l}(u_{n_l}(t);t) $,  and  
    \[  \liminf_{l \to \infty} \mathcal{E}_{n_l}(u_{n_l}(t);t) \ge  \liminf_{l \to \infty} \mathcal{E}_{n_l}^{ {\rm{ela}}}(u_{n_l}(t);t)  + \liminf_{l \to \infty} \mathcal{E}_{n_l}^{ {\rm{cra}}}(u_{n_l}(t);t) \geq \mathcal{E}(u(t), K(t))\,.\]
\end{proposition}

\begin{proof}%[Proof of Theorem %\ref{liminfineq}%]
The compactness of the displacements follows from  Proposition \ref{bounded-jump},  by using Corollary \ref{cor: energy bound} and the $L^\infty$-bound on the displacements given by Lemma \ref{lemma: dm}. In particular, \BBB  the subsequence $(n_l)_l$ can be chosen \EEE such that $\liminf_{n \to \infty} \mathcal{E}_{n}(u_{n}(t);t) =  \liminf_{l \to \infty} \mathcal{E}_{n_l}(u_{n_l}(t);t) $ holds. \EEE Then, in view of \eqref{energy:main}, Lemma \ref{elastic-lowersim},    and the nonnegativity of  $\mathcal{E}_n^{{\rm{rem}}}$ and $\mathcal{E}_n^{ {\rm{bdy}}}$,  we only need to show that 
 $$\liminf_{n \to \infty} \mathcal{E}_n^{ {\rm{cra}}}(u_{n}(t);t) \ge  \int_{ K(t)  }\varphi(\nu_{K(t)})\, {\rm d} \mathcal{H}^{1}\,.  $$   
 Observe that $\sup_k \# \mathcal{C}^k_n(u_n^k) \le C \eps_n^{-1}$ for a constant independent of $n \in \N$, by Corollary \ref{cor: energy bound} and \eqref{eq: jump bound}. Then, using Proposition \ref{prop: repri}, $R_n \to \infty$ by \eqref{rnepsn}, and $\lim_{ s  \to \infty}\Psi(s)=\kappa$, we have
$$
\liminf_{n \to \infty} \mathcal{E}_n^{{\rm{cra}}}(u_{n}(t);t)  = \liminf_{n \to \infty}   \frac{\kappa}{\sqrt{3}}\sum_{i=1}^{3}\int_{{K}^{i}_n(t)}  \varphi_i(\nu_{{K}^{i}_n(t)}) \,{\rm d}\mathcal{H}^1  ,   
$$
where $\varphi_\alpha(\nu):=(|\mathbf{v}_\beta\cdot\nu|+|\mathbf{v}_\gamma\cdot \nu|)$ for $\nu \in \mathbb{S}^1$, \BBB and \EEE for $\alpha,\beta,\gamma \in \lbrace 1,2,3 \rbrace $ pairwise distinct.  We now  employ Proposition \ref{sigmaconvergence} and Theorem \ref{helly-sigma-p}, i.e., the lower semicontinuity with respect to $\sigma^{\OOO 3/2 \EEE}$-convergence, to obtain
\[\liminf_{n \to \infty} \mathcal{E}_n^{{\rm{cra}}}(u_{n}(t);t)  \geq\frac{\kappa}{\sqrt{3}}\sum_{i=1}^{3}\int_{K(t)} \varphi_i(\nu_{K(t)})\,{\rm d}\mathcal{H}^1\,.\]
By the definition of $\varphi$ in \eqref{varphidef}, the proof is concluded. 
\end{proof}

\subsection{Proof of the main result}\label{sec: 5.3}

In this subsection we give the proof of Theorem \ref{maintheorem}. As a last preparation, we need the following stability result for our time-discrete atomistic energy. Recall the splitting of the energy in \eqref{energy:main}. \BBB For $t \in [0,T]$ and $n \in \N$, we let $k(t)$ be the ($n$-dependent) index such that $t \in [t^{k(t)}_n,t^{k(t)+1}_n)$. \EEE

\begin{theorem}[Stability]\label{stability}
    Let $t \mapsto u_n(t)$ be the evolution defined by \eqref{definitionofun} and let $K(t)$ be the $\sigma^{\OOO 3/2 \EEE}$-limit of $K_n^{{\rm L}}(t)$. For each $\psi\in SBV^{\MMM 2 \EEE}(\Omega)$ with $\psi = g(t)$ on $\Omega \setminus \overline{U}$ there exists a sequence $(\psi_{n})_{n}$ of discrete displacements with $\psi_n \in   \mathcal{A}_{n}\BBB (g(t_n^{k(t)})) \EEE$  such that $\psi_n$ AC-converges to $\psi$ and
    \begin{equation}\label{stab1}
    \limsup_{n\to \infty} \, \Big(  \big(
   \mathcal{E}^{ {\rm cra}}_{n}(\psi_n;t) +  \mathcal{E}^{ {\rm rem}}_{n}(\psi_n;t) \big) - \big(   \mathcal{E}^{ {\rm cra}}_{n}(u_n(t);t) + \mathcal{E}^{ {\rm rem}}_{n}(u_n(t);t)  \big) \Big)
    \leq \int_{S({\psi})\setminus K(t) } \varphi(\nu_\psi)\,d \mathcal{H}^{1}\,,
    \end{equation}
  where  $\nu_\psi$ denotes the normal to the jump set $S(\psi)$  of $\psi$. Moreover,
   \begin{equation}\label{do not remove}
        \lim_{n\to \infty}\,  \mathcal{E}^{ {\rm ela}}_{n}(\psi_n;t)  =    \int_{\Omega } \Phi(\nabla \psi) \, {\rm d}x,
    \end{equation}
    where $\Phi$ is given in \eqref{eq: lineraition}, and we have 
    \begin{equation}\label{stab3}
        \lim_{n\to \infty} \mathcal{E}^{\rm{bdy}}_{n}(\psi_n; t )=0\,.
    \end{equation}
\end{theorem}

This result is essential to pass from the minimality condition \eqref{minimizing-scheme} in the atomistic setting to the global stability in \eqref{finalstability}. For this reason, estimates of this kind are often referred to as \emph{stability of unilateral minimizers}, see e.g.\ \cite{GiacPonsi}. The proof will be deferred to   Section \ref{stabilitysection}. The strategy follows the one of the \emph{jump transfer lemma} introduced in {\sc Francfort and Larsen} \cite{Francfort-Larsen:2003}, see the  works \cite{dMasoFranToad, Lazzaroni, FriedrichSolombrino, Giacomini:2005b} for several variants. Our situation is more delicate for two reasons. First, due to the underlying atomistic nature, special care is needed in the construction of the sequence $(\psi_n)_n$, in particular concerning the corresponding `broken' springs. Secondly, whereas in \cite{Francfort-Larsen:2003} only the location of the crack is relevant, which in our notation means \BBB   $K^{\rm L}_n(t)$, \EEE  our atomistic energy also takes into account \emph{how much each spring} was deformed at former times. Therefore, we need to verify that the exact value of the memory variable is indeed negligible in the limit, see also \cite{Giacomini:2005b} for a similar problem.

For the proof, we will also need the following corollary, which  formally  follows from Theorem \ref{stability} by choosing $u_n(t) \equiv 0$ and $K(t) = 0 $. \BBB Recall the energy $\mathscr{E}_n$ defined in \eqref{basicenergy-neu}. \EEE

\begin{corollary}[Recovery sequence]\label{cor: stability}
 For each $\psi\in SBV^{\MMM 2 \EEE}(\Omega)$ with $\psi = g(t)$ on $\Omega \setminus \overline{U}$ there exists a sequence $(\psi_{n})_{n}$ with $\psi_n \in   \mathcal{A}_{n} \BBB (g(t_n^{k(t)})) \EEE $ such that $\psi_n$ AC-converges to  $\psi$  and
    \begin{equation*}%\label{stab1.111}%
    \limsup_{n\to \infty} \mathscr{E}_n(\psi_n) \le \mathcal{E}(\psi,S(\psi)).
    \end{equation*}
\end{corollary}
We are now in \BBB the \EEE position to prove the main result of this paper. 

\begin{proof}[Proof of Theorem \ref{maintheorem}]
We split the proof into  five   steps.

\emph{Step 1: Limiting evolution.} Let  $K(t)$ be the $\sigma^{\OOO 3/2 \EEE}$-limit of $K_n^{{\rm L}}(t)$ for $t \in [0,T]$ given by Proposition~\ref{sigmaconvergence}. Since   $t \mapsto K(t)$ is an increasing set function, \BBB i.e., \EEE $K(t_1) \, \tilde{\subset} \,  K(t_2)$ for all $0\leq t_1\leq t_2\leq T$, \BBB the \EEE irreversibility condition (a),  is satisfied. For each fixed time $t\in [0,T]$, by virtue of Proposition~\ref{liminfineq}, there exists $u(t) \in SBV^{\MMM 2 \EEE}(\Omega)$ and a subsequence depending on $t$ (not relabeled) such that $u_n(t)$ AC-converges to $u(t)$. In particular, as $u_n(t)\in \mathcal{A}_n\BBB (g(t_n^{k(t)})) \EEE$, this also shows that $u(t) = g(t)$ on $\Omega \setminus \overline{U}$. In view of Corollary \ref{cor: energy bound} and Lemma~\ref{lemma: dm} we can apply Proposition \ref{crucial-inclusion} to deduce $S(u(t))\, \tilde{\subset} \, K(t)$. This shows that $u(t) \in AD(g(t),K(t))$ for all $t \in [0,T]$.

Next, we check that  the  mapping  $t\to (u(t),K(t))$ is an irreversible quasi-static crack evolution. For this, we need to confirm \eqref{finalstability} and \eqref{energybalance} which is content of Step 2 and Step 3--4, respectively. In the last step, we then prove the convergence of energies and show that the crack sets also converge in the sense of $\sigma$-convergence, see Definition \ref{def: sigma conv}.
 
\emph{Step 2: Stability  \eqref{finalstability}.} Fix $t \in [0,T]$. Let $\psi\in AD(g(t),H)$ be an admissible competitor for a set $H$ with $K(t) \, \tilde{\subset} \, H$.  We employ Theorem \ref{stability} to obtain a sequence of discrete displacements $\psi_n \in \mathcal{A}_{n}\BBB (g(t_n^{k(t)})) \EEE$ such that \eqref{stab1}--\eqref{stab3} are valid. By the minimality property of the solution   $u_n(t)$, see \eqref{minimizing-scheme0}--\eqref{minimizing-scheme}, and the shorthand notation in \eqref{energy along evo2} we have  
    \[\mathcal
    {E}_{n}(u_n(t); t)\leq \mathcal{E}_{n}(\psi_{n} ;t)\,.\]
    We now split the energy on both sides like in \eqref{energy:main}: Subtracting the crack energy and the remainder term on the left-hand side and using that $ \mathcal{E}_n^{{\rm{bdy}}}$ is nonnegative, we obtain   
    \begin{equation}\label{Passtolimit}
    \begin{aligned}    
   \mathcal{E}_n^{{\rm{ela}}}(u_{n}(t);t)\leq  & \;\mathcal{E}_n^{{\rm{ela}}}(\psi_{n}; t) + \mathcal{E}_n^{{\rm{cra}}}(\psi_{n};t)-\mathcal{E}_n^{{\rm{cra}}}(u_{n}(t);t) + \mathcal{E}_n^{{\rm{rem}}}(\psi_{n};t) -  \mathcal{E}_n^{{\rm{rem}}}(u_{n}(t);t)    + \mathcal{E}_n^{{\rm{bdy}}}(\psi_{n};t)
       \,.
    \end{aligned}
\end{equation}
Passing to the limit in \eqref{Passtolimit}, by employing \eqref{stab1}--\eqref{stab3}, we find 
$$\limsup_{n\to \infty}  \mathcal{E}_n^{{\rm{ela}}}(u_{n}(t);t)  \le  \int_{\Omega } \Phi(\nabla \psi) \, {\rm d}x +  \int_{S({\psi})\setminus K(t) } \varphi(\nu_\psi)\,{\rm d} \mathcal{H}^{1}.$$
Since \BBB ${u}_{n}(t)$ AC-converges to $u(t)$ \EEE (for a $t$-dependent subsequence, not relabeled), \BBB employing also Lemma~\ref{lemma: dm}~and Corollary~\ref{cor: energy bound}, \EEE we can use Lemma~\ref{elastic-lowersim} and derive  
    \begin{equation}\label{estimateelastic}
        \int_{\Omega}\Phi(\nabla u(t))\, {\rm d}x \leq \limsup_{n\to \infty} \mathcal{E}_n^{{\rm{ela}}}(u_{n}(t);t)  \le  \int_{\Omega } \Phi(\nabla \psi) \, {\rm d}x +  \int_{S(\psi)\setminus K(t) } \varphi(\nu_\psi)\,{\rm d} \mathcal{H}^{1}.
    \end{equation}
    Since $\psi\in AD(g(t),H)$, we have $S(\psi) \BBB \, \tilde{\subset} \, \EEE H$ and thus 
    \[\int_{\Omega} \Phi(\nabla u(t))\, {\rm d}x \leq \int_{\Omega} \Psi(\nabla \psi)\,{\rm d}x+ \int_{H\setminus K(t) } \varphi(\nu_H)\,{\rm d} \mathcal{H}^{1}\,.\]   
Since $K(t) \ \tilde{\subset} \,  H$, we hence  conclude 
    \[\int_{\Omega}\Phi(\nabla u(t))\,   {\rm d}x + \int_{K(t)} \varphi(\nu_{K(t)})\,{\rm d} \mathcal{H}^{1}\,\leq \int_{\Omega} \Phi(\nabla \psi)\,{\rm d}x+ \int_{H }\varphi (\nu_H)\,{\rm d} \mathcal{H}^{1}\,.\]
Recalling the definition of the limiting energy \eqref{eq: lim-en}, this shows \eqref{finalstability}.

In particular, for $t = 0$, denoting by $(\psi^{0}_n)_n$   a recovery sequence for $u(0)$ from Corollary \ref{cor: stability},  the minimality property in \eqref{minimizing-scheme0} and  Proposition \ref{liminfineq}  yield
$$ \mathcal{E}(u(0), K(0)) \le \liminf_{n \to \infty} \mathscr{E}_n(u_{n}(0)) \le  \limsup_{n \to \infty} \mathscr{E}_n(\psi^{0}_n) \le   \mathcal{E}(u(0), S(u(0))) \,.$$
This along with $u(0)\in AD(g(0),K(0))$, i.e., $S(u(0)) \, \tilde{\subset} \, K(0)$, shows
\begin{align}\label{energy at 0 time}
  \lim_{n \to \infty} \mathscr{E}_n(u_{n}(0))  =  \mathcal{E}(u(0), K(0)) =\mathcal{E}(u(0), S(u(0))). 
\end{align}

\emph{Step 3: Convergence of work by  applied boundary load}.
Next, we check that for all $t \in [0,T]$ 
\begin{align}\label{eq: to checcc}
\lim_{n \to \infty}   \int_{0}^{t} \int_{\Omega_n} D \Psi_n^{\rm cell }(\nabla \hat{u}_n(\tau)) \cdot \nabla \partial_{t}\tilde{g}_n(\tau)\,{\rm d} x\,{\rm d}\tau =  \int_{0}^{t} \int_{\Omega} D\Phi(\nabla u(\tau)) \cdot \nabla \partial_{t} g(\tau)\, {\rm d}x \, {\rm d}\tau\,.
\end{align}
\BBB Let $\tau\in [0,t)$. \EEE Using the stability estimate \eqref{estimateelastic} for $\psi =  u(\tau) \in AD(g(\tau), K(\tau))$, and employing $S(u(\tau)) \, \tilde{\subset} \, K(\tau)$ we get 
$$ \lim_{n\to \infty}  \mathcal{E}_n^{{\rm{ela}}}(u_{n}(\tau);  \tau   )  =    \int_{\Omega}\Phi(\nabla u(\tau))\, {\rm d}x \quad \text{for all $\tau \in [0,t)$\,.} $$
Then \eqref{eq: ela en} and Lemma \ref{rem: the new one} imply  $\chi_{\Omega_n} \nabla  \hat{u}_n(\tau) \to \nabla u(\tau) $ strongly in $L^{\MMM 2 \EEE}(\Omega;\R^2)$. 
Thus, as in Step 1 of the proof of Lemma \ref{e-balance-lemma1}, we can employ   \cite[Lemma 2.4]{GiacPonsi}   to find  
$$ \Big|  \int_{\Omega_{n}}   D \Psi^{\rm cell}_n  \big(\nabla \hat{u}_{n}(\tau)  \big)   \cdot \nabla \partial_{t} g(\tau)\,{\rm d}x -   \int_{\Omega_{n}}   D \Psi^{\rm cell}_n  \big(\nabla {u}(\tau)  \big)   \cdot \nabla \partial_{t} g(\tau)\,{\rm d}x \Big| \to 0, $$
where we also used Lemma \ref{lem:help2}(iii),(iv)  and the fact that $\nabla \partial_{t} g(\tau)\in L^{\MMM 2 \EEE}(\Omega; \R^2)$ for almost every $\tau \in [0,t)$.  
Since $\nabla {u}(\tau)  \in L^{\MMM 2 \EEE}(\Omega;\R^2)$ and $D \Psi^{\rm cell}_n \to D\Phi $ pointwise, we also get by dominated convergence
$$ \Big|  \int_{\Omega_{n}}   D \Psi^{\rm cell}_n  \big(\nabla \hat{u}_{n}(\tau)  \big)   \cdot \nabla \partial_{t} g(\tau)\,{\rm d}x -   \int_{\Omega_{n}}   D \Phi  \big(\nabla {u}(\tau)  \big)   \cdot \nabla \partial_{t} g(\tau)\,{\rm d}x \Big| \to 0. $$
At this stage, we can repeat the argumentation in \eqref{eq: nochmal}, in particular using Lemma \ref{lem:help2}(iii),   the fact that $\chi_{\Omega_n}\nabla \partial_t \tilde{g}_n$ converges strongly to $ \nabla \partial_t  g$ in $L^1([0,T];L^{\infty}(\Omega;\R^2))$,   %$g \in W^{1,1}([0,T]; W^{2,\infty}(\Omega))$%
and the uniform bound on $\nabla \hat{u}_n(\tau)$ in $L^{\MMM 2 \EEE}(\Omega_n;\R^2)$, see \eqref{eq: hat finally}, to obtain \eqref{eq: to checcc}.

\emph{Step 4: Energy balance \eqref{energybalance}.}   For $t \in [t_n^k,t_n^{k+1})$, by Lemma \ref{e-balance-lemma1} we find $(\eta_n)_n$ with $\eta_n \to 0$ such that
\begin{align*}
      \mathcal{E}_n^k(u_n(t)) - \mathcal{E}_{n}^{0}(u^0_n) & \leq \int_{0}^{t_{n}^{k}} \int_{\Omega_n} D \Psi_n^{\rm cell }(\nabla \hat{u}_n(\tau)) \cdot \nabla \partial_{t}\tilde{g}_n(\tau)\,{\rm d} x\,{\rm d}\tau  +  ( 1  + \eps_n \# \mathcal{C}^k_n) \eta_n   \\
      & \le    \int_{0}^{t} \int_{\Omega_n} D \Psi_n^{\rm cell }(\nabla \hat{u}_n(\tau)) \cdot \nabla \partial_{t}\tilde{g}_n(\tau)\,{\rm d} x\,{\rm d}\tau +     C\eta_n +  C\delta_n, 
      \end{align*}      
where in the second step we have used  Proposition \ref{bounded-jump} and Corollary \ref{cor: energy bound}  to get $\# \mathcal{C}^k_n \le C\eps_n^{-1}$.   Furthermore, we used the first estimate in \eqref{eq: nochmal} and $g\in W^{1,1}([0,T]; W^{2,\infty}(\Omega))$ to obtain a uniform control \BBB on \EEE the integrand, which allows to estimate the integral from $t_n^k$ to $t$ in terms of $C\delta_n$. 
Combining this with Proposition \ref{liminfineq} and \eqref{energy at 0 time}, and using that $\mathcal{E}^0_n = \mathscr{E}_n$, we get by \eqref{eq: to checcc}
\begin{align}\label{eq: for en ba}
        \mathcal{E}(u(t), K(t)) & \leq \liminf_{n\to \infty} \mathcal{E}_{n}(u_n(t);t) \leq \limsup_{n\to \infty} \mathcal{E}_{n}^{0}(u_{n}^0) +  \limsup_{n\to \infty} \int_{0}^{t} \int_{\Omega_n} D \Psi_n^{\rm cell }(\nabla \hat{u}_n(\tau)) \cdot \nabla \partial_{t}\tilde{g}_n(\tau)\,{\rm d} x\,{\rm d}\tau \notag \\ &=  \mathcal{E}(u(0), K(0)) +  \int_0^t  \int_{\Omega} D\Phi(\nabla u(\tau)) \cdot \nabla \partial_{t} g(\tau)\, {\rm d}x \, {\rm d}\tau.
\end{align}
Thus, it remains to prove the reverse inequality
\begin{equation*}
    \mathcal{E}(u(t), K(t))\geq \,\mathcal{E}(u(0),K(0))+\int_0^t  \int_{\Omega} D\Phi(\nabla u(\tau)) \cdot \nabla \partial_{t} g(\tau)\, {\rm d}x \, {\rm d}\tau.
\end{equation*} \BBB
This can be done by repeating the arguments in \cite[Lemma 7.1]{dMasoFranToad} which are based on the  global stability result \eqref{finalstability} proved in \emph{Step 2}, as well as approximation of integrals by Riemann sums and their convergence, see \cite[Lemmas 5.12 and 5.7]{dMasoFranToad}. \EEE

\emph{Step 5: Energy convergence and $\sigma$-convergence of crack sets}.    Combining \eqref{eq: for en ba} with the energy balance \eqref{energybalance} we find
\begin{align}\label{en convo3}
     \mathcal{E}(u(t), K(t))  = \lim_{n\to \infty} \mathcal{E}_{n}(u_n(t);t) \quad \quad \text{ for all $t \in [0,T]$,} 
     \end{align}
which shows \eqref{energ convi}. To conclude the proof it remains to show that \eqref{eq: sigma coniiiii} holds.

Consider the crack sets $K_n(t)$ defined in \eqref{eq: first K def}. We note that $\mathcal{H}^1(K_n(t)) \le C$ for all $t \in [0,T]$ and $n \in \N$ by \eqref{eq: jump bound},  \BBB Lemma~\ref{lemma: dm}, \EEE  and Corollary \ref{cor: energy bound}.  As  $K_n(t)$ is increasing in $t$, Remark \ref{wurschel remark} implies that there is a subsequence (not relabeled) and an increasing set function $K^\sigma(t)$ such that $K_n(t)$ $\sigma$-converges to  $K^\sigma(t)$ for each $t \in [0,T]$. To conclude the proof, it suffices to check that $K^\sigma(t) \, \tilde{=} \,  K(t)$ for all $t \in [0,T]$, where $K(t)$ is given in Proposition \ref{sigmaconvergence}. We note that $K^\sigma(t) \, \tilde{\supset} \,  K(t)$ by Remark \ref{sigma-p-remark}(d), so it suffices to check that $K^\sigma(t) \, \tilde{\subset} \,  K(t)$. 

Recall the \BBB splitting \EEE of $\mathcal{E}_{n}(u_n(t);t)$ in \eqref{energy:main}. The energy convergence in \eqref{en convo3} along with the lower-semicontinuity for $\mathcal{E}_{n}^{ {\rm{cra}}}$ and $\mathcal{E}_{n}^{ {\rm{ela}}}$ from   Lemma \ref{elastic-lowersim}   and Proposition \ref{liminfineq}, and the fact that $\mathcal{E}_{n}^{ {\rm{rem}}}$ and $\mathcal{E}_{n}^{ {\rm{bdy}}}$ are nonnegative show 
$$ \lim_{n \to \infty} \mathcal{E}_{n}^{ {\rm{rem}}}(u_{n}(t);t) = \lim_{n \to \infty} \mathcal{E}_{n}^{ {\rm{bdy}}}(u_{n}(t);t)  = 0 $$
and 
\begin{align}\label{en convo1}
 \liminf_{n \to \infty} \mathcal{E}_{n}^{ {\rm{ela}}}(u_{n}(t);t) =  \int_{\Omega}\Phi(\nabla u(t))\,{\rm d}x, \quad \quad 
 \liminf_{n \to \infty} \mathcal{E}_{n}^{ {\rm{cra}}}(u_{n}(t);t)  = \int_{K(t)}\varphi(\nu_{K(t)})\, {\rm d} \mathcal{H}^{1}\, . 
 \end{align}
\BBB In view of  Remark \ref{remark-about-2Ks}, we obtain 
\[\mathcal{H}^{1}\big(K_n(t) \setminus K^{\rm L}_n(t)\big)\leq   \frac{3}{2} \varepsilon_n \, \# \big(\Cn{k}{u^{k}_n} \setminus \Cl{k}{u^{k}_n}\big)  +   \varepsilon_n \, \# \big\{ \mathcal{C}^{k,2}_n ({u^{k}_n}) \colon M_{n,\triangle}^{k,\mathbf{v}_i} (u^k_n) > R \text{ for all $i=1,2,3$}\big\}  \,,\]
where the second sum is related to the triangles with (in the notation of Remark \ref{remark-about-2Ks})   $\mu_{n,\triangle}^k = 3$  and $\nu_{n,\triangle}^k = 2$.  Recalling the definition of $\mathcal{E}_{n}^{ {\rm{rem}}}$ in \eqref{eq: different enegies} and using $\Psi(R) >0$ 
$$    \# \big(\Cn{k}{u^{k}_n} \setminus \Cl{k}{u^{k}_n}\big)  +    \# \big\{ \mathcal{C}^{k,2}_n( {u^{k}_n}) \colon M_{n,\triangle}^{k,\mathbf{v}_i} (u^k_n) > R \text{ for all $i=1,2,3$}\big\} \le C \eps_n^{-1}\mathcal{E}_{n}^{ {\rm{rem}}}(u_{n}(t);t).$$
This shows $\mathcal{H}^{1}(K_n(t) \setminus K^{\rm L}_n(t)) \to 0$ as $n \to \infty$. \EEE By Remark \ref{sigma-p-remark}(b) we therefore get that the $\sigma$-limit of  $K_n(t)$ and $K^{\rm L}_n(t)$ coincide, i.e., $K^{\rm L}_n(t)$ $\sigma$-converges to $K^\sigma(t)$ for all $t \in [0,T]$. Arguing as in the proofs of Propositions \ref{sigmaconvergence} and \ref{liminfineq},  and using Remark \ref{wurschel remark} in place of Theorem \ref{helly-sigma-p}, we get
\begin{align*}
 \liminf_{n \to \infty} \mathcal{E}_n^{{\rm{cra}}}(u_{n}(t);t)  \geq \int_{K^\sigma(t)} \varphi(\nu_{K^\sigma(t)})\,{\rm d}\mathcal{H}^1.
 \end{align*}
This along with \eqref{en convo1} and the inclusion $K(t) \, \tilde{\subset} \, K^\sigma(t)$ shows $K^\sigma(t) \tilde{=} K(t) $ and concludes the proof. 
\end{proof}

\section{Proof of the stability result}\label{stabilitysection}
This section is entirely devoted to  the proof of the stability result in Theorem \ref{stability} and Corollary \ref{cor: stability}. 

\subsection*{Density argument} 
We start by observing that it suffices to prove the statement for functions $\psi$ with more regularity, employing a suitable density argument. 
Let $\mathcal{W}(\Omega ) \subset SBV(\Omega )$ be the collection of functions $v$ such that $S(v)$ is closed and included in a finite union of closed and  connected  pieces of $C^1$-curves and $v$ lies in $W^{2,\infty}(\Omega  \setminus S({v}))$. For each $v \in SBV^{\MMM 2 \EEE}(\Omega )$ with $v = g(t)$ on $\Omega  \setminus \overline{U}$ we can choose a sequence $(v_n)_n \subset\mathcal{W}(\Omega ) $ with   $v_n = g(t)$ \OOO on \EEE  $\Omega  \setminus \overline{U}$   such that 
$$ v_n \to  v  \text{ in } L^1(\Omega ), \quad \Vert \nabla v_n  - \nabla v \Vert_{L^{\MMM 2 \EEE}(\Omega )} \to 0, \quad  \mathcal{H}^{1}\big(  S(v_n) \setminus  S(v)\big)  + \mathcal{H}^1\big(S(v) \setminus  S(v_n)   \big) \to  0.
$$
This follows from \cite[Theorem C]{Pratelli} (see also \cite[Theorem 1.1]{Crismale3} for the control on the Sobolev norms of $\nabla v_n$), where a construction similar to the one performed in \cite[Proposition 2.5]{Giacomini:2005} ensures that the boundary values are attained along the sequence. \OOO We also refer to \cite[Theorem 3.2]{steinke} for an analogous statement and proof for $GSBD$-functions. \EEE

With this density result at hand, we observe that it suffices to construct a sequence as in Theorem~\ref{stability} for a function $\psi \in \mathcal{W}(\Omega )$ with  $\psi = g(t)$ in a neighborhood of  $\Omega  \setminus \overline{U}$ in $\Omega$. The general case then follows by a standard diagonal argument. Further, we fix $\theta>0$ and observe that it suffices to construct a sequence $(\psi_n)_n $ of  discrete displacements with $\psi_n \in   \mathcal{A}_{n} \BBB (g(t)) \EEE $ such that  
  \begin{equation}\label{stab0neu}
    \limsup_{n\to \infty}  \Vert  \chi_{\Omega_n} \tilde{\psi}_n  -   \psi  \Vert_{L^1(\Omega )} \le C\theta\,,
    \end{equation}
    \begin{equation}\label{stab1neu}
    \limsup_{n\to \infty} \, \Big(  \big(
   \mathcal{E}^{ {\rm cra}}_{n}(\psi_n;t) +  \mathcal{E}^{ {\rm rem}}_{n}(\psi_n;t) \big) - \big(   \mathcal{E}^{ {\rm cra}}_{n}(u_n(t);t) + \mathcal{E}^{ {\rm rem}}_{n}(u_n(t);t)  \big) \Big)
    \leq \int_{S(\psi)\setminus K(t) } \varphi(\nu_\psi  )\,d \mathcal{H}^{1} + C\theta\,,
    \end{equation}
   \begin{equation}\label{stab2neu}
        \limsup_{n\to \infty}\,  \mathcal{E}^{ {\rm ela}}_{n}(\psi_n;t)  \le     \int_{\Omega  } \Phi(\nabla \psi) \, {\rm d}x + C\theta,
    \end{equation}
       \begin{equation}\label{stab3neu}
        \lim_{n\to \infty} \mathcal{E}^{\rm{bdy}}_{n}(\psi_n; t )=0\,,
    \end{equation}
    where $\Phi$ and $\varphi$ are given in \eqref{eq: lineraition} and \eqref{varphidef}, respectively, and $C>0$ is a universal constant. Then, the statement follows again by a diagonal argument, sending $\theta \to 0$.  \BBB Note that, strictly speaking, proceeding in this way, the  boundary values $g(t^{k(t)}_n)$ are not satisfied. Therefore, we eventually need to replace the sequence $(\psi_n)_n $ by $\psi_n(x) -g(t,x) + g(t^{k(t)}_n,x)$ for $x \in \mathcal{L}_n(\Omega)$. Due to the regularity of $g$ and the Lipschitz continuity of $\Psi$, this still leads to  \eqref{stab0neu}--\eqref{stab3neu}, cf.\ also \eqref{wanttoprove} for a similar estimate.

     Now, we fix $\psi \in \mathcal{W}(\Omega )$      with $\psi = g(t)$ \OOO on \EEE $\Omega  \setminus \overline{U}$, \EEE and construct a sequence $(\psi_n)_n$ satisfying \eqref{stab0neu}--\eqref{stab3neu} \BBB and  $\psi_n \in   \mathcal{A}_{n}  (g(t))  $. \EEE 

\subsection*{Besicovitch covering}

Following closely the procedure in \cite{Francfort-Larsen:2003}, we introduce a fine cover of the jump set $S(\psi)$ with closed squares satisfying   certain additional properties. \OOO For simplicity we only treat the case that $\mathcal{H}^1(S(\psi) \cap \partial_D U)=0$ (no jump along the boundary), for the general case follows by minor adaptations of the construction at the boundary (see \cite{Francfort-Larsen:2003}) which would merely overburden notation in the sequel. As a consequence, \EEE we can assume that all these squares are contained in $U$. We furthermore denote by $\nu_\psi$ a measure-theoretic normal at $S(\psi)$.

Because of the $\sigma^{\OOO 3/2 \EEE}$-convergence of $K_n^{\rm L}(t)$ to $K(t)$, we have a function $v \in SBV^{\OOO 3/2 \EEE}(\Omega )$ and a sequence $(v_n)_n \subset SBV^{\OOO 3/2 \EEE}(\Omega )$ such that $v_n \to v$ in $L^1(\Omega )$, $\nabla v_n \rightharpoonup \nabla v$ weakly in $L^{\OOO 3/2 \EEE}(\Omega ;\R^2)$, $S(v) \BBB \, \tilde{=} \,  \EEE K(t)$, and $S(v_n) \BBB \,  \tilde{\subset} \, \EEE   K_n^{{\rm L}}(t)$. We can choose a suitable subset $G_j \subset S(v)  = K(t)$ as done preceding  \cite[(2.2)]{Francfort-Larsen:2003} such that
\begin{align}\label{eq gj}
\mathcal{H}^1( S(v) \setminus G_j  ) = \mathcal{H}^1( K(t) \setminus G_j  ) \le \theta
\end{align}  
and for each $x\in G_j \cap S(\psi) $ we consider closed squares $Q_r(x)$ with sidelength $2r$ and \BBB two sides \EEE orthogonal to $\nu_\psi(x)$ which are contained in $U $ and satisfy \cite[(2.3), (2.5)]{Francfort-Larsen:2003}. 
 
For closed squares $Q_r(x) \subset U $ with a center $x\in S(\psi)\setminus K(t)$, still oriented in direction $\nu_\psi(x)$,  we can assume that for $\mathcal{H}^1$-a.e.\ $x\in S(\psi) \setminus K(t)$ and for $r$ sufficiently small it holds
    \begin{equation}\label{almostnoK}
\mathcal{H}^{1}\big(K(t) \cap Q_{r}(x)\big)\leq \theta r\,.
    \end{equation}
This is possible by the fact that $K(t)$ has $\mathcal{H}^{1}$-density $0$ almost everywhere in $S(\psi)\setminus K(t)$. As $S(\psi)$ is contained in a finite union of closed $C^1$-curves, for a.e.\ $x \in S(\psi)$, possibly passing to smaller $r$, the above squares   can be chosen such that they also   satisfy 
\begin{equation}\label{normal-angle-estimate0}
2r \le \mathcal{H}^1\big(S(\psi) \cap Q_r(x) \big) \le 4r,
\end{equation}    
    \begin{equation}\label{normal-angle-estimate2}
S(\psi) \cap Q_{r \BBB (1+\theta)  }(x)  \subset \lbrace y \colon   |(y-x) \cdot \nu_\psi(x)| \le \theta r \rbrace,        
    \end{equation}
 \begin{equation}\label{normal-angle-estimate}
    |\nu_{\psi}(y)-\nu_{\psi}(x)|\leq \theta \quad \forall y \in S(\psi) \cap  Q_{r}(x) .     
    \end{equation}      
With this, we obtain a fine cover of  $\Gamma:= (G_j \cap  S(\psi)) \cup (S(\psi) \setminus K(t))$ to which we can apply the Besicovitch covering theorem  with respect to the Radon measure $\mathcal{L}^{2}+\mathcal{H}^{1}|_{\Gamma}$. For $\theta>0$ fixed as above,  we hereby find a finite and disjoint subcollection $\mathcal{B}:= (Q_{r_i}(x_i))_i$, or shortly denoted by $(Q_i)_i$, such that  $(Q_i)_i$ satisfy the properties mentioned above, in particular \eqref{normal-angle-estimate0}--\eqref{normal-angle-estimate} and for $x_i \notin K(t)$ also \eqref{almostnoK}, as well as 
    \begin{equation}\label{besicovitch-props}
            \mathcal{L}^2\big(\bigcup\nolimits_{\mathcal{B}}  Q_{i} \EEE \big)  \le   \theta^{2},    
            \quad \quad \quad 
            \mathcal{H}^{1}\big( \Gamma \setminus \bigcup\nolimits_{\mathcal{B}}Q_i\big)  \le \theta\,.
        \end{equation} 
Here and in the following, we use $\bigcup_{\mathcal{B}} Q_i$ as a shorthand for $\bigcup_{Q_{r_i}(x_i) \in \mathcal{B}} Q_{r_i}(x_i)$. Note that this implies
      \begin{equation}\label{besicovitch-props2}   
            \int_{\bigcup\nolimits_{\mathcal{B}}Q_i} |\nabla v|\, {\rm d}x    \le C\theta, \quad \quad \quad  \mathcal{H}^{1}\big( S(\psi)\setminus \bigcup\nolimits_{\mathcal{B}}Q_i\big) \le 2\theta.
              \end{equation} 
Indeed, the first property follows from \eqref{besicovitch-props},  H\"older's inequality,  and  \BBB $v \in SBV^{\OOO 3/2 \EEE}(\Omega)$, \EEE whereas the second is a consequence of \eqref{besicovitch-props}, together with the definition of $\Gamma$  and \eqref{eq gj}. Without further notice,  we will frequently use that the    squares are pairwise disjoint. 

By $\mathcal{B}_{\rm good}\subset \mathcal{B}$ we denote the collection of squares $Q_i = Q_{r_i}(x_i)$ with $x_i\in S(\psi)\setminus K(t)$, and similarly we let $\mathcal{B}_{\rm bad}\subset \mathcal{B}$ be the collection of all squares $Q_i =Q_{r_i}(x_i)$ with $x_i\in S(\psi)\cap \BBB G_j \EEE $. Clearly, we have $\mathcal{B}=\mathcal{B}_{\rm good}\cup \mathcal{B}_{\rm bad}$.  We also define the sets 
$$B_{\rm good} =  \bigcup\nolimits_{\mathcal{B}_{\rm good}} Q_i, \quad \quad \quad B_{\rm bad} =  \bigcup\nolimits_{\mathcal{B}_{\rm bad}} Q_i. $$

\subsection*{Good squares}

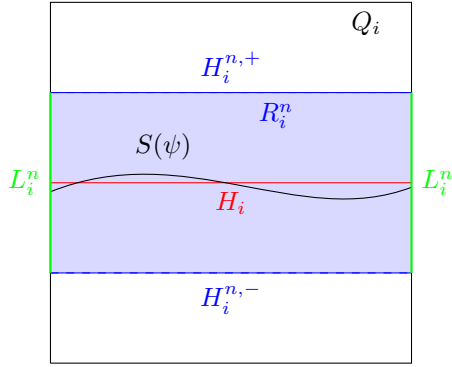
\begin{figure}
    \begin{tikzpicture}[scale=0.6]
        \filldraw[fill=blue!15!white,draw=blue]{} (0,6)-- (8,6) --(8,2)-- (0,2)--cycle;
        \node[anchor=north] at (7,8) {$Q_i$};
        \draw[color=black] (0,0) -- (8,0) --  (8,8) -- (0,8) -- (0,0);
        \node[anchor=south, color=blue] at (5,5) {$R^n_i$};
        \draw[color=blue, dashed] (0,2) -- node[anchor=north] {$H_i^{n,-}$}(8,2);
        \draw[color=blue, dashed] (0,6) --node[anchor=south] {$H_i^{n,+}$}(8,6);
        \draw[color=red] (0,4)--node[anchor=north] {$H_i$}(8,4);
        \draw (0,3.8) .. controls (3,5)  and (5.5,2.95)  .. (8,3.9);
        \node[anchor=south] at (2.5,4.3) {$S(\psi)$}; 
        \draw[color=green, thick](0,2) -- node[anchor=east] {$L^n_i$}(0,6);
        \draw[color=green, thick](8,2) -- node[anchor=west] {$L^n_i$} (8,6);
    \end{tikzpicture} \caption{Illustration of some notions used in the construction. Note  that $L^n_i$, $R^n_i$, and $H_i^{n,\pm}$ \BBB are independent of $n$ if $Q_i\in \mathcal{B}_{\rm good}$.\EEE}\label{fig: square}
\end{figure}

Due to \eqref{almostnoK} and \eqref{normal-angle-estimate0}, there is a \BBB universal \EEE constant $C>0$  such that  
    \begin{equation}\label{crucial_badgood}
        \begin{aligned}
            \mathcal{H}^{1}\big(  K(t)  \cap B_{\rm good}  \big) \leq C \BBB  \mathcal{H}^1(S(\psi)) \EEE \theta\,.
        \end{aligned}
    \end{equation}
Next, for technical reasons we straighten the jump set  $S(\psi)$ inside of $B_{\rm good}$, i.e., we construct a function $\phi$ which jumps only on straight segments \BBB inside $B_{\rm good}$. \EEE To this end, for $Q_i = Q_{r_i}(x_i) \in \mathcal{B}_{\rm good}$ we define
      \begin{equation}\label{Ri-estimates2-neu}
      H_i :=  \lbrace y \in Q_i \colon   (y-x_i) \cdot \nu_\psi(x_i)  =0 \rbrace, \quad  R_i :=    \lbrace y \in  Q_i \colon   |(y-x_i) \cdot \nu_\psi(x_i)| \le \theta r_i \rbrace, \quad     L_i := \partial Q_i \cap \partial R_i,
    \end{equation}
    see Figure~\ref{fig: square}. By an extension result we  can  modify $\psi$ in each square $Q_i \in \mathcal{B}_{\rm good}$ and obtain a function  $\phi \in W^{2,\infty}(\Omega)$ with $\phi = \psi$ on $\Omega \setminus B_{\rm good}$ and $S(\phi) \cap Q_i \subset L_i \cup H_i$ for each $Q_i \in \mathcal{B}_{\rm good}$. As  \BBB $\psi \in W^{2,\infty}(\Omega \setminus \overline{S(\psi)})$ \EEE and  $  \mathcal{L}^2(B_{\rm good}) \le C\theta$  by \eqref{besicovitch-props}, this can be done such that the estimate  
    \begin{align}\label{elastic for later}
\int_{\Omega} \Phi(\nabla \phi) \, {\rm d} x \le \int_{\Omega} \Phi(\nabla \psi) \, {\rm d} x + C\theta
\end{align} 
holds, for some $C$ depending on $\psi$. Moreover, recalling the definition $\varphi(\nu)=\frac{2\kappa}{\sqrt{3}}\sum_{\mathbf{v}\in\mathcal{V}}|\nu \cdot \mathbf{v}|$ for $\nu \in \mathbb{S}^1$, using the continuity of $\varphi$, and \eqref{normal-angle-estimate0}--\eqref{normal-angle-estimate} we find
\begin{align}\label{crack  for later}
\int_{S(\phi)} \varphi(\nu_\phi) \, {\rm d} \mathcal{H}^1 \le \int_{S(\psi)} \varphi(\nu_\psi) \, {\rm d} \mathcal{H}^1 + C\theta.
\end{align}

We now start to  construct the sequence of atomistic displacements $\psi_n \colon \mathcal{L}_n(\Omega )\to \R$ separately in the sets   $B_{\rm bad}$ and $\Omega _n \setminus B_{\rm bad}$. First, for   $x\in \mathcal{L}_n(\Omega )  \setminus B_{\rm bad}$ we will simply set  $\psi_n(x)=\phi(x)$. On $B_{\rm bad}$, however, the construction is more delicate as, due to the presence of the crack set $K(t)$, we cannot simply  discretize the function $\phi$. Instead, we will employ a \emph{transfer of jump sets} á la {\sc Francfort and Larsen} \cite{Francfort-Larsen:2003} to transfer the jump set of $\phi$ that lies in $K(t)$ onto $K_n^{{\rm L}}(t)$ defined in \eqref{def: sets-evol}. Afterwards, we discretize the resulting displacement.

%With some variant of the original construction, one can also get some additional properties: a  variant by {\sc Giacomini}, see \cite[(5.36)--(5.38)]{Giacomini:2005b}, shows that $\phi_n$ can be also constructed in the way that 
%\begin{equation}\label{limsup-transfer2}
%    \mathcal{H}^{1}\big((S(\phi_n)\setminus \lbrace x \in S(v_n) \colon [v_n](x) \ge \rho_n \rbrace     \cap    B_{\rm bad} \big)\leq C\theta\,\quad \, \text{ for all }  \, n\ge N.
%\end{equation}

\subsection*{Bad squares: jump transfer}

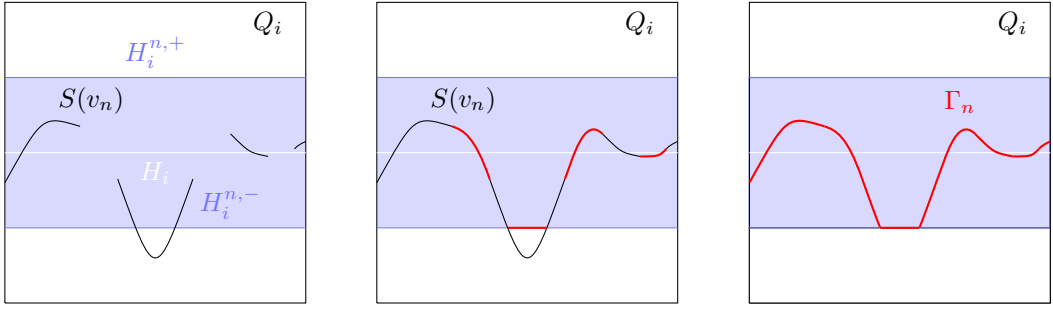
\begin{figure}
    \begin{minipage}{.3\textwidth}
    \begin{tikzpicture}[scale=0.5]
        \filldraw[fill=blue!15!white, draw=white]{} (0,6)-- (8,6) --(8,2)-- (0,2)--cycle;
        \node[anchor=north] at (7,8) {$Q_i$};
        \draw[color=black] (0,0) -- (8,0)  -- (8,8) --  (0,8) -- (0,0);
        \draw[color=blue!55!white] (0,2) -- (8,2);
        \draw[color=blue!55!white] (0,6) --node[anchor=south] {$H_i^{n,+}$}(8,6);
        \node[anchor=south,color=blue!55!white] at (6,2) {$H_i^{n,-}$};
        \draw[color=white] (0,4)--node[anchor=north] {$H_i$}
        (8,4);
        \draw (0,3.2) .. controls (1,5)  .. (2,4.7);
     %   \draw[color=red, dashed] (2,4.7) .. controls (2.3,4.6) and  (2.6,4.45)  .. (3,3.3);
        \draw (3,3.3) .. controls (4,0.5)  .. (5,3.3);
       % \draw[color=red, dashed] (5,3.3) .. controls (5.2,3.8) and (5.5,5) .. (6,4.5);
        \draw (6,4.5) .. controls (6.5,4) ..(7,3.9);
    %\draw[color=red, dashed] (7,3.9) .. controls (7.5,3.9) ..(7.7,4.1);
        \draw (7.7,4.1) .. controls (7.8,4.2) ..(8,4.3);
       % \draw[color=red, dashed] (4-0.52,2)--(4+0.52,2);
        \node[anchor=south] at (2.3,4.8) {$S(v_n)$}; 
       % \draw[color=green](0,2) -- node[anchor=east] {$L_i$}(0,6);
        %\draw[color=green](8,2) -- node[anchor=west] {$L_i$} (8,6);
    \end{tikzpicture}
\end{minipage}
\begin{minipage}{.3\textwidth}
     \begin{tikzpicture}[scale=0.5]
        \filldraw[fill=blue!15!white, draw=white]{} (0,6)-- (8,6) --(8,2)-- (0,2)--cycle;
        \node[anchor=north] at (7,8) {$Q_i$};
        \draw[color=black] (0,0) -- (8,0) --  (8,8) -- (0,8) -- (0,0);
        \draw[color=blue!55!white] (0,2) -- %node[anchor=south] {$H_i^{-}$}%
        (8,2);
        \draw[color=blue!55!white] (0,6) --%node[anchor=south] {$H_i^{-}$}%
        (8,6);
        \draw[color=white] (0,4)--%node[anchor=north] {$H_i$}%
        (8,4);
        \draw (0,3.2) .. controls (1,5)  .. (2,4.7);
        \draw[color=red, thick] (2,4.7) .. controls (2.3,4.6) and  (2.6,4.45)  .. (3,3.3);
        \draw (3,3.3) .. controls (4,0.5)  .. (5,3.3);
        \draw[color=red, thick] (5,3.3) .. controls (5.2,3.8) and (5.5,5) .. (6,4.5);
        \draw (6,4.5) .. controls (6.5,4) ..(7,3.9);
        \draw[color=red, thick] (7,3.9) .. controls (7.5,3.9) ..(7.7,4.1);
        \draw (7.7,4.1) .. controls (7.8,4.2) ..(8,4.3);
        \draw[color=red,thick] (4-0.52,2)--(4+0.52,2);
        \node[anchor=south] at (2.3,4.8) {$S(v_n)$}; 
        %\draw[color=green](0,2) -- node[anchor=east] {$L_i$}(0,6);
       % \draw[color=green](8,2) -- node[anchor=west] {$L_i$} (8,6);
        \end{tikzpicture}
    \end{minipage}
    \begin{minipage}{.3\textwidth}
      \begin{tikzpicture}[scale=0.5]
        \filldraw[fill=blue!15!white, draw=blue]{} (0,6)-- (8,6) --(8,2)-- (0,2)--cycle;
        \node[anchor=north] at (7,8) {$Q_i$};
        \draw[color=black] (0,0) -- (8,0) --  (8,8) -- (0,8) -- (0,0);
        \draw[color=white] (0,4)--(8,4);
        \draw[color=red, thick](0,3.2) .. controls (1,5)  .. (2,4.7);
        \draw[color=red, thick] (2,4.7) .. controls (2.3,4.6) and  (2.6,4.45)  .. (3,3.3);
        \draw[color=red, thick](3,3.3) .. controls (4,0.5)  .. (5,3.3);
        \draw[color=red, thick] (5,3.3) .. controls (5.2,3.8) and (5.5,5) .. (6,4.5);
        \draw[color=red, thick] (6,4.5) .. controls (6.5,4) ..(7,3.9);
        \draw[color=red, thick] (7,3.9) .. controls (7.5,3.9) ..(7.7,4.1);
        \draw[color=red, thick] (7.7,4.1) .. controls (7.8,4.2) ..(8,4.3);
        
        \node[anchor=south, color=red] at (5.6,4.8) {$\Gamma_n$}; 
       % \draw[color=green](0,2) -- node[anchor=east] {$L_i$}(0,6);
       % \draw[color=green](8,2) -- node[anchor=west] {$L_i$} (8,6);
        \filldraw[fill=white, draw=black]{} (0,2)-- (8,2) --(8,0)-- (0,0)--cycle;
        \draw[color=blue!55!white] (0,2) -- %node[anchor=south] {$H_i^{-}$}%
        (8,2);
        \draw[color=blue!55!white] (0,6) --%node[anchor=south] {$H_i^{-}$}%
        (8,6);
        \draw[color=red, thick] (4-0.52,2)--(4+0.52,2);
        \end{tikzpicture}
    \end{minipage}
    \caption{Construction of the continuous curve $ \Gamma^n_i \subset R_{i}^{n}$ from $S(v_n)$. The curve is related to a level set $v_n$ and cut at the segments $H_i^{n, \pm }$. Note that the figure is only a schematic illustration \BBB as \EEE the jump set $S(v_n)$ is actually contained in $K_n^{L}(t)$, which consists of straight segments \BBB of length $\frac{\eps_n}{2}$.\EEE } \label{Gamma-constr} 
\end{figure}

Let us now come to the essential points of the \emph{jump transfer}. For all details, however, we refer to \cite{Francfort-Larsen:2003}. Due to the construction of the collection $\mathcal{B}_{\rm bad}$, we can repeat the reasoning in \cite[Theorem 2.1]{Francfort-Larsen:2003}.  As before, let $(v_n)_{n}$ be the sequence with $S(v_n) \BBB \,  \tilde{\subset} \, \EEE  K_n^{L}(t)$ and  $v_n \to v$ in $L^1(\Omega )$, $\nabla v_n \rightharpoonup \nabla v$ weakly in $L^{\OOO 3/2 \EEE}(\Omega ;\R^2)$, where $S(v)  \BBB \, \tilde{=} \,  \EEE K(t)$. Then, it is shown that the jump set of $\phi$ can be transferred on to the ones of $(v_n)_n$: This means that there exist functions $\phi_n\in SBV^{\MMM 2 \EEE}(\Omega )$ with $\phi_n = \phi$ on $\Omega \setminus B_{\rm bad}$  such that 
 \begin{align}\label{a sup}
 \Vert  \phi_n\Vert_\infty  \le  C\Vert  \phi\Vert_\infty, \quad \quad     \Vert \nabla \phi_n\Vert_\infty \le C\Vert \nabla  \phi \Vert_\infty  \quad \text{ for all $n \in \N$}
 \end{align}
 for a universal $C>0$ and  such that there exists an index $N\in \N$ depending on $\theta$ with  
\begin{equation}\label{limsup-transfer1}
        \mathcal{H}^{1}\big((S(\phi_n)\setminus S(v_n)) \cap B_{\rm bad}\big) \leq C\theta\,\quad \, \text{ for all }  \, n\ge N.
\end{equation}
The main idea of the proof is to transfer the jump of $\phi$ in each $Q_i \in \mathcal{B}_{\rm bad} $ onto a continuous curve $\Gamma^n_i \subset Q_i$ related to the boundary of a certain  level set  for $v_n$,  see Figure \ref{Gamma-constr}  for the construction of $\Gamma^n_i$. \JJJ (The construction provides a function with discontinuities on the boundary of a set of finite perimeter $P^n_i \subset Q_i$. Without restriction, up to filling holes and removing components, $P^n_i$ can be chosen such that both $P^n_i$ and $Q_i \setminus P^n_i$ are connected sets (more precisely, indecomposable, see \cite[Sections 4,5]{Ambrosio-Morel}). Then, by \cite[Theorem 7]{Ambrosio-Morel} the boundary $\partial^* P^n_i$ can be represented by a Jordan curve and, as a consequence, $\Gamma^n_i := \partial^* P^n_i \cap Q_i$ is a continuous curve.) \EEE  Here, one uses \eqref{besicovitch-props2} and the $BV$ coarea formula to show that $\mathcal{H}^1(\Gamma^n_i \setminus S(v_n)) \le C\theta r_i$, see estimate \cite[(2.15)]{Francfort-Larsen:2003}.

More precisely, for each $n\in \N$ and each $Q_i:=Q_{r_i}(x_i)\in\mathcal{B}_{\rm bad}$ there are  two lines $H^{n,+}_i$ and $H^{n,-}_i$ with normal $\nu_\phi(x_i)$ which lie above and below the middle line $H_i$   containing the point $x_i$, also with normal $\nu_\phi(x_i)$, such that 
      \begin{equation}\label{Ri-estimates2}
R^{n}_i \supset Q_i \cap   \lbrace y \colon   |(y-x_i) \cdot \nu_\phi(x_i)| \le  2   \theta r_i \rbrace,  
    \end{equation}
    \begin{equation}\label{Ri-estimates}
%    \mathcal{H}^{1}\big(K(t)\setminus \bigcup\nolimits_{\mathcal{B}_{\rm bad}} R_i \big)\leq C\theta, \quad \quad 
      \mathcal{H}^{1}\big(\bigcup\nolimits_{\mathcal{B}_{\rm bad}}  L^n_i \big)\leq C\theta
    \end{equation}
for a universal constant $C>0$, where $R^{n}_i$ denotes  the rectangular subset that lies between $H^{n,+}_i$ and $H^{n,-}_i$, and $L^{n}_i:=\partial R^{n}_i \setminus (H^{n,+}_i \cup H^{n,-}_i)$ denotes its lateral boundaries,  cf.\ \cite[(2.10)--(2.11)]{Francfort-Larsen:2003}. We again refer to   Figure \ref{fig: square} for an illustration. 
%In \cite[(5.36)]{Giacomini:2005b}  it is observed that excluding jumps of small height $\rho_n$ with $\rho_n \to 0$ does essentially not affect the argument.   
  The function $\phi_n$ is then constructed such that $\phi_n=\phi$ outside of $\bigcup_{\mathcal{B}_{\rm bad }} R^{n}_i $ and \BBB inside $\bigcup_{\mathcal{B}_{\rm bad }} R^{n}_i $  it is defined by reflection. More precisely,  \EEE  since $S(\phi) \cap Q_i \subset  R^{n}_i$ by \eqref{normal-angle-estimate2} and \eqref{Ri-estimates2}, the construction in \cite{Francfort-Larsen:2003} actually yields  
\begin{align}\label{franfort again}
S(\phi_n) \cap Q_i \, \tilde{ \subset} \,  \Gamma^n_i \BBB \cup L_i^n \EEE  \quad \quad \text{for all $Q_i \in \mathcal{B}_{\rm bad}$}.
\end{align}

We point out that, \BBB as \EEE the construction of $\phi_n$ in \cite{Francfort-Larsen:2003} relies on a reflection, this leads only to a Lipschitz-regularity of $\phi_n$  in the two components of $Q_i \setminus \Gamma^n_i$,    see  \eqref{a sup}. However, for our purposes we will need that also $\nabla \phi_n$ is Lipschitz continuous. Therefore, we use a variant of the construction  which relies on a   cut-off argument and allows us to construct functions $\phi_n$   such that,  for a constant $C > 0$,  %depending on $\theta$%
we have 
\begin{align}\label{eq: Lipschitz1}
\nabla \phi_n \in SBV( \Omega ;\R^2) \  \text{ with }   \Vert  \nabla (\nabla \phi_n) \Vert_{L^\infty(\Omega )} \le   C \theta^{-1} (\min\nolimits_i r_i)^{-1}\Vert    \nabla \phi \Vert_{\infty} + \EEE  C\Vert    \nabla^2 \phi \Vert_{\infty}. 
\end{align}
\JJJ (Cut-off constructions for $SBV$-functions are by now standard, see e.g.\ \cite[Sections 2.3, 2.4]{BFLM}, and have also already been used in related second-order problems, see \cite[Lemma 4.4]{FHP} or \cite[Proposition 4.2]{donnarumma}.) \EEE More precisely, supposing $\nu_\phi   (x_i) = e_2$ and $x_i = 0$ for notational convenience, we let $T^\pm_i := \lbrace y \in  Q_i \colon   \pm y_2      \in  (\theta r_i,2\theta r_i)   \rbrace $ and     $b_\pm(y_1) := \fint_{\theta r_i}^{2\theta r_i} \phi(y_1,\BBB \pm \EEE s) \, {\rm d}s $ for $y_1 \in (-r_i,r_i)$. In the connected component $P^{n,+}_i$ of $Q_i \setminus \Gamma^n_i$ containing $H_i^{n,+}$, we set  $\phi_n(y_1,y_2) = \phi(y_1,y_2) \varphi^+_i(y_1,y_2) + (1-\varphi^+_i(y_1,y_2))b_+(y_1)$, where  $\varphi^+_i \in C^\infty(Q_i)$ with  $\varphi_i^+=1$ on $\lbrace y_2 >2\theta r_i \rbrace $ and $\varphi_i^+= 0$ on $\lbrace y_2 < \theta r_i \rbrace $.  In the other connected components of $Q_i \setminus \Gamma^n_i$ we perform a similar cutoff with $b_-(y_1)$. Omitting exact details, we mention that,  to see the bounds \eqref{a sup} and \eqref{eq: Lipschitz1}, we particularly exploit $\Vert \partial_1 b_{\pm}  \Vert_{L^\infty(-r_i,r_i) } \le C\Vert  \nabla \phi \Vert_\infty$ and $\Vert \phi(y_1,\BBB \pm \EEE \cdot) - b_{\pm}(y_1) \Vert_{L^\infty(\theta r_i,2\theta r_i  )} \le C\theta r_i \Vert \nabla \phi \Vert_\infty   $    for all $y_1 \in (-r_i,r_i)$, and $\Vert \nabla \varphi^\pm_i \Vert_\infty \le C(\theta r_i)^{-1}$  For later convenience, up to performing a suitable small modification of the functions $b_\pm$, by using \eqref{Ri-estimates} one can also achieve that, \BBB for a sufficiently small constant $c_\theta$ depending on $\theta$, \EEE    
\begin{align}\label{eq: Lipschitz2}
{\rm (i)} \ \ \bigcup\nolimits_{\mathcal{B}_{\rm bad}} L^n_i \cup \Gamma^n_i \, \tilde{ \subset} \, S(\phi_n), \quad \quad {\rm (ii)}  \ \  \mathcal{H}^1\big( \big\{ x \in S(\phi_n) \cap B_{\rm bad} \colon|[\phi_n](x)| \le \BBB c_\theta \EEE  \big\} \big)    \le C\theta.
\end{align}
In the following, we suppose that the additional properties \eqref{franfort again}--\eqref{eq: Lipschitz2} are satisfied. As by definition of the sequence $(v_n)_n$ we have $S(v_n) \BBB \,  \tilde{\subset} \, \EEE K_n^{\rm L}(t)$, \eqref{limsup-transfer1} yields
    \begin{equation}\label{limsup-transfer3}
        \mathcal{H}^{1}\big((S(\phi_n)\setminus  K_n^{{\rm L}}(t)) \cap B_{\rm bad} \big) \leq C\theta\,\quad \, \text{ for all }  \, n\ge N.
    \end{equation}   
Moreover, $S(\phi_n)$ is essentially closed, i.e., $\mathcal{H}^1(\overline{S(\phi_n)} \setminus S(\phi_n)) = 0$, and by  \eqref{franfort again} and  \eqref{eq: Lipschitz2} we have
\begin{align}\label{jumpiiii}
{S(\phi_n)}  \, \tilde{=} \,   \big(  S(\phi) \setminus \overline{B_{\rm bad}} \big) \cup       \bigcup_{Q_i \in \mathcal{B}_{\rm bad}} ( L^n_i  \cup  \Gamma^n_i).  
\end{align}
With these preparations, we can define the sequence of atomistic displacements  by   
\begin{align}\label{psin2}
\psi_{n}(x):= \phi_n(x)\quad \text{ for all }  x\in \mathcal{L}_n(\Omega ) \,.
\end{align}
This means in particular that $\psi_n(x) = \phi(x)$ for $x \in \mathcal{L}_n(\Omega ) \setminus {B}_{\rm bad}$ and $\psi_n(x) = \psi(x)$ for $x \in \mathcal{L}_n(\Omega ) \setminus (B_{\rm good} \cup {B}_{\rm bad})$ by the definition of $\phi$ below \eqref{Ri-estimates2-neu}.   As $B_{\rm good} \cup B_{\rm bad} \subset U$, this gives $\psi_n(x) =  \BBB g(t,x) \EEE $ for all $x \in \L_n(\Omega \setminus U)$, i.e., $\psi_n \in   \mathcal{A}_{n}\BBB (g(t)) \EEE$.

\subsection*{Some auxiliary lemmas}
In the following,  we choose $k \in  \BBB \N\EEE$ such that   $t \in [t^{k}_{n},t_n^{k+1})$. (Clearly, $k$ depends on $t$ and $n$ but we do not include this in the notation for simplicity.) 
We  recall the definition of the energies  $\mathcal{E}_n^{k, {\rm{ela}}}$, $\mathcal{E}_n^{k, {\rm{cra}}}$, $ \mathcal{E}_n^{k, {\rm{rem}}}$, and $\mathcal{E}_n^{k, {\rm{bdy}}}$   in  \eqref{energy:main}, see also \eqref{energy along evo2}. We also recall the definition of   $M_{n,\triangle}^{k, \mathbf{v}}(v_n)$,   $\mathcal{C}_n^{k}(v_n)$, and $\mathcal{C}_n^{k, {\rm L}}(v_n)$ in  \eqref{eq: noacons-neu}, \eqref{nom jump tri}, and \eqref{big jump tri}. We denote by  $\mathcal{D}^{\mathbf{v}}_{n}(\phi_n)$ the collection of triangles such that   $\overline{S(\phi_n)}$  intersects the side of the triangle that is oriented in direction $\mathbf{v}\in \mathcal{V}$. We also define 
% we use the notation  
\begin{align}\label{DDD}
\mathcal{D}_{n}(\phi_n):= \bigcup_{j=1,2,3}  \mathcal{D}^{\mathbf{v}_j}_{n}(\phi_n)= \big\{\triangle \colon  \overline{S(\phi_n)}\cap\partial \triangle \neq \emptyset \big\} \,,
\end{align}
  which is the collection of triangles that are intersected by $\overline{S(\phi_n)}$. \EEE Note that the sets $\BBB\mathcal{D}^{\mathbf{v}_j}_{n}\EEE (\phi_n)$ are obviously not disjoint in general. We now formulate three auxiliary results that will help us to prove Theorem \ref{stability}. 

\begin{lemma}[Small stretching of springs if there is no  jump]\label{stablemma2}
For each $\mathbf{v} \in \mathcal{V}$ it holds
 \begin{align}\label{auch das nochc}
\,|({\psi}_n)_{\triangle}  \cdot \mathbf{v}| \le C \quad \text{ for all }   \triangle \in  \mathcal{T}_n \setminus \mathcal{D}^{\mathbf{v}}_n(\phi_n) 
\end{align}
 for a constant $C>0$ depending on  $\phi$. Moreover, for $\eps_n$ small enough, we have  
 \begin{align}\label{auch das nochb}
M_{n,\triangle}^{k, \mathbf{v}}(\psi_n) \le   M_{n,\triangle}^{k, \mathbf{v}} \vee R \le   M_{n,\triangle}^{k, \mathbf{v}}(u^k_n) \vee R  \ \ \  \text{ for all }   \triangle \in  \mathcal{T}_n \setminus \mathcal{D}^{\mathbf{v}}_n(\phi_n).
 \end{align}
\end{lemma}
\noindent
Whenever a spring does not intersect the jump set of $\phi_n$, the stretching is necessarily controlled, see \eqref{auch das nochc}, and thus does effectively not increase the memory variable in the sense of \eqref{auch das nochb}.

\begin{lemma}[Small stretching of springs despite of jump]\label{stablemma3}
For $\eps_n$ small enough, it holds that 
\begin{align}\label{auch das noch3}
\varepsilon_n \sum_{\mathbf{v} \in \mathcal{V}}  \# \, \big\{   \triangle \in \mathcal{D}^{\mathbf{v}}_{n}(\phi_n)  \colon \,   M_{n,\triangle}^{k, \mathbf{v}}(\psi_n)  \leq R     \big\} \le C\theta
\end{align}
and
\begin{align}\label{auch das noch2}
\varepsilon_n  \# \, \big(\mathcal{D}_n(\phi_n)  \setminus \Cn{k}{\psi_n}  \big)  \le C\theta. 
\end{align}
\end{lemma}

In principle, it can happen that a spring intersects the jump set of $S(\phi_n)$, but still the stretching is smaller than $R$. This is the case if the jump height $[\phi_n](x)$ for some $x \in S(\phi_n)$ is small. Due to the fine control on the amount of points where this can happen (use $\phi_n = \phi$ outside of $B_{\rm bad}$ and  \eqref{eq: Lipschitz2}(ii)), the energy contribution of such springs can be controlled in terms of $\theta$.

For the third lemma, we define the collection of triangles inside   $B_{\rm bad}$ by $\mathcal{T}_n^{\rm bad}  =  \lbrace \triangle \in \mathcal{T}_n \colon \triangle \subset   B_{\rm bad} \rbrace\,$.  
\begin{lemma}[Bounds in $B_{\rm bad}$ and complement]\label{stablemma}
It holds that 
    \begin{align}\label{crucial-limsup}
    \limsup_{n \to \infty}  \frac{\varepsilon_n}{2}  \sum_{ ( \mathcal{C}_n^{k}(\psi_n) \cap  \mathcal{D}_n(\phi_n)) \setminus  \mathcal{T}_n^{\rm bad}   }  \ \ \biggl( \sum_{M_{n,\triangle}^{k, \mathbf{v}}(\psi_n)>R}     \Psi(M_{n,\triangle}^{k, \mathbf{v}}(\psi_n))-  & \sum_{   \Met{k}{u^k_n}  > R}  \Psi(M_{n,\triangle}^{k, \mathbf{v}}(u^{k}_n))     \biggr)  \notag \\ & \leq \int_{S(\phi)\setminus K(t) } \varphi(\nu_\phi)\,d \mathcal{H}^{1} + C\theta
    \end{align}
     and 
\begin{equation}\label{2equals2}
\varepsilon_n\, \#    \Big\{   (\triangle,\mathbf{v}) \in \big(  \mathcal{C}_n^{k}(\psi_n) \cap   \mathcal{T}_n^{\rm bad}  \cap \mathcal{D}_n(\phi_n) \big) \times \mathcal{V}  \colon M_{n,\triangle}^{k, \mathbf{v}}(\psi_n)> M_{n,\triangle}^{k, \mathbf{v}}(u^k_n)  \vee R \text{ and } M_{n,\triangle}^{k, \mathbf{v}}(u^k_n)\le R_n \Big\}   \le C\theta.
\end{equation}  
In particular, we have
\begin{align}\label{oje}
    \varepsilon_n \, \#   \mathcal{C}_{n}^{k}(\psi_n)  + \eps_n \# \mathcal{C}_{n}^{k}(u^k_n)  \leq   C\,.
\end{align}     
\end{lemma}
Here, we recall the shorthand notation  in \eqref{crucial-limsup}, namely that the first sum  runs over triangles $\triangle$ and the other sums run over directions $\mathbf{v}$ satisfying the corresponding inequality. The two estimates \eqref{crucial-limsup} and \eqref{2equals2} deal with 'broken springs' that are intersected by $S(\phi_n)$ in $\Omega \setminus B_{\rm bad}$ and $B_{\rm bad}$, respectively. The proof of the first estimate \eqref{crucial-limsup} resembles the construction of a recovery sequence in the atomistic-to-continuum $\Gamma$-convergence result \cite{FS152}. In particular, outside $B_{\rm bad}$ no transfer of jump enters the argument. Yet, the latter is at the core of proving \eqref{2equals2} which indeed fundamentally relies on \eqref{limsup-transfer1}.  Eventually, estimate \eqref{oje} follows from   \eqref{auch das nochc}--\eqref{2equals2}.\EEE
\subsection*{Proof of the stability result.} 
We defer the proofs of the three lemmas to the end of the section and proceed with the proof of Theorem \ref{stability}.

\begin{proof}[Proof of Theorem \ref{stability}]
     Let $t\in [0,T]$ be given and, for each $n \in \N$, choose $k$ such that $t \in [t^k_n,t_n^{k+1})$. Let $(\psi_n)_n$ be the sequence defined in \eqref{psin2}.  By construction, we clearly have \BBB $\psi_n(x)   = g(t,x)$ \EEE  for $x \in \mathcal{L}_n (\Omega  \setminus U)$ as $\phi_n$ coincides with $\psi$ outside of $B_{\rm good} \cup B_{\rm bad}$ and $B_{\rm good} \cup B_{\rm bad} \subset U$. Recalling the discussion at the beginning of the section, it suffices to show \eqref{stab0neu}--\eqref{stab3neu} for an arbitrary but fixed $\theta>0$.

 \emph{Step 1: Proof of \eqref{stab0neu}}.
Consider the interpolation  $\tilde{\psi}_{n}$ related to $\psi_n$ which is affine on each triangle of $\mathcal{T}_n$. In view of \eqref{psin2} and the regularity of $\psi$, one can check that 
$$\chi_{\Omega _n}  \tilde{\psi}_{n}  \to \psi \quad \text{ on $L^1(\Omega  \setminus \big(B_{\rm good} \cup B_{\rm bad})\big)$}. $$ 
By  \eqref{a sup},  \eqref{psin2}, and the first property in  \eqref{besicovitch-props},  we also find  
$$\Vert \chi_{\Omega _n}  \tilde{\psi}_{n}  - \psi \Vert_{L^1(B_{\rm good} \cup B_{\rm bad})}    \le C  \big(  \Vert \phi \Vert_{L^{\infty}(\Omega  )} +   \Vert \psi \Vert_{L^{\infty}(\Omega  )} \big)  \theta^2 \leq C \theta^2. $$
This shows \eqref{stab0neu}.

\emph{Step 2: Proof of \eqref{stab1neu}}. 
Recall the definition of $\mathcal{E}^{k,{\rm cra}}_{n}(\psi_{n})$ and $\mathcal{E}_n^{k, {\rm{rem}}}(\psi_{n})$ from \eqref{eq: different enegies}, namely  
% In view of \eqref{eq: different enegies}, we recall that 
\begin{equation*}\label{helpestimate}
    \mathcal{E}^{k,{\rm cra}}_{n}(\psi_{n})  =  \frac{\varepsilon_n}{2}\sum_{\Cl{k}{\psi_n}} \sum_{\Met{k}{\BBB \psi_n}>R_n} \Psi(\Met{k}{\psi_n})               
\end{equation*}    
  and  
\begin{align*}
\mathcal{E}_n^{k, {\rm{rem}}}(\psi_{n})  &   =\frac{\varepsilon_n}{2}\sum_{\Cn{k}{\psi_n}}\Big( \sum_{R< \Met{k}{\psi_n}\le R_n} \Psi(\Met{k}{\psi_n}) + \sum_{\Met{k}{\psi_n}  \leq R} \Psi(\OOO \varepsilon_n^{1/2} \EEE \,|({\psi}_n)_{\triangle}  \cdot \mathbf{v}|) \Big)   \\
     & \ \ \  + \frac{\varepsilon_n}{2}\sum_{\Cn{k}{\psi_n} \setminus \Cl{k}{\psi_n}}   \sum_{\Met{k}{\psi_n}> R_n} \Psi(\Met{k}{\psi_n})\,.
\end{align*}    
% \RRR hier summieren wir eigentlich gerade alle beitrage in gebrochenen zellen auf. Wenn man das bemerkt, kann man sich das aufschreiben oben eigentlich sparen... \EEE 
The same definition holds for $u_n^k$ in place of $\psi_n$. We then estimate the contributions by 
 $$   \mathcal{E}^{k,{\rm cra}}_{n}(\psi_{n})-\mathcal{E}^{k,{\rm cra}}_{n}(u^{k}_{n}) +  \mathcal{E}_n^{k, {\rm{rem}}}(\psi_{n}) - \mathcal{E}_n^{k, {\rm{rem}}}(u^{k}_{n}) \le  \mathcal{F}_n^{k,1} + \mathcal{F}_n^{k,2} +\mathcal{F}_n^{k,3}  +\mathcal{F}_n^{k,4} , $$
 where 
\begin{align*}
\mathcal{F}_n^{k,1} &:=  \frac{\varepsilon_n}{2}\sum_{ ( \mathcal{C}_n^{k}(\psi_n) \cap  \mathcal{D}_n(\phi_n)) \setminus  \mathcal{T}_n^{\rm bad}   }  \ \ \biggl( \sum_{M_{n,\triangle}^{k, \mathbf{v}}(\psi_n)>R}     \Psi(M_{n,\triangle}^{k, \mathbf{v}}(\psi_n))-  \sum_{   \Met{k}{u^k_n}  > R}  \Psi(M_{n,\triangle}^{k, \mathbf{v}}(u^{k}_n))     \biggr), \\
  \mathcal{F}_n^{k,2} &:= \frac{\varepsilon_n}{2}\sum_{ \mathcal{C}_n^{k}(\psi_n)  \cap \mathcal{D}_n(\phi_n) \cap   \mathcal{T}_n^{\rm bad}  } \biggl( \sum_{M_{n,\triangle}^{k, \mathbf{v}}(\psi_n)>R}     \Psi(M_{n,\triangle}^{k, \mathbf{v}}(\psi_n))-  \sum_{\Met{k}{u^{k}_n}  > R  }  \Psi(M_{n,\triangle}^{k, \mathbf{v}}(u^{k}_n))     \biggr),\\
    \mathcal{F}_n^{k,3} &:= \frac{\varepsilon_n}{2}\sum_{\mathcal{C}_n^{k}(\psi_n) \setminus \mathcal{D}_n(\phi_n)} \biggl( \sum_{M_{n,\triangle}^{k, \mathbf{v}}(\psi_n)>R}     \Psi(M_{n,\triangle}^{k, \mathbf{v}}(\psi_n))-  \sum_{\Met{k}{u^{k}_n}  > R \,  }  \Psi(M_{n,\triangle}^{k, \mathbf{v}}(u^{k}_n))     \biggr),\\  
    \mathcal{F}_n^{k,4} &:=      \frac{\varepsilon_n}{2}\sum_{\Cn{k}{\psi_n}}    \sum_{\Met{k}{\psi_n}  \leq R} \Psi(\varepsilon_n^{\MMM 1/2 \EEE}\,|({\psi}_n)_{\triangle}  \cdot \mathbf{v}|)\,.
   \end{align*}
    Here, we summed up all the contributions in 'broken' triangles $\triangle\in \Cn{k}{\psi_n}$. However, we dropped the contributions of  related to $u^k_n$  in triangles $\Cn{k}{u^k_n} \setminus \Cn{k}{\psi_n}$ and the contributions $ \Psi(\varepsilon_n^{\MMM 1/2 \EEE}\,|({u}^{k}_n)_{\triangle}  \cdot \mathbf{v}|) $ for $(\triangle, \mathbf{v}) \in \Cn{k}{u^k_n} \times \mathcal{V}$ with $\Met{k}{u^k_n}  \leq R$, which are all clearly nonnegative.\\
We now estimate the various terms. First, \eqref{crucial-limsup}   together with  \eqref{crack  for later}   implies
\begin{align}\label{to conclude1}
\limsup_{n\to \infty} \mathcal{F}_n^{k,1}  \leq \int_{S(\psi)\setminus K(t) } \varphi(\nu_\psi)\,d \mathcal{H}^{1} + C\theta\,.   
\end{align}  
        To estimate $  \mathcal{F}_n^{k,2} $ we use \eqref{2equals2}. More precisely, we denote the set on the left-hand side of \eqref{2equals2} by $\Lambda_n$  for shorthand.  Then, recalling that $\Psi$ is increasing  and $\Psi \le \kappa$ by \ref{phiprop3} and \ref{phiprop4}, we compute
\begin{align*}
 \mathcal{F}_n^{k,2} \le \frac{3\eps_n}{2} \# \big( (\mathcal{C}_n^{k}(\psi_n)  \cap \mathcal{D}_n(\phi_n) \cap   \mathcal{T}_n^{\rm bad}) \setminus  \Lambda_n \big) \big( \kappa - \Psi(R_n) \big) +  \frac{3\eps_n}{2}  \# \Lambda_n  \kappa.         
\end{align*}
Then, by $\lim_{t \to \infty} \Psi(t) = \kappa$, \eqref{2equals2}, and the fact that $\# \Cn{k}{\psi_n} \le C\eps_n^{-1}$ (see \eqref{oje}) we conclude        
\begin{align}\label{to conclude2}
\limsup_{n\to \infty} \mathcal{F}_n^{k,2} \le C \theta\,.
      \end{align}       
    By \eqref{auch das nochb} \BBB and the monotonicity of $\Psi$ \EEE we directly obtain       $  \mathcal{F}_n^{k,3}  \le 0$.             Eventually, by \eqref{auch das nochc} and \eqref{auch das noch3}  we get 
\begin{align*}\limsup_{n\to \infty}\mathcal{F}_n^{k,4}  &  \le   \limsup_{n\to \infty}  \frac{\varepsilon_n}{2}  \sum_{  \mathbf{v}  \in \mathcal{V}}\# \big(\Cn{k}{\psi_n} \setminus \mathcal{D}^{\mathbf{v}}_n(\phi_n) \big)   \Psi\big(\varepsilon_n^{\MMM 1/2 \EEE}C\big)  \\ & \ \ \ +   \limsup_{n\to \infty}  \frac{\varepsilon_n}{2} \sum_{  \mathbf{v}  \in \mathcal{V}}   \# \,  \big\{ \triangle \in \mathcal{D}^{\mathbf{v}}_n(\phi_n)   \colon \,    \Met{k}{\psi_n}  \leq R     \big\} \Psi(R) \le C\theta,
\end{align*}
where we also used  $\# \Cn{k}{\psi_n} \le C\eps_n^{-1}$ by \eqref{oje} and  $\lim_{t \to 0} \Psi(t) = 0$. This along with \eqref{to conclude1}--\eqref{to conclude2} concludes the proof of \eqref{stab1neu}.

  \emph{Step 3: Proof of \eqref{stab2neu}}.   The essential point is to prove that there exists a sequence of sets $(\Upsilon_n)_n \subset \Omega $ such that 
            \begin{align}\label{strongconv}
 & \chi_{\Omega _n}  \nabla \BBB \hat{\psi}^k_n  \EEE  - \nabla \phi_n \to 0 \text{ pointwise a.e.\ in }  \Omega , \quad \mathcal{L}^2(\Upsilon_n) \le C\theta\eps_n, \quad   \Vert         \nabla \BBB \hat{\psi}^k_n  \EEE   \Vert_{L^\infty(\Omega_{n} \setminus \Upsilon_n)} \le  C.
            \end{align}
Assume for the moment that \eqref{strongconv} is true and recall the definition of $\Psi_n^{\rm cell}$ in \eqref{psitriangle} with $ \Psi_n^{\rm cell}  \le C\eps_n^{-1}$ on $\R^2$. Thus, using \eqref{eq: ela en} and \eqref{strongconv} it follows 
\begin{align}\label{last estimate}
    \limsup_{n\to \infty}\,\mathcal{E}_{n}^{k,{\rm ela}}(\psi_{n}) & = \limsup_{n\to \infty} \int_{\Omega _{n}} \Psi_{n}^{\rm cell}(\nabla \BBB \hat{\psi}^k_n  \EEE(x)) \,{\rm d}x  \le \limsup_{n\to \infty} \int_{\Omega _{n} \setminus \Upsilon_n} \Psi_{n}^{\rm cell}(\nabla \BBB \hat{\psi}^k_n  \EEE(x)) \,{\rm d}x  +   C\mathcal{L}^2(\Upsilon_n) \eps_n^{-1}  \notag  \\ 
    & \le  \limsup_{n\to \infty} \int_{\Omega _n\setminus (B_{\rm bad} \cup \Upsilon_n )   } \Psi_{n}^{\rm cell}(\nabla \BBB \hat{\psi}^k_n  \EEE(x)) \,{\rm d}x + \limsup_{n \to \infty}  \int_{B_{\rm bad}  \setminus \Upsilon_n  } \Psi_{n}^{\rm cell}(\nabla \BBB \hat{\psi}^k_n  \EEE(x)) \,{\rm d}x + C\theta \,.
\end{align}
Note that the function $x\mapsto \Psi_n^{\rm cell}(\nabla \BBB \hat{\psi}^k_n  \EEE(x) )$ is uniformly bounded on $\Omega_n \setminus \Upsilon_n$ because of \eqref{strongconv} and the fact that $\Psi_n^{\rm cell}$ is bounded on $\lbrace \BBB  |y| \EEE \le C\rbrace$, see Lemma \ref{lem:help2}(iii). Hence, by  $\mathcal{L}^2(B_{\rm bad }) \le \theta^2 $   (see the first item in  \eqref{besicovitch-props}), the second integral in \eqref{last estimate} can be controlled in terms of $C\theta$. For the first integral we argue similarly to the proof of Lemma \ref{elastic-lowersim}. The fact that $\phi_n = \phi$ in $\Omega  \setminus B_{\rm bad}$, the reverse Fatou lemma, and the first item of \eqref{strongconv} lead to   
\begin{align*}
    \limsup_{n\to \infty} \int_{ \Omega_n \setminus (B_{\rm bad} \cup \Upsilon_n  )}  \Psi_{n}^{\rm cell}(\nabla \BBB \hat{\psi}^k_n  \EEE(x)) \,{\rm d}x  & \leq  \int_{\Omega \setminus B_{\rm bad}}  \hspace{-0.25cm} \Phi(\nabla \phi) \, {\rm d}x +   \limsup_{n\to \infty} \int_{\BBB \Omega_n\EEE \setminus (B_{\rm bad} \cup \Upsilon_n)}   |\nabla\BBB \hat{\psi}^k_n  \EEE|^{\MMM 2 \EEE} \frac{\omega(\varepsilon_n^{\MMM 1/2 \EEE}\nabla\BBB \hat{\psi}^k_n  \EEE)}{(\varepsilon_n^{\MMM 1/2 \EEE}|\nabla\BBB \hat{\psi}^k_n  \EEE|)^{\MMM 2 \EEE}}   \,  {\rm d}x  \\ & \le \int_{\Omega \setminus B_{\rm bad}}   \Phi(\nabla \phi) \, {\rm d}x\,.
\end{align*}
Here we have used Lemma \ref{lem:help2}(i) and the fact $\|\nabla \BBB \hat{\psi}^k_n  \EEE\|_{L^{\MMM 2 \EEE}(\Omega_n\setminus\Upsilon_n)}\leq C$ as well as $\varepsilon_n^{\MMM 1/2 \EEE}|\nabla\BBB \hat{\psi}^k_n  \EEE|^{\MMM 2 \EEE}\to 0$ uniformly, which both follow from the last item in \eqref{strongconv}.   In view of \eqref{elastic for later}, this concludes the proof of \eqref{stab2neu}.

Now, we proceed with the proof of \eqref{strongconv}. First, since $\nabla \BBB \hat{\psi}^k_n  \EEE$ is zero  on each $\triangle \in \mathcal{C}_n^k(\psi_n)$ by construction, and $\nabla \phi_n$ is uniformly bounded on $\Omega $, see \eqref{a sup}, we get 
   \begin{equation}\label{strongconv1}
\sum_{\triangle \in  { \mathcal{C}_n^k(\psi_n)}} \chi_{ \triangle }(\nabla \BBB \hat{\psi}^k_n  \EEE- \nabla \phi_n) \to 0 \quad \text{strongly in}\, L^1(\Omega ;\R^2)\,, 
   \end{equation}           
   where we used $\mathcal{L}^2(\triangle) = \frac{\sqrt{3}}{4} \eps_n^2$ and $\# \mathcal{C}_n^k(\psi_n) \le C \eps_n^{-1}$ by \eqref{oje}.  Now, we consider triangles  $\triangle \in \mathcal{D}_{n}(\phi_n) \setminus  \mathcal{C}_n^k(\psi_n)$. Here, we have no good control on  $\nabla \BBB \hat{\psi}^k_n  \EEE= \nabla \tilde{\psi}_{n}$, but we can control the volume of the set in terms of $\theta\varepsilon_n$. Letting   $\Upsilon_n = \bigcup_{\triangle \in  \mathcal{D}_{n}(\phi_n) \setminus   \mathcal{C}_n^k(\psi_n)} \triangle$, by    \eqref{auch das noch2} we compute  
   \begin{equation}\label{dont you remove}
   \mathcal{L}^2(\Upsilon_n) \le \,  C \eps_n^2 \# (\mathcal{D}_{n}(\phi_n) \setminus   \mathcal{C}_n^k(\psi_n)) \le C\theta \eps_n.  
   \end{equation}   
Eventually, we consider triangles $\triangle\in\mathcal{T}_n \setminus (\mathcal{D}_{n}(\phi_n) \cup \Cn{k}{\psi_n})$. In view of \eqref{jumpiiii}, this implies  that either (a) $\triangle \not\subset B_{\rm bad}$ with $\triangle \cap \bigcup_{\mathcal{B}_{\rm bad}}  L_i^n  = \emptyset$ or (b) $\triangle  \subset B_{\rm bad}$. First, assume that  $\triangle \not\subset B_{\rm bad}$ and $\triangle \cap \bigcup_{\mathcal{B}_{\rm bad}}  L^n_i  = \emptyset$. Then, by construction of $\phi_n$ below \eqref{eq: Lipschitz1}, we have $\psi_n(x)= \phi_n (x) = \phi(x)$ for all vertices $x$ of $\triangle$. Since $\phi \in W^{2,\infty}(\Omega  \setminus \overline{S(\phi)})$, $\nabla \phi$ is Lipschitz on $\triangle$ and thus $\sup_{x,x' \in \triangle} |\nabla \phi(x) - \nabla \phi(x')| \le \Vert \nabla^2 \phi \Vert_{L^\infty(\triangle)} \eps_n$. By the mean value theorem this yields
$$\Vert \nabla \tilde{\psi}_n  - \nabla \phi \Vert_{L^\infty(\triangle)} \le C\eps_n \Vert \nabla^2 \phi \Vert_{L^\infty(\Omega  )}.      $$  
In a similar fashion,  if  $\triangle  \subset B_{\rm bad}$, we have $\psi_n(x)=\phi_n(x)$ for the vertices of $\triangle$, and as $\nabla \phi_n$ is Lipschitz on $\triangle$, we get
$$\Vert \nabla \tilde{\psi}_n  - \nabla \phi_n \Vert_{L^\infty(\triangle)} \le   C    \eps_n \Vert \nabla^2 \phi_n \Vert_{L^\infty( \triangle \EEE)} \le C\eps_n \big(  \theta^{-1} (\min\nolimits_i r_i)^{-1}\Vert    \nabla \phi \Vert_\infty +   \Vert    \nabla^2 \phi \Vert_{\infty}  \big),       $$   
where we also used \eqref{eq: Lipschitz1}. As $\triangle \notin \Cn{k}{\psi_n}$, we have $\BBB \hat{\psi}^k_n  \EEE= \tilde{\psi}_n$ on such triangles. Therefore, for fixed $\theta>0$, we obtain \EEE
\begin{equation}\label{strongconv3}
     \sum_{\triangle \in \mathcal{T}_n \setminus ( \mathcal{D}_{n}(\phi_n) \cup\Cn{k}{\psi_n} ) }\chi_{ \triangle }(\nabla \BBB \hat{\psi}^k_n  \EEE- \nabla \phi_n)  \to 0 \quad \text{strongly in}\, L^1(\Omega ;\R^2)\,.   
\end{equation} 
Combining \eqref{strongconv1}--\eqref{strongconv3} with the fact that $\Vert \nabla \phi_n  \Vert_\infty \le  C  $ yields \eqref{strongconv}.

\emph{Step 4: Proof of \eqref{stab3neu}}. Recall the definition of $\mathcal{E}_n^{\rm{bdy}}$ in \eqref{bdy energy}. Let $\mathcal{T}_n^{\rm bdy} = \lbrace \triangle \in \mathcal{T}_n \colon \, \partial \triangle \cap \partial \Omega _n   \neq \emptyset\rbrace$ and $\Omega_n^{\rm bdy} = \bigcup_{\triangle \in \mathcal{T}_n^{\rm bdy}} \triangle $. First, since $ B_{\rm good} \cup B_{\rm bad} \subset \subset \Omega $ by construction, we get $\triangle \cap (B_{\rm good} \cup B_{\rm bad}) = \emptyset$ for all $\triangle \in \mathcal{T}_n^{\rm bdy}$, for $\eps_n$ small enough. By  $\mathcal{H}^1(S(\psi) \cap \Omega_n^{\rm bdy}) \to 0$ as $n \to \infty$ and  the regularity of $S(\psi)$,   this shows that 
\begin{align}\label{last estimate2}
\eps_n \#  (\mathcal{T}_n^{\rm bdy} \cap \mathcal{D}_n(\phi_n))= \eps_n \# (\mathcal{T}_n^{\rm bdy} \cap \mathcal{D}_n(\psi))  \to 0,
\end{align}
where $\mathcal{D}_n(\psi)$ is defined analogous to \eqref{DDD}, i.e., as the collection of triangles that are intersected by $\overline{S(\psi)}$. \EEE Now, the contribution of triangles in $\mathcal{T}^{\rm bdy}_n \cap  (\mathcal{C}_n^k(\psi_n) \setminus  \mathcal{D}_n(  \phi_n) \EEE )$ is nonpositive by \eqref{auch das nochb} and repeating the argument for $\mathcal{F}^{k,3}_n$ above.  Furthermore, the contribution of triangles in   $\mathcal{T}^{\rm bdy}_n \cap    \mathcal{D}_n(\BBB \phi_n) \EEE $ is negligible by \eqref{last estimate2} and the boundedness of $\Psi$. Concerning the triangles $\triangle \notin \mathcal{C}^k_n(\psi_n)$, we can follow the lines of Step~3, in particular \eqref{last estimate}, where we use that also the integral over $ \Omega_n^{\rm bdy}  \setminus  \Upsilon_n$ vanishes since $\mathcal{L}^2(\Omega_n^{\rm bdy}) \to 0$ and we have $\mathcal{L}^2(\Upsilon_n \cap \Omega_n^{\rm bdy}) \eps_n^{-1} \to 0 $ by   \eqref{last estimate2}. This shows \eqref{stab3neu}.
\end{proof}  

\subsection*{Proofs of the auxiliary lemmas}
We proceed to prove the three auxiliary results.

\begin{proof}[Proof of Lemma \ref{stablemma2}]
 Fix $\mathbf{v} \in \mathcal{V}$ and $\triangle \in  \mathcal{T}_n \setminus \mathcal{D}^\mathbf{v}_n(\phi_n)$. Let $x,x'$ be the pair of vertices with $x-x'$ parallel to $\mathbf{v}$. Since $\phi_n\in W^{1,\infty}(\Omega  \setminus \overline{S(\phi_n)})$, more precisely  $L:=  \sup_n \Vert \nabla \phi_n \Vert_\infty<\infty$  by \eqref{a sup}, and the segment between $x$ and $x'$ does not intersect $\overline{S(\phi_n)}$, we get  
    \[|\phi_n(x)-\phi_n(x')|\leq \|\nabla \phi_n\|_\infty \eps_n \le L\eps_n\,.\] 
As $\psi_n = \phi_n$ on $\L_n(\Omega )$,  we get  
$$|({\psi}_n)_{\triangle} \cdot \mathbf{v}|\leq L\,, $$  
which shows   \eqref{auch das nochc}. In particular, for  $\varepsilon_{n}$ small enough, we have $\varepsilon_{n}^{\MMM 1/2 \EEE} L \le R$, and thus by \eqref{eq: noacons}--\eqref{eq: noacons-neu}, we find   
$$\Met{k}{\psi_n}  \le \Me{k}  \vee R \le \Met{k}{u_n^k} \vee R.$$  
This concludes the proof of \eqref{auch das nochb}.
\end{proof}
 
For convenience we proceed with the proof of Lemma \ref{stablemma} and address Lemma  \ref{stablemma3} at the end.
\begin{proof}[Proof of Lemma  \ref{stablemma}]
We first prove \eqref{crucial-limsup} in Steps 1--3. Step 4  addresses \eqref{2equals2}, and finally in Step 5 we show \eqref{oje}. As a preparation for the proof of \eqref{crucial-limsup}, we split the energy contribution  on the left hand side   into three terms,  namely
\begin{align}\label{splittttt}
 \frac{\varepsilon_n}{2}\sum_{ ( \mathcal{C}_n^{k}(\psi_n) \cap  \mathcal{D}_n(\phi_n)) \setminus  \mathcal{T}_n^{\rm bad}   }  \hspace{-0.1cm} \biggl( \sum_{M_{n,\triangle}^{k, \mathbf{v}}(\psi_n)>R}     \Psi(M_{n,\triangle}^{k, \mathbf{v}}(\psi_n))-  \sum_{   \Met{k}{u^k_n}  > R}  \Psi(M_{n,\triangle}^{k, \mathbf{v}}(u^{k}_n))     \biggr)  =  \mathcal{F}^{\rm good}_{n}+  \mathcal{F}^{\rm rest}_{n} +\mathcal{F}^{\rm glue}_{n}\,,
\end{align}
where in each of the three terms the summation is restricted to corresponding subsets of $\mathcal{T}_n$, namely
\begin{align*}
\mathcal{T}_n^{\rm good} = \lbrace \triangle \colon \triangle \cap B_{\rm good} \neq \emptyset  \rbrace, \quad   \mathcal{T}_n^{\rm rest} =  \lbrace \triangle \colon \triangle \cap (B_{\rm good} \cup B_{\rm bad }) = \emptyset   \rbrace,   \quad 
\mathcal{T}_n^{\rm glue}  =  \lbrace \triangle \colon \triangle \cap \partial B_{\rm bad} \neq \emptyset  \rbrace.
\end{align*}
We now treat the three terms separately.  

\emph{Step 1: $\mathcal{F}^{\rm good}_{n}$}.  First, by  $\Psi \le \kappa$ we get
\begin{align*}
   \mathcal{F}^{\rm good}_{n}  \le \sum_{\mathbf{v} \in \mathcal{V}} \frac{\kappa \eps_n}{2}  \#  \big(   \mathcal{D}^{\mathbf{v}}_n(\phi_n)\cap \mathcal{T}_n^{\rm good} \big).
\end{align*}
 Recall that $\psi_n = \phi$ on $\L_n(\Omega )  \setminus B_{\rm bad}$ and $\triangle \cap B_{\rm bad} = \emptyset$ for all $\triangle \in \mathcal{T}_n^{\rm good}$ (for $\eps_n$ small enough). In particular, we get $\mathcal{D}_{n}(\phi) \cap \mathcal{T}_n^{\rm good} = \mathcal{D}_{n}(\phi_n) \cap \mathcal{T}_n^{\rm good}$, where  we use the notation introduced in \eqref{DDD} also for $\phi$. Thus, we have  
    \begin{equation}\label{estimate-goodcubes}
\mathcal{F}_n^{\rm good}  \leq \sum_{j=1}^3 \frac{\kappa  \eps_n}{2}   \# \big( \mathcal{D}^{\mathbf{v}_j}_{n}(\phi) \cap \mathcal{T}_n^{\rm good} \big).
    \end{equation}
    Denote by $\nu_{\phi}$ again the measure-theoretic unit normal to $S(\phi)$. For a square $Q_i=Q_{r_i}(x_i)\in \mathcal{B}_{\rm good}$ we now prove that
    %\RRR changed from here, also tried to improve the slicing text : 
    \EEE
    %that is centered at $x\in S(\psi)\setminus K(t)$ 
    
    \begin{equation}\label{jumpsetsize-estimate}
    \#\big\{\triangle\in \mathcal{D}^{\mathbf{v}_j}_{n}(\phi)\,\text{with}\; \triangle\cap Q_i \neq \emptyset \big\} \leq 2r_i \frac{4|\nu_{\phi}(x_i)\cdot \mathbf{v}_j|}{\sqrt{3}\varepsilon_n} + \frac{C r_i  }{\varepsilon_n}\,\theta \,. \end{equation}  
By \BBB \eqref{normal-angle-estimate2} and \EEE the construction below \eqref{Ri-estimates2-neu} we observe that, \BBB for each $\triangle$ with $\triangle \cap Q_i \neq \emptyset$, we have   that $S(\phi) \cap \triangle$ \EEE is contained in the two segments forming   $ L_i$  and in the segment $H_i$. As $\mathcal{H}^1(L_i) = 4\theta r_i$, the number of triangles intersecting   $L_i$  can clearly be controlled by $C\theta r_i/\eps_n$. Therefore, it remains to consider triangles intersecting $H_i$. Recall that $H_i$ is a segment with normal vector $\nu_\phi(x_i)$. Let $\mathbf{v}_j^\bot$ be a unit vector orthogonal to $\mathbf{v}_j$, and denote by $\Pi_j (H_i)$ the orthogonal projection of $H_i$ onto $\R \mathbf{v}_j^\bot$. We compute $\mathcal{H}^1(\Pi_j (H_i)) = \mathcal{H}^1(H_i)|\nu_{\phi}(x_i)\cdot \mathbf{v}_j| = 2r_i |\nu_{\phi}(x_i)\cdot \mathbf{v}_j|$, and then by basic geometric considerations, one can estimate the amount of springs in direction $\mathbf{v}_j$ that are crossed by $H_i$ as $ (\frac{\sqrt{3}\eps_n}{2} )^{-1} 2r_i |\nu_{\phi}(x_i)\cdot \mathbf{v}_j| + 1$. As each spring is contained in two triangles, this implies 
$$ { \#\big\{\triangle\in \mathcal{D}^{\mathbf{v}_j}_{n}(\phi)\,\text{with}\; \triangle  \cap  H_i \neq \emptyset \big\} \le   2r_i   \frac{4|\nu_{\phi}(x_i)\cdot \mathbf{v}_j|}{\sqrt{3}\varepsilon_n} + 2.}$$ 
For $\eps_n$ small enough with respect to $\theta>0$   this implies \eqref{jumpsetsize-estimate}.  
Then, using  $\mathcal{H}^1(Q_i \cap S(\phi)) \ge 2r_i $, recalling the definition $\varphi(\nu)=\frac{2\kappa}{\sqrt{3}}\sum_{\mathbf{v}\in\mathcal{V}}|\nu \cdot \mathbf{v}|$ for $\nu \in \mathbb{S}^1$, and putting together \eqref{estimate-goodcubes}--\eqref{jumpsetsize-estimate}, we obtain
\begin{equation*}
   \mathcal{F}_n^{\rm good} \leq \kappa \sum_{Q_i \in \mathcal{B}_{\rm good}} \sum_{j=1}^3 \mathcal{H}^{1}(S(\phi)\cap Q_i)  \Big(\frac{2|\nu_{\phi}(x_i)\cdot \mathbf{v}_j|}{\sqrt{3} } + C\,\theta\Big) \le \sum_{Q_i \in \mathcal{B}_{\rm good}}   \mathcal{H}^{1}(S(\phi)\cap Q_i)\,\varphi(\nu_{\phi}(x_i))+ C\, \theta 
\end{equation*}
   for a constant $C>0$ also depending on $\mathcal{H}^1(S(\phi))$. Then, due to \eqref{crucial_badgood}, we also get
$$
  \mathcal{F}_n^{\rm good} \leq  \sum_{Q_i \in \mathcal{B}_{\rm good}}   \mathcal{H}^{1}\big( (S(\phi) \setminus K(t) )\cap Q_i\big)\,\varphi(\nu_{\phi}(x_i))+ C\, \theta.
$$  
Eventually, using  that $\nu_{\phi}(x_i)$ is the normal vector to $H_i$ and   $\mathcal{H}^1(L_i) = 4\theta r_i$  for each $Q_i \in \mathcal{B}_{\rm good}$, we find
 \begin{equation}\label{Bgood-finalestimate}  
  \mathcal{F}_n^{\rm good} \leq \int_{ (S(\phi) \setminus K(t) ) \cap B_{\rm good}} \varphi(\nu_\phi) \, {\rm d} \mathcal{H}^1 + C\theta.   
   \end{equation}

 \emph{Step 2:  $\mathcal{F}^{\rm rest}_{n}$}. We recall  that   $\psi_n = \psi$ on $\L_n(\Omega ) \setminus  (B_{\rm good}\cup B_{\rm bad}) $.  Repeating the argument leading to \eqref{estimate-goodcubes} we find 
    \begin{equation}\label{estimate-rest}
       \mathcal{F}^{\rm rest}_{n}\leq    \frac{3\kappa \, \varepsilon_n}{2}   \; 
        \# \big( \mathcal{D}_{n}(\psi) \cap \mathcal{T}_n^{\rm rest} \big).
    \end{equation} 
    In view of the second item of \eqref{besicovitch-props2} and the regularity of $S(\psi)$ we have  $\# ( \mathcal{D}_{n}(\psi) \cap \mathcal{T}_n^{\rm rest} ) \le C\theta \eps_n^{-1}$.  Therefore, we have
        \begin{equation}\label{Brest-finalestimae}
 \mathcal{F}^{\rm rest}_{n} \leq C\theta\,.
 \end{equation}

 \emph{Step 3:  $\mathcal{F}^{\rm glue}_{n}$}.
    Let $Q_i\in \mathcal{B}_{\rm bad}$ and consider a triangle   $\triangle \in \mathcal{T}^{\rm glue}_n \cap \mathcal{D}_n(\phi_n) \cap \mathcal{C}^k_n(\psi_n)$ with $\triangle \cap \partial Q_i \neq \emptyset$. Recalling the sets  $L^n_i$  introduced in \eqref{Ri-estimates},  by \eqref{normal-angle-estimate2} and \eqref{jumpiiii} we get  $\triangle \cap \overline{S(\phi_n)} = \emptyset$  whenever $\triangle \cap L^n_i  = \emptyset$.  It is thus left to consider  $\triangle$  with $\triangle \cap L^n_i\neq \emptyset$.  For this, we use \eqref{Ri-estimates}:  Because   $\bigcup_{Q_i \in \mathcal{B}_{\rm bad}}  L^n_i   $ consists of straight segments, we obtain
     \begin{equation*}
 \sum_{ Q_i \in \mathcal{B}_{\rm bad}}    \# \big\{\triangle\in \mathcal{T}_n \colon \triangle\cap  L^n_i  \neq \emptyset\big\}\leq \frac{C}{\varepsilon_n} \sum\nolimits_{Q_i \in \mathcal{B}_{\rm bad}} \mathcal{H}^{1}( L^n_i  ) \leq \frac{C\,\theta}{\varepsilon_n}\,.
    \end{equation*} 
   Since an estimate analogous to  \eqref{estimate-goodcubes}  or \eqref{estimate-rest} also holds for $\mathcal{F}^{\rm glue}_{n}$, this finally leads to
    \begin{equation}\label{theta-glue}
 \mathcal{F}^{\rm glue}_{n} \le C\theta\,.
 \end{equation}
 Recalling the splitting \eqref{splittttt} and  combining \eqref{Bgood-finalestimate}, \eqref{Brest-finalestimae},  and \eqref{theta-glue}  we conclude the proof of \eqref{crucial-limsup}.

 \emph{Step 4: Proof of  \eqref{2equals2}}. We define 
 $ {\mathcal{D}}^*_{n}:= \{\triangle\in \mathcal{T}^{\rm bad}_n\colon (S(\phi_n)\setminus K^{{\rm L}}_n(t)) \cap \triangle \neq \emptyset \}$
 and first argue that it suffices to show that 
\begin{align}\label{eq: to checkii}
\#  {\mathcal{D}}^*_{n} \le C\theta \eps_n^{-1}.
\end{align} 
To see this, we show that $\mathcal{T}_n^{\rm bad }\setminus \mathcal{D}_n^*$ does not intersect with the collection from \eqref{2equals2}. For each $\triangle \in \mathcal{T}_n^{\rm bad }\setminus \mathcal{D}_n^*$ and some $\mathbf{v} \in \mathcal{V}$ we either have $M_{n,\triangle}^{k, \mathbf{v}}(u^k_n)> R_n $ or $M_{n,\triangle}^{k, \mathbf{v}}(u^k_n) \le R_n$. In the first case, the conditions in \eqref{2equals2} are obviously violated. In the latter case, the segment of $\partial \triangle$ in direction $\mathbf{v}$ does not intersect $K^{{\rm L}}_n(t)$, cf.\ \eqref{eq: hat sets} \BBB and \eqref{def: sets-evol},  \EEE and therefore, by $\triangle \cap (S(\phi_n)\setminus K_n^{L}(t))=\emptyset$, we conclude $\triangle \notin \mathcal{D}^{\mathbf{v}}_n(\phi_n)$. Then, \eqref{auch das nochb} implies $M_{n,\triangle}^{k, \mathbf{v}}(\psi_n) \le   M_{n,\triangle}^{k, \mathbf{v}}(u^k_n) \vee R $, and thus the triangle does not lie in the collection given in \eqref{2equals2}. This shows that \eqref{eq: to checkii} induces the desired estimate.
 
Let us now show \eqref{eq: to checkii}. To this end, we claim that by the construction of $\psi_n$ we can find a \BBB universal \EEE constant $C>0$ such that for every $Q_i\in \mathcal{B}_{\rm bad}$  we have
    \begin{equation}\label{estimate_badcubes}
    \#\big\{\triangle \in   \mathcal{D}^*_{n}  \colon \triangle \subset Q_i \big\}\leq \frac{C}{\varepsilon_n}\, \mathcal{H}^{1}\big((S(\phi_n)\setminus K^{{\rm L}}_n(t))\cap Q_i\big).
    \end{equation}
Recall the curve $\Gamma^n_i$ satisfying \eqref{franfort again}. Note also that $K^{ {\rm L}}_n(t)$ consists of segments with length $\frac{\varepsilon_n}{2}$ and each two of such segments either touch or have distance at least $\frac{\eps_n}{2}$. Therefore, a connected component  $\Gamma_{i,l}^n$ of $(\Gamma^n_i\setminus K^{{\rm L}}_n(t))\cap Q_i$  satisfies  either   $\mathcal{H}^{1}(\Gamma_{i,l}^n)\geq \frac{\varepsilon_n}{2}$ or that the start and endpoint of $ \Gamma^n_{i,l}$ intersect the same or two consecutive segments contained in $K^{ {\rm L}}_n(t)$. In the latter case, in the construction of $\phi_n$, the part $\Gamma^n_{i,l}$ of the curve $\Gamma^n_i$ could be replaced by a piecewise affine curve inside $K^{ {\rm L}}_n(t)$ such that the estimate \eqref{limsup-transfer1} still holds. Therefore, in the following, we will suppose without restriction that all connected components $(\Gamma_{i,l}^n)_l$ of  $(\Gamma^n_i\setminus K^{{\rm L}}_n(t))\cap Q_i$  fulfill $\mathcal{H}^{1}(\Gamma_{i,l}^n)\geq \frac{\varepsilon_n}{2}$. We now aim to prove   for each of such components  that 
    \begin{equation}\label{estimate-component}
        \#\big\{\triangle\colon \, \Gamma^n_{i,l}\cap \triangle\neq \emptyset  \big\} \leq \frac{C}{\varepsilon_n}\mathcal{H}^{1}( \Gamma^n_{i,l})\,.
    \end{equation}  
From this, we then deduce \eqref{estimate_badcubes} as 
 \[ \#\big\{\triangle \in   \mathcal{D}^*_{n}  \colon \triangle \subset Q_i \big\}\leq \sum\nolimits_l \#\big\{\triangle\colon \Gamma^n_{i,l} \cap \triangle \neq \emptyset\big\} \leq \frac{C}{\varepsilon_n}\sum\nolimits_{l}\mathcal{H}^{1}(\Gamma^n_{i,l}) \leq \frac{C}{\varepsilon_n}\mathcal{H}^{1}\big((S(\phi_n)\setminus K^{{\rm L}}_n(t))\cap Q_i\big)\,.\] 
    \EEE In order to prove \eqref{estimate-component}, we look at the $  \varepsilon_n $-neighborhood of $ \Gamma^n_{i,l}$, i.e.,
 \begin{align}\label{consider neih}
 N_{\varepsilon_n}( \Gamma^n_{i,l}):=\{x \in \R^2\colon {\rm dist}(x, \Gamma^n_{i,l})\leq \varepsilon_n\}\,.
 \end{align} 
As $\mathcal{H}^1 ( \Gamma^n_{i,l})\geq \frac{\varepsilon_n}{2}$, an elementary geometric argument shows that there exists a universal constant $C>0$ such that 
\[ \mathcal{L}^2\big(N_{\varepsilon_n}( \Gamma^n_{i,l})\big)\leq C \varepsilon_n \mathcal{H}^{1}( \Gamma^n_{i,l})\,. \] 
Since $\mathcal{L}^2(\triangle) = \frac{\sqrt{3}}{4} \varepsilon_n^2$, we hence obtain  
    \[\#\big\{\triangle\colon \, \Gamma^n_{i,l}\cap \triangle\neq \emptyset  \big\}\le \# \big\{\triangle\colon \, \triangle\subset N_{\varepsilon_n}(\Gamma^n_{i,l})\big\} \leq \frac{\mathcal{L}^2\big(N_{\varepsilon_n}( \Gamma^n_{i,l})\big)}{\mathcal{L}^2(\triangle)}\leq \frac{C}{\varepsilon_n}\mathcal{H}^{1}( \Gamma^n_{i,l}) \,.\] 
This yields \eqref{estimate-component} and therefore \eqref{estimate_badcubes}. In view of  \eqref{limsup-transfer3} and \eqref{estimate_badcubes}, we eventually get \eqref{eq: to checkii} by summing over the finite number of squares in $\mathcal{B}_{\rm bad}$.

\emph{Step 5: Proof of  \eqref{oje}}. By Lemma \ref{lemma: dm} and Corollary \ref{cor: energy bound}, together with Proposition \ref{bounded-jump}, we obtain  $ \varepsilon_n \,\#  \mathcal{C}_n^{k}(u_n^{k})\leq C$. It remains to estimate $\#\mathcal{C}_{n}^{k}(\psi_n)$. Because of Lemma \ref{stablemma2}, we \BBB know  \EEE that $\mathcal{C}_n^{k}(\psi_n)\setminus \mathcal{D}^{\mathbf{v}}_n(\phi_n)\subset \mathcal{C}_n^{k}(u^{k}_n)$  for all $\mathbf{v}\in \mathcal{V}$   and thus $\#(\mathcal{C}_n^{k}(\psi_n)\setminus \mathcal{D}_n(\phi_n))\leq \#\mathcal{C}_n^{k}(u_n^{k})\leq C\varepsilon_n^{-1}$.   It hence remains to estimate $\#(\mathcal{C}_n^{k}(\psi_n)\cap \mathcal{D}_n(\phi_n))$. Due to \eqref{crucial-limsup}, we obtain   
\begin{equation*}
 \varepsilon_n \# \big( \mathcal{C}_n^{k}(\psi_n) \cap  \mathcal{D}_n(\phi_n) \big) \setminus  \mathcal{T}_n^{\rm bad}\leq C + \kappa \varepsilon_n \#\mathcal{C}_n^{k}(u_n^{k}) + C \theta  \leq C  \,. 
 \end{equation*}
Here, we have used that contributions of triangles with $\Met{k}{u^k_n}  > R$ for some $\mathbf{v}\in \mathcal{V}$ can be controlled by $\kappa \varepsilon_n \#\mathcal{C}_n^{k}(u_n^{k})$ and the fact that $\int_{S(\phi)\setminus K(t) } \varphi(\nu_\phi)\,d \mathcal{H}^{1} \le C$.     
Finally, if we denote the set on the left-hand side of \eqref{2equals2} by $\Lambda_n$, we obtain  $\big(\mathcal{C}_n^{k}(\psi_n) \cap  \mathcal{D}_n(\phi_n) \cap \mathcal{T}_n^{\rm bad}\big)\setminus  \Lambda_n \subset  \mathcal{C}_n^{k}(  u^{k}_n ) $. We hence can deduce from \eqref{2equals2} that 
 \begin{equation*}
    \varepsilon_n \#\big(\mathcal{C}_n^{k}(\psi_n) \cap  \mathcal{D}_n(\phi_n) \cap \mathcal{T}_n^{\rm bad}\big) \leq  \varepsilon_n  \#\mathcal{C}_n^{k}(  u^{k}_n  )   + C\theta \leq C\,.
\end{equation*}
Putting all bounds together we conclude the proof of \eqref{oje}. 
\EEE
\end{proof}

\begin{proof}[Proof of Lemma  \ref{stablemma3}]
 We start with  \eqref{auch das noch3}. To this end, recall the construction of $\phi_n$, see  particularly \eqref{jumpiiii}, and note that it coincides with $\phi$ outside of $B_{\rm bad}$. Therefore, for $c_{\theta}$ small enough, we find that $\mathcal{H}^1\big(\lbrace x \in S(\phi)\setminus B_{\rm bad} \colon    |[\phi](x)| \le c_\theta \rbrace \big) \le \theta $. Then, using \eqref{eq: Lipschitz2}(ii), possibly passing to a smaller $c_\theta$ depending only on $\theta$, we get 
\begin{align}\label{YYY}
\mathcal{H}^1\big(\lbrace x \in S(\phi_n) \colon    |[\phi_n](x)| \le c_\theta \rbrace \big) \le C\theta \quad \text{ for all $n \in \N$}, 
\end{align}  
where $C>0$ is a universal constant.  

As we already found the control  $\# ( \mathcal{D}_{n}(\psi) \cap \mathcal{T}_n^{\rm rest} ) \le C\theta \eps_n^{-1}$ preceding \eqref{Brest-finalestimae},  to show \eqref{auch das noch3} it thus suffices to consider triangles intersecting $B_{\rm good} \cup B_{\rm bad}$. We recall that $S(\phi_n) \BBB \cap (B_{\rm good} \cup B_{\rm bad}) $ is contained  in a finite number of curves whose number is bounded uniformly in $n$.  Therefore, there is only a bounded number of triangles intersecting more than one of these curves. Consequently, in order to show \eqref{auch das noch3} for $\eps_n$ small enough, it suffices to consider triangles $\triangle$ such that $N_{\varepsilon_n}(\triangle)$ intersects exactly one of these curves, where similarly to \eqref{consider neih} we define 
\[N_{\varepsilon_n}(\triangle):=\{x \in \R^2\colon {\rm dist}(x,\triangle)\leq \varepsilon_n\}\,.\]
Further, in view of \eqref{eq: to checkii},  the triangles that intersect any component $(\Gamma_{i,l}^n)$ of $\Gamma^n_i\setminus K_n^{k,L}$ are negligible. Hence, in order to show \eqref{auch das noch3}, it suffices to consider triangles $\triangle$ such that $N_{\varepsilon_n}(\triangle)$ does not intersect any component $( \Gamma^n_{i,l})_{i,l}$. \EEE Now, except for the parts $( \Gamma^n_{i,l})_{i,l}$ considered in \eqref{estimate-component}, the curves are either the straight segments $H_i=S(\phi)\cap Q_i$ and $L_i$ for squares $ Q_i\in \mathcal{B}_{\rm good}\EEE$, the segments $L^n_i$ related to $Q_i \in \mathcal{B}_{\rm bad}$, or a union of intervals contained in $K^{k,{\rm L}}_n$, see \eqref{jumpiiii}. We denote the collection of triangles  that intersect one of these curves   by $\mathcal{G}_{n}$ and define
$$ \mathcal{G}^{\mathbf{v}}_{n}:= \big\{   \triangle \in \mathcal{D}^{\mathbf{v}}_{n}(\phi_n) \cap \mathcal{G}_{n}  \colon \,   M_{n,\triangle}^{k, \mathbf{v}}(\psi_n)  \leq R   \big\}.$$
In view of  this preliminary discussion,  it suffices to check $\eps_n\# \mathcal{G}^{\mathbf{v}}_{n} \le  C\theta$ for all
 $\mathbf{v}\in \mathcal{V}$. Let us consider a triangle $\triangle \in \mathcal{G}^{\mathbf{v}}_{n}$. We note that $S(\phi_n) \cap N_{\varepsilon_n}(\triangle)$ is contained in a straight segment or in a union of intervals of length $\frac{\eps_n}{2}$ oriented parallel to the lattice directions.  Let $x$ and  $x'$ be the vertices of $\triangle$ such that $x-x'$ is in direction $\mathbf{v}$ and let $x_0$ be the unique point in \BBB    $\overline{S(\phi_n)}$ on the segment connecting $x$ and $x'$. \EEE By  the fact that  $ M_{n,\triangle}^{k, \mathbf{v}}(\psi_n)  \leq R$, \eqref{psin2}, and the fundamental theorem of calculus we find
     \[    R \ge   \frac{\eps_n^{\MMM 1/2 \EEE}|\phi_n(x)-\phi_n(x')|}{\eps_n} \ge  \eps_n^{\MMM -1/2 \EEE} \big(  |[\phi_n](x_0)|  -       \|\nabla \phi_n\|_\infty \eps_n \big) \ge   \eps_n^{\MMM -1/2 \EEE} |[\phi_n](x_0)| -   L   \eps^{1/2}_n\,,\]
     where   we set $L :=\sup_n \Vert \nabla \phi_n \Vert_\infty < \infty$, see \eqref{a sup}.   
This shows  $|[\phi_n](x_0)| \le   \eps_n^{\MMM 1/2 \EEE} R + L\eps_n.$ Since $S(\phi_n) \cap N_{\varepsilon_n}(\triangle)$ is contained in only one of the curves of $S(\phi_n)$, again using $ \Vert \nabla \phi_n \Vert_\infty \le L$, we find
$$ |[\phi_n](z)| \le  \eps_n^{\MMM 1/2 \EEE} R + CL \eps_n    \quad \text{ for all $z \in S(\phi_n) \cap N_{\varepsilon_n}(\triangle)$.} $$
\BBB (In fact, $x_0$ and each $z$ can be connected by two different curves on different components of  $N_{\varepsilon_n}(\triangle)) \setminus \overline{S(\phi_n)}$ of length $\sim \eps_n$.) \EEE 
 As $\mathcal{H}^1(S(\phi_n) \cap N_{\varepsilon_n}(\triangle))$ is at least $\eps_n/2$, for $\eps_n$ small enough depending on $\theta$  we get that 
$$\mathcal{H}^1\big(   \lbrace x \in  S(\phi_n) \cap N_{\varepsilon_n}(\triangle)\colon   |[\phi_n](x)| \le c_\theta  \rbrace \big) \ge \frac{\eps_n}{2}. $$ 
Summing  over all $\triangle$ in $\mathcal{G}^{\mathbf{v}}_{n} $ and observing that each neighborhood $N_{\varepsilon_n}(\triangle)$ intersects only a bounded number of neighborhoods $N_{\varepsilon_n}(\triangle')$,  $\triangle'\neq \triangle$, we get
$${\eps_n \# \mathcal{G}^{\mathbf{v}}_{n} \le C\mathcal{H}^1\big(   \lbrace x \in  S(\phi_n) \colon     |[\phi_n](x)| \le c_\theta  \rbrace \big)   \le C\theta,}$$
where in the last step we used \eqref{YYY}. This concludes the proof of \eqref{auch das noch3}. 

Eventually, \eqref{auch das noch2} is a simple consequence of \eqref{auch das noch3} and the definition of $\mathcal{D}_n(\phi_n)$ and $\mathcal{C}^k_n(\psi_n)$. 
\end{proof}

We close the section with the proof of    Corollary \ref{cor: stability}.

\begin{proof}[Proof of Corollary \ref{cor: stability}]
Inspection of the proof of Theorem \ref{stability} shows that it is irrelevant that $u_n(t)$ satisfies the boundary values $g(t_n^{k(t)})$. In fact, it only matters that $\psi = g(t)$ on $\Omega \setminus \overline{U}$. Therefore, we can consider another sequence of solutions $u^*_n(t)$ related to boundary values $g^*(t) \defas 0$ for all $t \in [0,T]$. We obtain $u^{*}_n(t) = 0$ for all $t\in [0,T]$ and the corresponding limiting evolution $u^*(t) = 0$ and $K^*(t) = \emptyset$ for all $t \in [0,T]$. Applying Theorem \ref{stability} in this situation for some  $\psi\in SBV^{\MMM 2 \EEE}(\Omega)$ with $\psi = g(t)$ on $\Omega \setminus \overline{U}$  yields $(\psi_{n})_{n}$   with $\psi_n \in   \mathcal{A}_{n} (g(t_n^{k(t)})) $  such that $\psi_n$ AC-converges to $\psi$ and 
$$        \limsup_{n\to \infty}\, \big(  \mathcal{E}^{ {\rm ela}}_{n}(\psi_n;t) +  \mathcal{E}^{ {\rm cra}}_{n}(\psi_n;t) +  \mathcal{E}^{ {\rm rem}}_{n}(\psi_n;t) + \mathcal{E}^{\rm{bdy}}_{n}(\psi_n; t ) \big)    \le     \int_{\Omega } \Phi(\nabla \psi) \, {\rm d}x + \int_{S({\psi})  } \varphi(\nu_\psi)\,d \mathcal{H}^{1}, $$
where we used $\mathcal{E}^{ {\rm cra}}_{n}(u^*_n(t);t) + \mathcal{E}^{ {\rm rem}}_{n}(u^*_n(t);t) = 0$ and $K^*(t) = \emptyset$.  Now it suffices to observe that $ \mathcal{E}_{n}(\psi_n;t) =  \mathcal{E}_n^{ {\rm{ela}}}(\psi_n;t) + \mathcal{E}_n^{ {\rm{cra}}}(\psi_n;t) + \mathcal{E}_n^{ {\rm{rem}}}(\psi_n;t) +\mathcal{E}_n^{ {\rm{bdy}}}(\psi_n;t)$ by \eqref{energy:main} and the fact that $ \mathscr{E}_{n}(\psi_n) =  \mathcal{E}_{n}(\psi_n;t)$ since the memory variable is not active due to $ u^*_n(t)  = 0$. 
\end{proof}

 \section*{Acknowledgements} 
 This research was funded by the Deutsche Forschungsgemeinschaft (DFG, German Research Foundation) - 377472739/GRK 2423/2-2023. The authors are very grateful for this support.

%  \RRR war mir nicht so klar: \EEE
 
%  Corollary \ref{cor: stability} follows immediatly from the preceeding arguments. 
% \begin{proof}[Proof of Corollary \ref{cor: stability}]
%     We choose $\psi_n(x) := \psi(x)$ for all $\mathcal{L}_{n}(\Omega )$. Note that we have \[\mathscr{E}_{n}(\psi_n)\mathscr{E}_{n}(\psi_n)=\mathcal{E}_n^{\rm ela}(\psi_n,0)+\mathcal{E}_n^{\rm cra}(\psi_n,0)\,.\]  Since $\mathcal{E}_n^{\rm rem}(u_n(0),0)=\mathcal{E}_n^{\rm cra}(u_n(0),0)=0$,  \RRR Warum? \EEE the assertion follows directly from \eqref{stab1neu} and \eqref{stab2neu}.   \RRR Warum sollte $K(0)=0$ sein? \EEE \end{proof}

% \RRR add a short proof for the case of recovery sequences, cf.\ . \EEE 
 
%\RRR Brezis 3.32 
%[1] Brezis, H., & Brézis, H. (2011). Functional analysis, Sobolev spaces and partial differential equations (Vol. 2, No. 3, p. 5). New York: Springer.

\end{document}